\documentclass[11pt]{book}
\usepackage{epsfig}

\newenvironment{proof}{{\it Proof:\/}}{$\Box$\vskip 0.08in}
\usepackage{amssymb}

\newtheorem{theorem}{Theorem}[section]
\newtheorem{lemma}[theorem]{Lemma}
\newtheorem{corollary}[theorem]{Corollary}
\newtheorem{proposition}[theorem]{Proposition}

\newtheorem{problem}[theorem]{Problem}
\newtheorem{remark}[theorem]{Remark}

\newtheorem{definition}[theorem]{Definition}
\newtheorem{example}[theorem]{Example}
\newtheorem{exercise}[theorem]{Exercise}
\newtheorem{property}[theorem]{Property}

\newcommand{\Z}{\mathbb Z}
\newcommand{\Q}{\mathbb Q}

\begin{document}
\renewcommand{\thechapter}{\Roman{chapter}}
\thispagestyle{empty}

\
\vspace{0.5in}
 \begin{center}
 {\LARGE\bf KNOTS}\\
{\bf From combinatorics of knot diagrams to combinatorial topology
based on knots}
\end{center}

\vspace*{0.3in}

\centerline{Warszawa, November 30, 1984 -- Bethesda, October 31, 2004}
\vspace*{0.3in}

 \begin{center}
                      J\'ozef H.~Przytycki
\end{center}

\vspace*{0.5in}
{\Large\bf Introduction}\\
\ \\
This book is
about classical Knot Theory, that is, about
the position of a circle (a knot) or of a number of disjoint circles
(a link) in the space $R^3$ or in the sphere $S^3$.
We also venture into Knot Theory in general 3-dimensional
manifolds.

The book has its predecessor in Lecture Notes on Knot Theory,
 published in Polish\footnote{The
Polish edition was prepared for the ``Knot Theory" mini-semester
at the Stefan Banach Center, Warsaw, Poland, July-August, 1995.}
in 1995  \cite{P-18}.
A rough translation of the Notes (by J.Wi\'sniewski) was
ready by the summer of 1995. It differed from the Polish edition
with the addition of
the full proof of Reidemeister's theorem. While I couldn't find
time to refine the translation and prepare the final manuscript,
I was adding new material and rewriting existing
chapters. In this way I created a new book based on the Polish
Lecture Notes
but expanded 3-fold.
Only the first part of Chapter III (formerly Chapter II),
on Conway's algebras is essentially unchanged from the Polish book
and is based on preprints \cite{P-1}.

As to the origin of the Lecture Notes, I was teaching an advanced course
in theory of 3-manifolds and Knot Theory at Warsaw University and it
was only natural to write down my talks (see Introduction to  (Polish)
Lecture Notes).
I wrote the proposal for the Lecture Notes by the  December 1, 1984
deadline.
In fact I had to stop for a while our work on generalization of
the Alexander-Conway and Jones polynomial in order to submit
the proposal. From that time several excellent books on Knot Theory
have been published on various level and for various readership.
This is reflected in my choice of material for the book -- knot theory
is too broad to cover every aspect in one volume. I decided to
concentrate on topics on which I was/am doing an active research.
Even with this choice the full account of skein module theory is
relegated to a separated book (but broad outline is given in
Chapter IX).

In the first Chapter we offer historical perspective to
the mathematical theory of knots, starting from the first
precise approach to Knot Theory
by Max Dehn
and Poul Heegaard in the Mathematical Encyclopedia [D-H,1907].
We start the chapter by introducing
lattice knots and polygonal knots.
The main part of the chapter is devoted to the proof of
Reidemeister's theorem which allows combinatorial treatment of
Knot Theory.

In the second Chapter we offer the history of Knot
Theory starting from the ancient Greek tract on surgeon's slings,
through Heegaard's thesis
relating knots with the field of {\it analysis situs}
(modern algebraic topology) newly developed by Poincar\`e,
 and ending with the Jones polynomial and related knot invariants.

In the third Chapter we discuss invariants of Conway type; that is,
invariants which have the following property: the values of the invariant
for oriented links $L_0$ and $L_-$  determine
its value for the link $L_+$ (similarly, the values of the invariant
for $L_0$ and $L_+$ determine its value for $L_-$).
The diagrams of oriented links $L_0$, $L_-$ and $L_+$ are different
only at small disks as pictured in Fig.~0.1.

\ \\

\centerline{\psfig{figure=L+L-L0.eps}}
\begin{center}
                    Fig.~0.1
\end{center}

Some classical invariants of knots turn out to be invariants
of Conway type.  \ \ ... \ \ SEE Introduction before CHAPTER I.

\chapter{Preliminaries}\label{I}
\section{Knot Theory: Intuition versus precision}\label{1}

In the XIX century Knot Theory was an experimental science.
Topology (or geometria situs) had not developed enough to
offer tools allowing precise definitions and proofs (here
Gaussian linking number is an exception). Furthermore, in
the second half of this century Knot Theory was developed
mostly by physicists and one can argue that the high
level of precision was not appreciated\footnote{This may be
a controversial statement. The precision of Maxwell was different
than that of Tait and both were physicists.}.
We outline the global history of the Knot Theory in Chapter II.
 In this chapter we deal with the 
struggle of mathematicians to understand precisely the phenomenon of
knotting.

Throughout the XIX century knots were understood as closed curves in a space
up to a natural deformation, which was described as a movement in
space without cutting and pasting. This understanding allowed scientists
to build tables of knots but didn't lead to precise methods
 allowing one to distinguish
knots which could not be practically deformed from one to another.
In a letter to O. Veblen, written in 1919, young J. Alexander expressed
his disappointment\footnote{We should remember that it was written by
a young revolutionary mathematician forgetting that he is ``standing on
the shoulders of giants." \cite{New}.}:
``When looking over Tait {\it On Knots} among other things,
He really doesn't get very far. He merely writes down all the plane projections
of knots with a limited number of crossings, tries out a few transformations
that he happen to think of and assumes without proof that if he is unable
to reduce one knot to another with a reasonable number of tries, the two
are distinct. His invariant, the generalization of the Gaussian invariant ...
for links is an invariant merely of the particular projection of the knot
that you are dealing with, - the very thing I kept running up against
in trying to get an integral that would apply. The same is true of his
`Beknottednes'."

In the famous Mathematical Encyclopedia Max Dehn and
Poul Heegaard outlined a systematic approach to topology,
in particular they precisely formulated the subject of
the Knot Theory \cite{D-H} 1907. To omit the notion
of deformation of a curve in a space (then not yet well defined)
they introduced lattice knots and the precise definition of
their (lattice) equivalence.
Later Reidemeister and Alexander considered more general
polygonal knots in a space with equivalent knots related by
a sequence of $\Delta$-moves. The definition of Dehn and Heegaard
was long ignored and only recently lattice knots are again studied.
It is a folklore result, probably never written down in detail, that the
two concepts, {\it lattice knots} and {\it polygonal knots}, are
equivalent.

\subsection{Lattice knots and Polygonal knots}

In this part we discuss two early XX century
definitions of knots and their
equivalence, by Dehn-Heegaard and by Reidemeister. In the XIX century
knots were treated from the intuitive point of view and as we mention
in Chapter II it was Heegaard in his 1898 thesis who came close to a
proof that there are nontrivial knots.

Dehn and Heegaard gave the following definition of a knot (or curve in
their terminology) and of equivalence of knots (which they call isotopy
of curves)\footnote{Translation from German due to Chris Lamm.}.

\begin{definition}[\cite{D-H}]\label{I.1.1}\ \\
A curve is a simple closed polygon on a cubical lattice.
It has coordinates $x_i,y_i,z_i$ and an {\it isotopy} of these curves
is given by:
\begin{enumerate}
\item[(i)]
Multiplication of every coordinate by a natural number,
\item[(ii)]
Insertion of an elementary square, when it does not interfere
with the rest of the polygon.
\item[(iii)]
Deletion of the elementary square.
\end{enumerate}
\end{definition}

Elementary moves of Dehn and Heegaard can be  summarized/explained as
follows: \ \ ... \ \ SEE  CHAPTER I.

\chapter{History of Knot Theory}\label{II}
  The goal of this Chapter is to present the history of ideas which lead up
to the development of modern knot theory. We will try to be more
detailed when pre-XX
century history is reported. With more recent times we are more selective,
stressing developments related to Jones type invariants of links.
Additional historical information on specific topics of Knot Theory
is  given in other chapters of the book\footnote{There are books
which treat the history of topics related to knot theory
\cite{B-L-W,Ch-M,Crowe,Die}.
J.Stillwell's textbook \cite{Stil} contains very interesting historical
digressions.}.

Knots  have fascinated people from the dawn of the human history.
We can wonder what caused a merchant living about 1700 BC. in Anatolia
and exchanging goods with Mesopotamians, to choose braids and knots as
his seal sign; Fig. 1. We can guess however that knots on seals or
cylinders appeared before proper writing was invented about 3500 BC.
\ \\
\ \\
\centerline{\psfig{figure=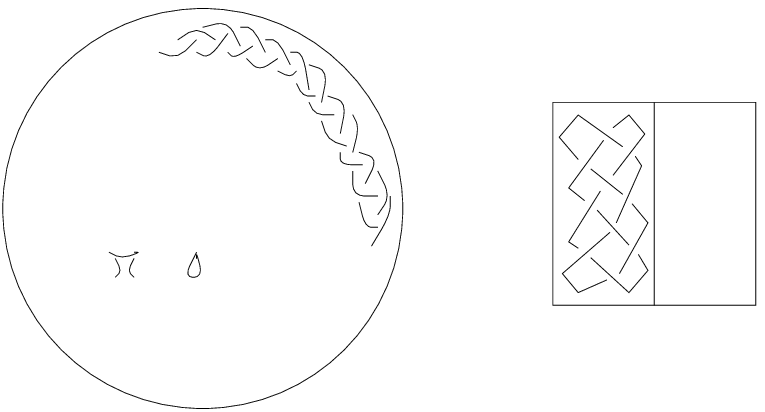}}
\begin{center}
 Figure 1; Stamp seal, about 1700 BC (the British Museum).
\end{center}
{\footnotesize {\it On the octagonal base [of hammer-handled haematite seal]
 are patterns surrounding a hieroglyphic inscription (largely erased).
Four of the sides are blank and the other four are engraved with elaborate
patterns typical of the period (and also popular in Syria) alternating
with cult scenes...}(\cite{Col}, p.93). }

   It is tempting to look for the origin of knot theory in Ancient
Greek mathematics (if not earlier). There is some justification to do so:
a Greek physician named Heraklas, who lived during the first century A.D.
and who was likely a pupil or associate of Heliodorus, wrote an essay
on surgeon's slings\footnote{
Heliodorus, who lived at the time of Trajan, also mentions in his work
knots and loops \cite{Sar}}.
Heraklas explains,  giving step-by-step instructions,
eighteen ways to tie orthopedic slings. His work survived because
Oribasius of Pergamum (ca. 325-400; physician of the emperor
 Julian the Apostate) included it toward the end of the fourth century
in his ``Medical Collections". The oldest extant manuscript of ``Medical
Collections" was made in the tenth century by the Byzantine physician
Nicetas. The Codex of Nicetas was brought to Italy in the fifteenth century by
an eminent Greek scholar,   J.~Lascaris,  a refugee from Constantinople.
Heraklas' part of the Codex of Nicetas has no illustrations,
and around 1500 an anonymous artist depicted Heraklas' knots in one of
the Greek manuscripts of Oribasus ``Medical Collections" (in Figure 2 we
reproduce the drawing of the third Heraklas knot together with its
original,   Heraklas',  description). Vidus Vidius (1500-1569),
a Florentine who became physician to Francis I (king of France, 1515-1547)
 and professor of medicine in the Coll\`ege
de France,   translated the Codex of Nicetas into Latin; it contains also
drawings of Heraklas' surgeon's slings by the Italian painter,
sculptor and architect Francesco Primaticcio (1504-1570); \cite{Da,Ra}.
\ \\ \ \\
\centerline{\psfig{figure=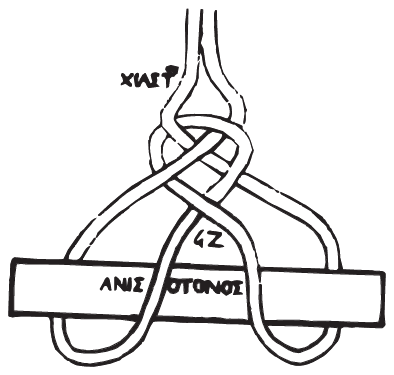, height=5.1cm}}
\begin{center}
                     Figure 2.2; The crossed noose
\end{center}

SEE CHAPTER II.

\chapter{Conway type invariants}\label{III}

\section{Conway algebras.}\label{3.1}

While considering quick methods of computing Alexander polynomial
(a classical invariant of links, see Chapter $IV$), Conway \cite{Co-1}
suggested a normalized form of it (now called the Conway or
Alexander-Conway polynomial)
and he showed that the polynomial, $\bigtriangleup_L(z)$,
satisfies the following two conditions:

\begin{enumerate}
\item [(i)] (Initial condition) \ If $T_1$ is the trivial knot then
$\bigtriangleup_{T_1}(z) = 1$.

\item [(ii)] (Conway's skein relation)\
$\bigtriangleup_{L_+}(z) - \bigtriangleup_{L_-}(z) = z\bigtriangleup_{L_0},$
where $L_+, L_-$ and $L_0$ are diagrams of oriented links
which are identical except for the part presented in Fig.~1.1.

\end{enumerate}
The conditions $(i)$ and $(ii)$ define the Conway polynomial (or, maybe
more properly, Alexander-Conway polynomial)
$\bigtriangleup_L(z)$ uniquely, see \cite{Co-1,K-1,Gi,B-M}.

In fact the un-normalized version of the
 formula $(ii)$ was noted by J.W.~Alexander in
his original paper
introducing the polynomial \cite{Al-3} 1928.
Alexander
polynomial was defined up to invertible elements, $\pm t^i$,
 in the ring of Laurent polynomials, $Z[t^{\pm 1}]$),
so the formula (ii) was not easily
available for a computation of the polynomial.

In the Spring of 1984 V.~Jones, \cite{Jo-1,Jo-2},
showed that there exists an invariant $V$ of links which is a Laurent
polynomial with respect to the variable $\sqrt{t}$ which satisfies
the following conditions:
\begin{enumerate}
\item[(i)]
$V_{T_1}(t) = 1,$

\item[(ii)]
${\displaystyle\frac{1}{t}V_{L_+}(t) - tV_{L_-}(t) = (\sqrt{t} -
\frac{1}{\sqrt{t}})V_{L_0}(t).}$
\end{enumerate}

These two examples of invariants were a base for a thought
that there exists an invariant (of ambient isotopy) of
oriented links which is a Laurent polynomial
$P_L(x,y)$ of two variables and which satisfies the following conditions:

\begin{enumerate}
\item[(i)]
$P_{T_1}(x,y) = 1$

\item[(ii)]
$xP_{L_+}(x,y) + yP_{L_-}(x,y) = zP_{L_0}(x,y).$
\end{enumerate}

Indeed, such an invariant exists and it was discovered a few months after
the Jones polynomial, in July-September of 1984, by four groups of
mathematicians:
R.~Lickorish and K.~Millett,
J.~Hoste, A.~Ocneanu as well as by  P.~Freyd and D.~Yetter \cite{FYHLMO}
(independently, it was discovered in November-December
of 1984 by J.~Przytycki
and P.~Traczyk \cite{P-T-1}).
We call this polynomial the Jones-Conway or Homflypt
 polynomial\footnote{HOMFLYPT is the acronym after the initials
of the inventors: Hoste,
Ocneanu, Millett, Freyd, Lickorish, Yetter, Przytycki and Traczyk.
We note also some other names that
are used for this invariant:
FLYPMOTH, HOMFLY, the generalized Jones polynomial,
two variable Jones polynomial,
twisted Alexander polynomial and skein-polynomial.}.
\ \\
\ \\
\begin{center}
\begin{tabular}{c}
\includegraphics[trim=0mm 0mm 0mm 0mm, width=.5\linewidth]
{L+L-L0.eps}\\
\end{tabular}
 \ \\
\ \\
Fig.~1.1
\end{center}

Instead of looking for polynomial invariants of links related as in Fig.~1.1
we can approach the problem from a more general point of view.
Namely, we can look for universal invariants of links which
have the following property:
a given value of the invariant for $L_+$ and $L_0$ determines
the value of the invariant for $L_-$, and similarly:
if we know the value of the invariant
for $L_-$ and $L_0$ we can find its value for $L_+$.

The invariants with this
property are called Conway type invariants.
We will develop these ideas in the present chapter of the book
which is based mainly on a joint paper of Traczyk and the author
\cite{P-T-1,P-1}.

SEE CHAPTER III.

\chapter{Goeritz and Seifert matrices}\label{IV}
The renessance of combinatorial methods in Knot Theory which can
be traced back to Conway's paper \cite{Co-1} and which bloomed after
the Jones breakthrough \cite{Jo-1} with Conway type invariants and
Kauffman approach (Chapter III), had its predecessor
in 1930th \cite{Goe,Se}.
Goeritz matrix of a link can be defined purely combinatorially and
is closely related to Kirchhoff matrix of an electrical network.
Seifert matrix is a generalization of the Goeritz matrix and, even
historically, its development was mixing combinatorial and topological
methods.
\section{Goeritz matrix and signature of a link.}\label{secIV1}

Apart from Conway type invariants and Kauffman approach, there are
other combinatorial methods of examining knots. One of them was discovered
in the 30-ties by L.~Goeritz and we present it in this section.
Goeritz \cite{Goe} showed how to associate a quadratic form
to a diagram of a link and moreover how to use this form
to get algebraic invariants of the knot
(the signature of this form, however, was not an invariant of the knot).
Later, Trotter \cite{Tro-1} applied a form of Seifert (see
Section 2) to introduce another quadratic form, the signature of which
was an invariant of links.

Gordon and Litherland \cite{G-L} provided a unified
approach to Goeritz and Seifert/Trotter forms. They showed how to use
the form of Goeritz to get (after a slight modification)
the signature of a link.

We begin with a purely combinatorial description
of the matrix of Goeritz and of the signature of a link.
This description is based on \cite{G-L} and
\cite{Tral-1}.

 \begin{definition}\label{c3:1.1}
Assume that $L$ is a diagram of a link.
Let us checkerboard color the complement of the diagram in the projection
plane $R^2$, that is, color in black and white
the regions into which the plane is divided by the diagram\footnote{This
 (checkerboard) coloring was first used
by P. G. Tait in 1876/7, see Chapter II.}.
 We assume that the unbounded region of $R^2\setminus L$
is colored white and it is denoted by $X_0$ while the other
white regions are denoted by $X_1, \ldots, X_n$.
Now, to any crossing, $p$, of $L$ we associate either $+1$ or $-1$
according to the convention in Fig.~1.1.
We denote this number by $\eta (p)$.  \\
\ \\
\begin{center}
\begin{tabular}{c}
\includegraphics[trim=0mm 0mm 0mm 0mm, width=.5\linewidth]
{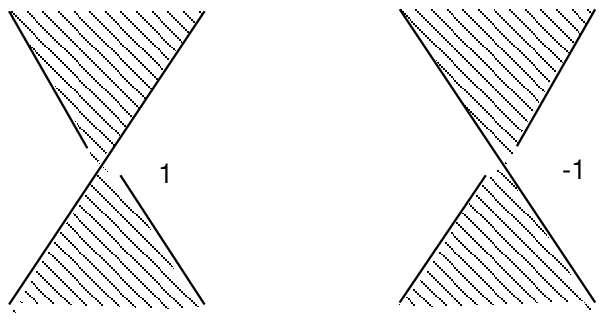}\\
\end{tabular}
\\
Fig. 1.1
\end{center}
Let $G' = \{ g_{i,j}\}^n_{i,j=0}$, where
$$g_{i,j} = \left\{
\begin{array}{ll}
-\sum_{p}\eta (p)
& \mbox{ for $i\neq j$, where the summation extends }\\
 \ \ &\mbox{ over crossings which connect $X_i$ and $X_j$}\\
-\sum_{k=0,1,\ldots,n; k\neq i} g_{i,k}&\mbox{ if $ i=j$}\\
\end{array}
\right.
$$

The matrix $G'=G'(L)$ is called the unreduced Goeritz matrix
of the diagram $L$.
The reduced Goeritz matrix (or shortly Goeritz matrix) associated to
the diagram $L$ is
the matrix $G=G(L)$ obtained by removing the first
row and the first column of $G'$.
\end{definition}
\ \\
{\bf Theorem}(\cite{Goe,K-P,Ky})\\
Let us assume that $L_1$ and $L_2$ are two diagrams of a given link.
Then the matrices $G(L_1)$ and $G(L_2)$ can be obtained one from the other
in a finite number of the following elementary equivalence
operations:... \\
 SEE CHAPTER IV.

\chapter{Graphs and links}\label{V}
We present in this Chapter several results which demonstrate a close
connection and useful exchange of ideas between graph theory and
knot theory. These disciplines were shown to be related from the
time of Tait (if not Listing) but the great flow of ideas started
only after Jones discoveries. The first deep relation in this new
trend was demonstrated by Morwen Thistlethwaite and we describe several
results by him in this Chapter. We also present
results from two preprints \cite{P-P-0,P-18},
in particular we sketch
two generalizations of the Tutte polynomial of graphs,
 $\chi(S)$, or, more precisely, the deletion-contraction method which Tutte
polynomial utilitize.
The first generalization considers, instead of graphs, general objects
called setoids or group systems. The  second one deals with completion
of the expansion of a graph with respect to subgraphs.
 We are motivated here by finite type
invariants of links developed by Vassiliev and Gusarov along the line
presented in \cite{P-9} (compare Chapter IX).
The dichromatic Hopf algebra, described in Section 2, have its
origin in Vassiliev-Gusarov theory mixed with work of G. Carlo-Rota and
his former student W. Schmitt.
\section{Knots, graphs and their polynomials}\label{V.1}
In this section we discuss relations between graph and knot theories.
We describe several applications of graphs to knots.
In particular we consider various
interpretations of the Tutte polynomial of graphs in knot theory.
This serves as an introduction to the subsequent sections
where we prove two of the classical conjectures of Tait.
In the present section we rely mostly on
\cite{This-1,This-5} and \cite{P-P-1}.

By a graph $G$ we understand a finite set $V(G)$ of vertices
together with a finite set of edges $E(G)$. To any edge we associate
a pair of (not necessarily distinct) vertices which we call
endpoints of the edge.
We allow that the graph $G$ has multiple edges and loops
(Fig.1.1)\footnote{In terms of algebraic topology a graph is
a 1-dimensional CW-complex. Often it is called a pseudograph and the
word ``graph" is reserved for a 1-dimensional simplicial complex,
that is,  loops and multiple edges are not allowed.
We will use in such a case the
term a {\it simple} (or {\it classical}) graph.
If multiple edges are allowed but loops
are not we use often the term a {\it multigraph}, \cite{Bo-1}}.
 A loop is an edge with one endpoint.

\centerline{{\psfig{figure=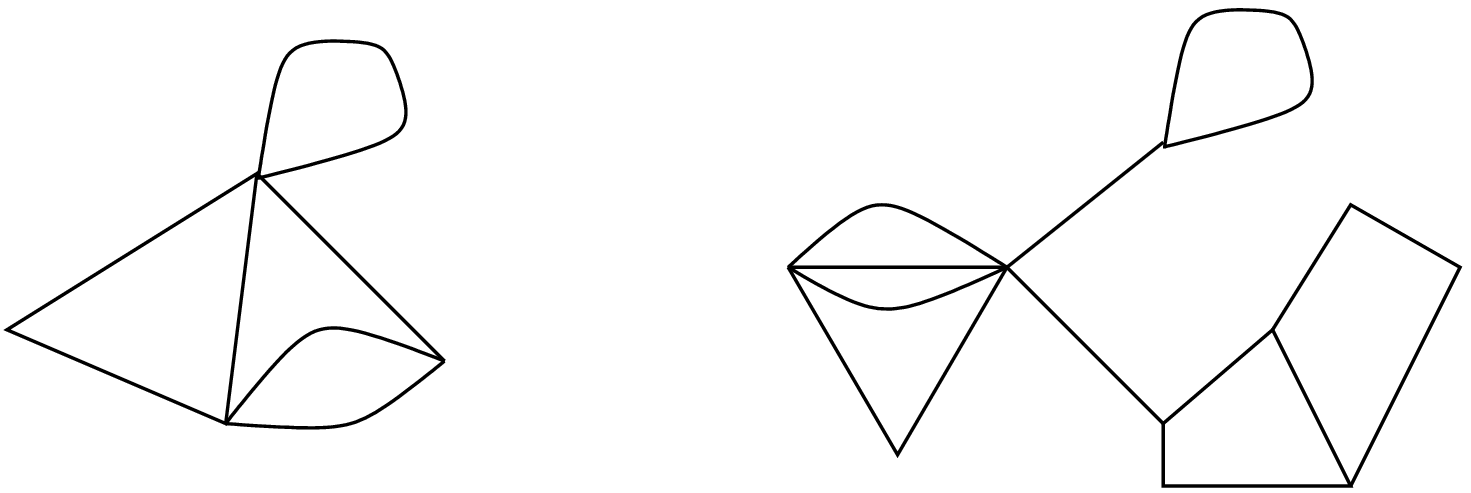,height=2.5cm}}}
\begin{center}
Fig.~1.1
\end{center}

By $p_0(G)$ we denote the number of components of the graph $G$
and by $p_1(G)$ we denote its cyclomatic number, i.e.~the minimal number
of edges which have to be removed
from the graph in order to get a graph without
cycles.\footnote{In terms of algebraic
topology $p_0(G)$ and $p_1(G)$ are equal to
dimensions of homology groups $H_0(G)$ and $H_1(G)$, respectively. In
this context the notation $b_0$ and $b_1$ is used and numbers are called
the Betti numbers.}
A connected graph without cycles (i.e.~$p_0=1$, $p_1=0$) is called a tree.
If $G$ has no cycles , i.e.~$p_1=0$, then the graph $G$ is called a forest.
By a spanning tree (resp. forest) of the graph $G$ we understand
a tree (resp. forest) in $G$ which contains all vertices of $G$.
By an isthmus of $G$ we understand an edge of $G$,
removal of which increases the number of components of the graph.

To a given graph we can associate a polynomial in various ways.
The first such a polynomial, called the chromatic polynomial of a graph,
was introduced by Birkhoff in 1912 \cite{Birk-1}\footnote
{J.B.Listing, in 1847\cite{Lis}, introduced polynomial of knot diagrams.
For a graph $G$, the Listing polynomial, denoted by $JBL(G)$,
 can be interpreted as follows:
$JBL(G)=\Sigma a_i(G)x^i$ where $a_i(G)$ is the number of vertices in
$G$ of valency $i$.}.
For a natural number $\lambda$, the chromatic polynomial, denoted
by $C(G,\lambda)$, counts the number
of possible ways of coloring the vertices of $G$ in $\lambda$ colors
in such a way that each edge has endpoints colored in different colors
(compare Exercise 1.14).
The chromatic polynomial was generalized
by Whitney and Tutte \cite{Tut-1}.\ ... \\
 SEE CHAPTER V.

\chapter{Fox $n$-colorings,Rational moves, Lagrangian tangles 
and Burnside groups}\label{VI}
\section{Fox $n$-colorings}\label{VI.1}
\markboth{\hfil{\sc Fox colorings...Burnside groups}\hfil}
{\hfil{\sc Fox $n$-colorings}\hfil}

Tricoloring of links, discussed in Chapter I, can be generalized, after Fox,
\cite{F-1}; Chapter 6, \cite{C-F}, Chapter VIII, Exercises 8-10,
\cite{F-2},
to n-coloring of links as follows:
\begin{definition}\label{2.1}
We say that a link diagram $D$ is n-colored if every
arc is colored by one of the numbers $0,1,...,n-1$ in such a way that
at each crossing the sum of the colors of the undercrossings is equal
to twice the color of the overcrossing modulo $n$.
\end{definition}

The following properties of n-colorings, can be proved in a similar way
as for the tricoloring properties. However, an elementary proof
of the part (g) is more involved and requires an interpretation
of n-colorings using the Goeritz matrix \cite{Ja-P}, or Lagrangian 
tangle idea (developed later in this chapter).
\begin{lemma}\label{2.2}
\begin{enumerate}
\item
[(a)] Reidemeister moves preserve the number of n-colorings, $col_n(D)$,
      thus it is a link invariant,
\item
[(b)] if $D$ and $D'$ are related by a finite sequence of $n$-moves, then
$col_n(D)=col_n(D')$,
\item
[(c)] $n$-colorings form an abelian group, $Col_n(D)$
(it is also a $Z_n$-module),
\item
[(d)] if $n$ is a prime number, then $col_n(D)$ is a power of $n$ and for a
link with $b$ bridges: \ \ $b\geq log_n(col_n(L))$,
\item
[(e)] $col_n(L_1) col_n(L_2)=n(col_n(L_1\# L_2))$,
\item
[(f)] if $n$ is a prime odd number then among the four numbers
$col_n(L_+), col_n(L_-),$ $col_n(L_0)$ and $col_n(L_{\infty}),$
three are equal one to another and the fourth is
either equal to them or $n$ times bigger,
\end{enumerate}
More generally: If $L_0, L_1, ..., L_{n-1},L_{\infty}$ are $n+1$ diagrams
generalizing the four diagrams from (f); see Fig.2.1 then:
\begin {enumerate}
\item
[(g)] if $n$ is a prime number then among the $n+1$ numbers
$col_n(L_0), col_n(L_1),...,col_n(L_{n-1})$ and $col_n(L_{\infty})$
$n$ are equal one to another and the $(n+1)$'th is $n$ times bigger,
\item
[(h)]
if $n$ is a prime number, then $u(K)\geq log_n(col_n(K))-1$.
\end{enumerate}
\end{lemma}
\ \\

\centerline{\psfig{figure=L0L1L2L3Linf.eps}}
\begin{center}
Fig. 2.1
\end{center}

\begin{corollary}\label{2.3}
\begin{enumerate}
\item
[(i)] For the figure eight knot, $4_1$, one has $col_5(4_1)=25$, so the figure
eight knot is a nontrivial knot; compare Fig.2.2.
\item
[(ii)] $u(4_1\# 4_1)=2$.
\end{enumerate}
\end{corollary}

\centerline{\psfig{figure=4-1fivecol.eps,height=3.6cm}}

\ \\ \ ... \ \ SEE CHAPTER VI.

\chapter{Symmetries of links}\label{VII}
In this  chapter we examine finite group actions on $S^3$
which map a given link onto itself.
For example, a torus link of type $(p,q)$ (we call it $T_{p,q}$)
is preserved by an action of a group $Z_p\oplus Z_q$ on $S^3$
(c.f.~ Exercise 0.1 and Figure 0.1, for the torus link of type $(3,6)$).\ \\
\centerline{{\psfig{figure=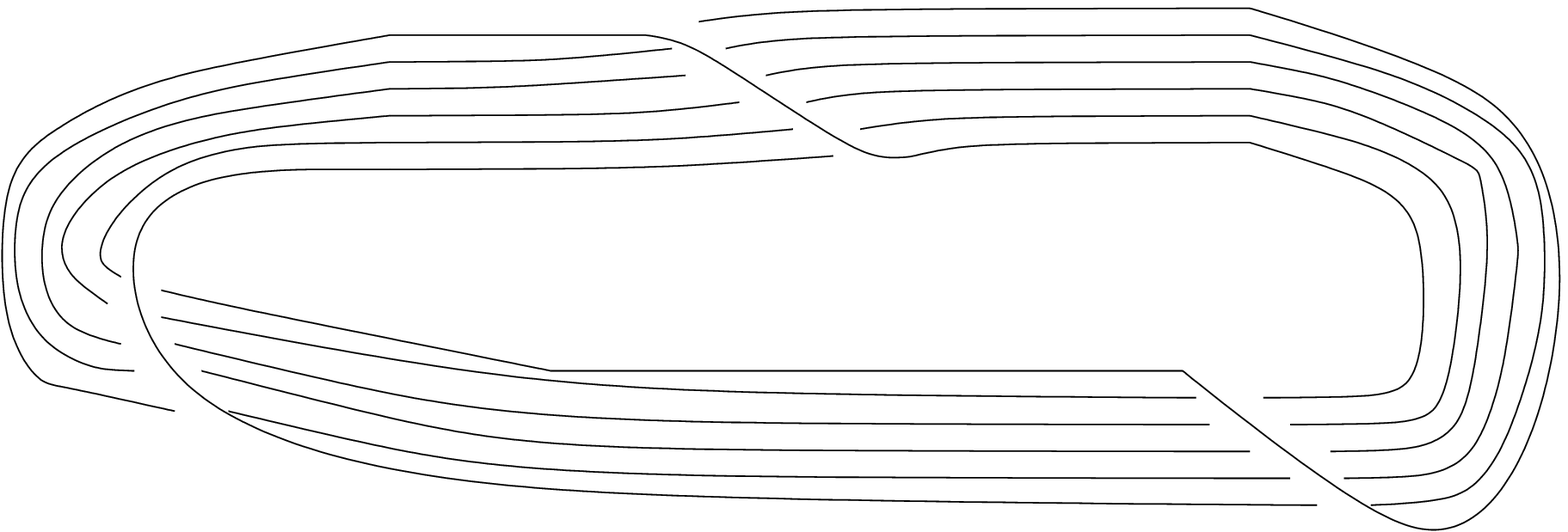,height=3.5cm}}}
\begin{center}
Fig.~0.1
\end{center}


Subsequently, we will focus on the action of a cyclic group
$Z_n$. We will mainly consider the case of an action on $S^3$
with a circle of fixed points.
The new link invariants, which we have discussed in previous chapters,
provide efficient criteria
for examining such actions.
The contents of this chapter is mostly based on
 papers of Murasugi \cite{M-6}, Traczyk \cite{T-1}
and of the author \cite{P-4,P-32}.

\ \\
{\bf Exercise VII.0.1}
Let $S^3 = \{ z_1,z_2\in C\times C:|z_1|^2+|z_2|^2 = 1\}$.
Let us consider an action of $Z_p\oplus Z_q$ on $S^3$
which is generated by  $T_p$ and $T_q$,
where $T_p(z_1,z_2) = (e^{2\pi i/p} z_1,z_2)$ and
$T_q(z_1,z_2) = (z_1, e^{2\pi i/q}z_2)$.
Show that this action preserves torus link of type
$(p,q)$.
This link can be described as the following set
$\{ (z_1,z_2)\in S^3:\ z_1 = e^{2\pi i(\frac{t}{p} +\frac{k}{p})}, z_2 =
e^{2\pi it/q}\}$, where $t$ is an arbitrary real number and $k$
is an arbitrary integer. \\
Show that if $p$ is co-prime to $q$ then $T_{p,q}$ is a knot that can be
parameterized by:
$$R\ni t\mapsto (e^{2\pi it/p},\ e^{2\pi it/q})\in S^3\subset C^2.$$

\section{Periodic links.}\label{VII.1}

\begin{definition}\label{1.1}
A link is called $n$-periodic if there exists an action of $Z_n$
on $S^3$ which preserves the link and the set of fixed points of the action
is a circle disjoint from the link.
If, moreover, the link is oriented then we assume that
the generator of $Z_n$ preserves the orientation of the link
or changes it globally (that is on every component).
\end{definition}

New polynomials of links provide strong periodicity criteria.

Let ${\cal R}$ be a subring\footnote{
${\cal R}$ is an example of a Rees algebra of $Z[a^{\mp 1}]$ with
respect to the ideal generated by $a+a^{-1}$; \cite{Ei}.}
 of the ring
$Z[a^{\mp 1},z^{\mp 1}]$ generated by $a^{\mp  1}$, $z$ and
$\frac{a+a^{-1}}{z}$.
Let us note that  $z$ is not invertible in ${\cal R}$.

\begin{lemma}\label{VII:1.2}
For any link $L$ its Jones-Conway polynomial $P_L(a,z)$
is in the ring ${\cal R}$.
\end{lemma}

\begin{proof}
For a trivial link $T_n$ with $n$ components we have
$P_{T_n}(a,z)=(\frac{a+a^{-1}}{z})^{n-1}\in {\cal R}$.
Further, if
$P_{L_+}(a,z)$ (respectively $P_{L_-}(a,z)$) and $P_{L_0}(a,z)$
are in ${\cal R}$ then $P_{L_-}(a,z)$ (respectively $P_{L_+}(a,z)$)
is in ${\cal R}$ as well.
This observation enables a standard induction to conclude Lemma 1.2.
Now we can formulate our criterion for $n$-periodic links. It has especially
simple form for a prime period (see Section 2 for a more general statement).
\end{proof}

\begin{theorem}\label{VII:1.3}
Let $L$ be an $r$-periodic oriented link and assume that $r$ is
a prime number. Then the Jones-Conway polynomial $P_L(a,z)$
satisfies the relation
$$P_L(a,z)\equiv P_L(a^{-1},z)\mbox{ mod }(r,z^r)$$
where $(r,z^r)$ is an ideal in ${\cal R}$ generated by $r$ and
$z^r$.
\end{theorem}

In order to apply Theorem 1.3 effectively, we need the following
fact.
\begin{lemma}\label{VI:1.4}
Suppose that $w(a,z)\in{\cal R}$ is written in the form
$w(a,z) = \sum_i v_i(a)z^i$, where
$v_i(a)\in Z[a^{\mp 1}]$.
Then $w(a,z)\in(r,z^r)$ if and only if
for any $i\leq r$ the coefficient $v_i(a)$ is in the ideal
$(r,(a+a^{-1})^{r-i})$.
\end{lemma}

\ \\ \ ... \ \ SEE CHAPTER VII.

\chapter{Different links with the same Jones type polynomials}\label{VIII}
\ \\
\begin{quotation}
ABSTRACT.
\baselineskip=10pt
We describe, in this chapter,
three methods of constructing different links with the same
Jones type invariant. All three can be thought as generalizations of mutation.
The first combines the satellite construction with mutation. The second
uses the notion of rotant, taken from the graph theory, the third, invented
by Jones, transplants into knot theory the idea of the Yang-Baxter equation
with the spectral parameter (idea employed by Baxter  in the theory of solvable
models in statistical mechanics). We extend the Jones result and
relate it  to Traczyk's work on rotors of links. We also show further
applications of the Jones idea, e.g. to 3-string links in the solid torus.
We stress the fact that ideas coming from various areas of mathematics (and
theoretical physics) has been fruitfully used in knot theory, and vice versa.
\\
\end{quotation}
\ \\
{\Large \bf 0\ \  Introduction}\\ \ \\
When at spring of 1984, Vaughan Jones introduced his
 Laurent) polynomial invariant of links, $V_L(t)$,
he checked immediately that it distinguishes many knots which were not
taken apart by the Alexander polynomial, e.g. the right handed trefoil
knot from the left handed trefoil knot, and the square knot from
the granny knot; Fig. 0.1.

\par\vspace{1cm}
\begin{center}
\begin{tabular}{c}
\includegraphics[trim=0mm 0mm 0mm 0mm, width=.65\linewidth]
{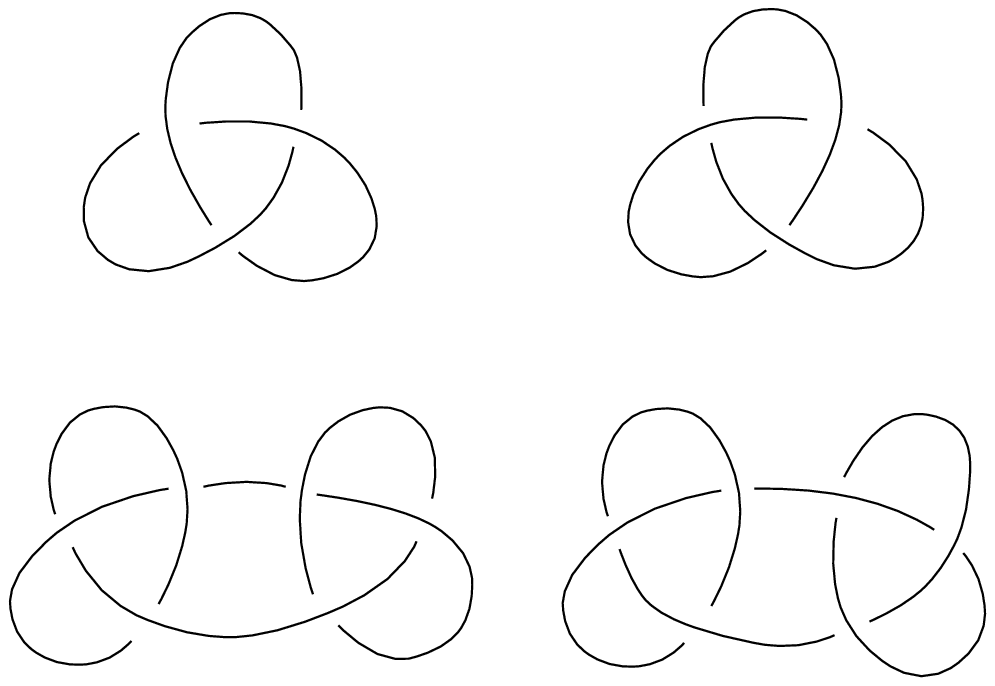}
\end{tabular}
\par\vspace{0.5cm}
Fig. 0.1
\end{center}

Jones also noticed that his polynomial is not universal. That is, there
are  different knots with the same polynomial; e.g. the Conway and
Kinoshita-Terasaka knots; Fig. 0.2.
\par\vspace{1cm}
\begin{center}
\begin{tabular}{cc}
\includegraphics[trim=0mm 0mm 0mm 0mm, width=.3\linewidth]
{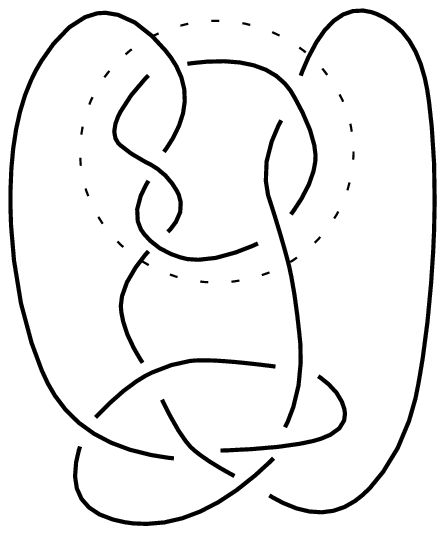}
&\includegraphics[trim=0mm 0mm 0mm 0mm, width=.3\linewidth]
{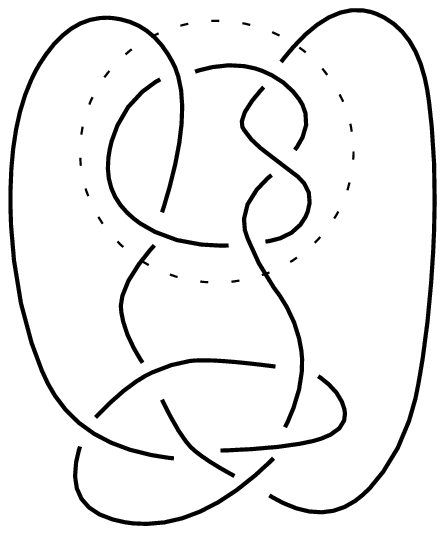}
\\
\end{tabular}
\par\vspace{0.5cm}
Fig. 0.2
\end{center}

Then Jones asked the fundamental
question whether there exists a nontrivial knot with the trivial
polynomial. Twenty years later this is still an open problem and specialists
differ in their opinion whether the answer is yes or no. In this section,
we concentrate on more accessible problem: how to construct different
links with the same Jones polynomial. It may shed some light into the
Jones question. We follow mostly \cite{P-16} in our exposition. 
\ \ \\ \ ... \ \ SEE CHAPTER VIII.

\chapter{Skein modules}\label{IX}
We  describe in this chapter the idea of building an
algebraic topology based on knots (or more generally on the position of
embedded objects). That is, our basic building blocks
are considered up to ambient isotopy (not homotopy or homology).
For example, one should start from knots in 3-manifolds, surfaces in
4-manifolds, etc. However our theory is, until now, developed only
for the case of links in 3-manifolds, with only a glance towards
4-manifolds. The main object of the theory is a {\it skein module} and
we devote this chapter mostly to description of skein modules in 3-dimensional
manifolds.
In this book we outline the theory of skein modules often giving
only ideas and outlines of proofs. The author is preparing a monograph
 devoted exclusively to skein modules and their ramifications
\cite{P-30}.


\section{ History of skein modules}\label{IX.1}

H.~Poincar\'e, in his paper ``Analysis situs" (1895), defined abstractly
homology groups starting from formal linear combinations of simplices,
choosing cycles and dividing them by relations coming from boundaries
\cite{Po}\footnote{Before Poincar\'e the only similar construction
was the formation of ``divisors" on an algebraic curve by Dedekind and
Weber \cite{D-W}, that is the idea of considering formal linear
combinations of points on an algebraic curve, modulo relations yielded
by rational functions on the curve.}.

The idea behind skein modules is to use links instead of cycles (in the
case of a 3-manifold).
More precisely we divide the free module
generated by links by properly chosen (local) skein relations.

Skein modules have their origin in the observation made by Alexander
(\cite{Al-3}, 1928)\footnote {We can argue further 
that Alexander was motivated by the
chromatic polynomial introduced in 1912 by George David Birkhoff
\cite{Birk-1}. Compare also the letter of Alexander to Veblen
(\cite {A-V}, 1919) discussed in Section 2.} that his polynomials of
three links $L_+, L_-$ and $L_0$ in $S^3$ are linearly related
(here $L_+,L_-$ and $L_0$ denote three
links which are identical except in a small ball as shown
in Fig. 1.1). Conway rediscovered the Alexander observation and
normalized the Alexander polynomial so that it satisfies the skein relation
$$\Delta_{L_+}(z) - \Delta_{L_-}(z) = z\Delta_{L_0}(z)$$
 (\cite{Co-1}, 1969). In the late seventies Conway advocated the idea of
considering  the free $Z[z]$-module over oriented
links in an oriented 3-manifold and
dividing it by the submodule generated by his skein relation
\cite{Co-2} (cited in \cite{Gi}) and \cite{Co-3}
(cited in \cite{Ka-1}). However, there is no published account of
the content of Conway's
talks except when $S^3$ or its submanifolds are analyzed.
The original name Conway used for this object was ``linear skein".
\centerline{\psfig{figure=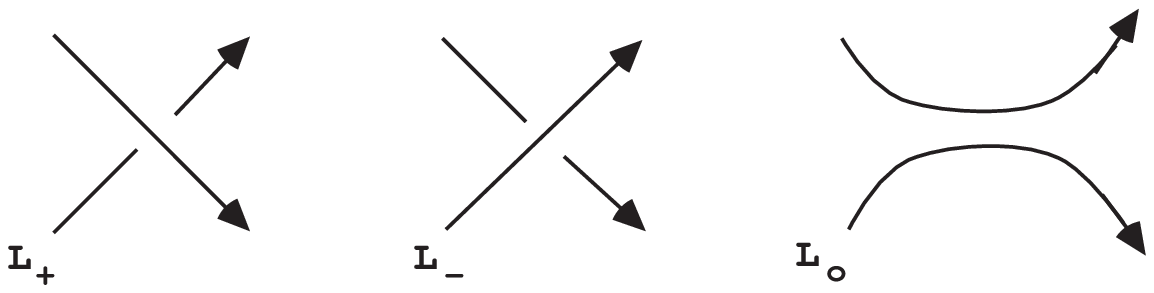,height=1.5cm}}
\begin{center}
Figure 1.1
\end{center}
Conway's idea was then pursued by Giller \cite{Gi}
(who computed the linear skein of a tangle), and Kauffman
\cite{Ka-1,Ka-8,Ka-3}, as well as Lickorish and Millett \cite{L-M-1}
(for subspaces of $S^3$).

In graph theory, the idea of forming a ring of graphs and dividing
it by an ideal generated by local relations was developed by W.~Tutte in
his 1946 PhD thesis \cite{Tut-1}, but the relation to knot theory
was observed much later.

The theory of Hecke algebras, as introduced by N.~Iwahori (\cite{Iw},1964),
is closely connected to the theory of skein modules, however the relation
of Hecke algebras to knot theory was noticed by V.~Jones in 1984,
20 years after Iwahori's and Conway's work
(it was crucial for Jones-Ocneanu construction of Markov traces).

The Temperley-Lieb algebra (\cite{T-L},1971) is related to the Kauffman
bracket skein module of the tangle, but any relation to knot theory was
again observed first by Jones in 1984.

At the time when I introduced skein modules, in April of 1987, I knew
the fundamental paper of Conway \cite{Co-1}, and \cite{Gi,Ka-8,L-M-1}
as well as \cite{Li-10}. However, the most stimulating paper
for me was one by J.~Hoste and M.~Kidwell \cite{Ho-K} about invariants
of colored links, which 
\ \\ \ ... \ \ SEE CHAPTER IX.

\chapter{Khovanov Homology: categorification of the Kauffman 
bracket relation}\label{X}
\centerline{Bethesda, October 31, 2004}

\baselineskip=14pt
\ \\
\begin{Large}{\bf Introduction} \end{Large}\\
Khovanov homology offers a nontrivial generalization of the Jones 
polynomials of links in $R^3$ (and of the Kauffman bracket skein modules 
of some 3-manifolds). 
In this chapter we define Khovanov homology of links in $R^3$
and generalize the construction into links in an $I$-bundle
over a surface.
We use Viro's approach to 
construction of Khovanov homology \cite{V-3}, and utilize 
the fact that one works 
with unoriented diagrams (unoriented framed links) in which case there
is a long exact sequence of Khovanov homology. Khovanov homology, over 
the field $\Q$, is a categorification of the Jones polynomial (i.e. one 
represents the Jones polynomial as the 
generating function of Euler characteristics). However, for integral 
coefficients Khovanov homology almost always has torsion. The first 
part of the chapter  is devoted to the construction of torsion in 
Khovanov homology. In the second part we analyze the 
thickness of Khovanov homology and reduced Khovanov homology.
In the last part we generalize Khovanov homology to links in 
$I$-bundles over  surfaces and demonstrate that for oriented surfaces 
we categorify the element of the Kauffman bracket skein module 
represented by the considered link.
We follow the exposition given in \cite{AP,APS-2,APS-3}.

Chapter X is organized as follows.\ 
In the first section we give the definition of Khovanov homology and
its basic properties. 

In the second section we prove that adequate link diagrams with 
an odd cycle property have $\Z_2$-torsion in Khovanov homology.

In the third section we discuss torsion in the Khovanov homology of an 
 adequate link diagram with an even cycle property.

In the fourth section we prove Shumakovitch's theorem that 
prime, non-split alternating links different from the trivial knot and 
the Hopf link have $\Z_2$-torsion in Khovanov homology.
We generalize this result to a class of adequate links.

In the fifth section we generalize result of E.S.Lee about 
the Khovanov homology of alternating links (they are 
$H$-thin\footnote{We also 
outline a simple proof of Lee's result\cite{Lee-1,Lee-2} 
that for alternating links Khovanov 
homology yields the classical signature, see Remark 1.6.}). 
We do not assume rational coefficients in this 
generalization.
We use Viro's exact sequence of Khovanov homology to extend Lee's 
results to  almost alternating diagrams and $H$-$k$-thick links.

In the sixth section we compute the Khovanov homology for a connected sum of 
$n$ copies of Hopf links and construct a short exact sequence 
of Khovanov homology 
involving a link and its connected sum with the Hopf link. By showing that 
this sequence splits, we answer the question asked by Shumakovitch. 

In the seventh section we notice that the results of sections 5 and 6 can be 
adapted to reduced and co-reduced Khovanov homology. Finally, we show that 
there is a long exact sequence
connecting reduced and co-reduced Khovanov homology with unreduced homology.
We illustrate our definitions by computing reduced and co-reduced 
Khovanov homology for the left handed trefoil knot.

In the eight section we define Khovanov homology for links in $I$-bundles 
over surfaces. For links in the product 
$F \times I$ we stratify Khovanov homology 
in such a way that it categorifies the Kauffman bracket skein module 
of $F \times I$. We end the section by computing stratified 
Khovanov homology for a trefoil knot diagram in an annulus.  

\section{Basic properties of Khovanov homology}\label{1}
\markboth{\hfil{\sc Khovanov Homology }\hfil}
{\hfil{\sc Basic properties of Khovanov homology of links}\hfil}

The first spectacular application of the Jones polynomial 
(via Kauffman bracket relation) was the solution of Tait conjectures on 
alternating diagrams and their generalizations to adequate diagrams 
(see Chapter V).
Our method of analyzing torsion in Khovanov homology has its root in 
work related to solutions of Tait conjectures \cite{K-5,M-4,This-3}. 

Recall that the Kauffman bracket polynomial $<D>$ of a link diagram $D$
is defined by the skein relations
$<\parbox{0.6cm}{\psfig{figure=L+nmaly.eps,height=0.6cm}}> =
A <\parbox{0.6cm}{\psfig{figure=L0nmaly.eps,height=0.6cm}}> + A^{-1}
<\parbox{0.6cm}{\psfig{figure=Linftynmaly.eps,height=0.6cm}}>$
and $<D\sqcup \bigcirc> =(-A^2-A^{-2})<D>$ and the normalization
$<\bigcirc> = 1$. The categorification of this invariant (named
by Khovanov {\it reduced homology}) is discussed in Section 7.
For the (unreduced) Khovanov homology we use the version of
the Kauffman bracket polynomial
normalized to be $1$ for the empty link (we use the notation $[D]$
in this case).

\begin{definition}[Kauffman States]\label{1.1} \ \\
Let $D$ be a diagram\footnote{We think of the 3-ball $B^3$ as $D^2\times I$ 
and the diagram is drawn on the disc $D^2$. In Section 8  
we demonstrate, after\cite{APS-2}, that 
the theory of Khovanov homology can be extended to links in an oriented 
3-manifold $M$ that is a bundle over a surface 
$F$ ($M= F \tilde \times I$). 
If $F$ is orientable then
$M= F\times I$. If $F$ is unorientable then $M$ is a twisted 
$I$ bundle over $F$ (denoted by $F \hat \times I$). Several results 
of this chapter are valid for the Khovanov homology of links 
in $M= F \tilde \times I$.} of an unoriented, framed link in a 3-ball 
$B^3$.  A Kauffman state $s$
of $D$ is a function from the set of crossings of $D$ to
the set $\{+1,-1\}$. Equivalently,
we assign to each crossing of $D$ a marker
according to the following convention:\\
\ \\
\centerline{\psfig{figure=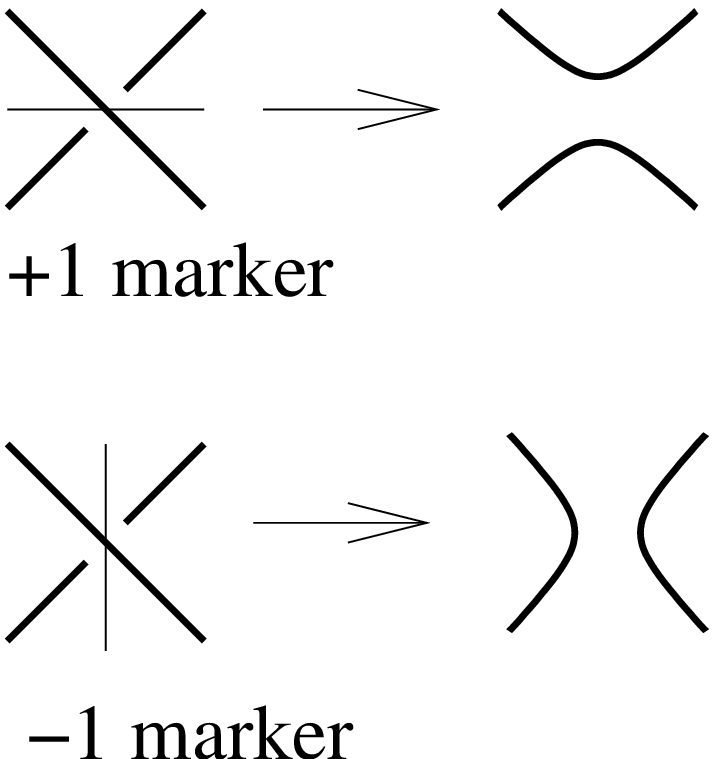,height=4.5cm}}
\begin{center}
Fig. 1.1; markers and associated smoothings
\end{center}

By $D_s$ we denote the system of circles in the diagram 
obtained by smoothing all
crossings of $D$ according to the markers of the state $s$, Fig. 1.1.\\
By $|s|$ we denote the number of components of $D_s$.
\end{definition}

In this notation the Kauffman bracket polynomial of $D$ is given by
the state sum formula:
$[D] = (-A^2-A^{-2})<D> = \sum_s A^{\sigma(s)}(-A^2-A^{-2})^{|s|}$,
where $\sigma(s)$ is the number
of positive markers minus the number of negative markers in the state $s$.


Convenient way of defining Khovanov homology is
 (as noticed by Viro \cite{V-3}) 
to consider enhanced Kauffman states\footnote{In Khovanov's original 
approach every circle of a Kauffman state was decorated by a 
2-dimensional module $A$ (with basis {\bf 1} and $X$) 
with the additional structure  of Frobenius algebra. As an algebra 
$A=\Z[X]/(X^2)$ and comultiplication is given by 
$\Delta({\mbox{\textbf{1}}})=X\otimes {\mbox{\textbf{1}}} + 
{\mbox{\bf{1}}}\otimes X$ and $\Delta(X) = 
X \otimes X$. Viro uses $-$ and $+$ in 
place of {\bf 1} and $X$. }.
\begin{definition}\label{1.2}
An enhanced Kauffman state $S$ of
an unoriented framed link diagram
$D$ is a Kauffman state $s$ with an additional assignment of
$+$ or $-$ sign to each circle of $D_s.$
\end{definition}
Using enhanced states we express the Kauffman bracket polynomial as
a (state) sum of monomials which is important in the definition
of Khovanov homology we use. We have
$[D] = (-A^2-A^{-2})<D> = \sum_S (-1)^{\tau(S)}A^{\sigma(s)+2\tau(S)}$,
where $\tau(S)$ is the number of positive circles minus 
the number of negative circles in the enhanced state $S$.

\begin{definition}[Khovanov chain complex]\label{1.3}\ \
\begin{enumerate}
\item[(i]) Let ${\cal S}(D)$ denote the set of enhanced Kauffman states 
of a diagram $D$, and let ${\cal S}_{i,j}(D)$ denote the set of enhanced 
Kauffman states $S$ such that $\sigma(S) = i$ and $\sigma(S) +2\tau(S) = j$, 
The group ${\cal C}(D)$ 
(resp. ${\cal C}_{i,j}(D)$) is defined to be the free abelian group 
spanned by ${\cal S}(D)$ (resp. ${\cal S}_{i,j}(D)$). ${\cal C}(D) = 
\bigoplus_{i,j\in \Z} {\cal C}_{i,j}(D)$ is a free abelian 
 group with (bi)-gradation.  
\item[(ii]) For a link diagram $D$ with ordered crossings, we define the 
chain complex $({\cal C}(D),d)$ where $d=\{d_{i,j}\}$ and 
the differential $d_{i,j}: {\cal C}_{i,j}(D) \to 
{\cal C}_{i-2,j}(D)$ satisfies $d(S) = \sum_{S'} (-1)^{t(S:S')}[S:S'] S'$ 
with $S\in {\cal S}_{i,j}(D)$, $S'\in {\cal S}_{i-2,j}(D)$, and 
 $[S:S']$  
equal to $0$ or $1$. $[S:S']=1$ if and only if markers of $S$ and $S'$ 
differ exactly at one crossing, call it $v$, and all the circles 
of $D_S$ and $D_{S'}$ 
not touching $v$ have the same sign\footnote{From our conditions 
it follows that at the crossing $v$ the marker of $S$ is positive, 
 the marker of $S'$ is negative, and 
that $\tau(S') = \tau(S)+1$.}. Furthermore, $t(S:S')$ is the number of 
negative markers assigned to crossings in $S$ bigger than $v$ in 
the chosen ordering. First two rows of Table 8.1 show all possible 
types of pairs of enhanced states for which $[S:S']=1$. $S$ and $S'$ have 
different markers at the crossing $v$ and $\varepsilon=+$ or $-$.

\item[(iii)] The Khovanov homology of the diagram $D$ is defined to be 
the homology of the chain complex $({\cal C}(D),d)$; 
$H_{i,j}(D) = ker(d_{i,j})/d_{i+2,j}({\cal C}_{i+2,j}(D))$. The 
Khovanov cohomology of the diagram $D$ is defined to be the cohomology 
of the chain complex $({\cal C}(D),d)$.
\end{enumerate}
\end{definition}

Khovanov proved that his homology is a topological invariant. 
$H_{i,j}(D)$ is preserved by the second and the third Reidemeister 
moves and $H_{i+1,j+3}(r_{+1}(D)) = H_{i,j}(D)=H_{i-1,j-3}(r_{-1}(D))$.
Where $r_{+1}(\parbox{0.8cm}{\psfig{figure=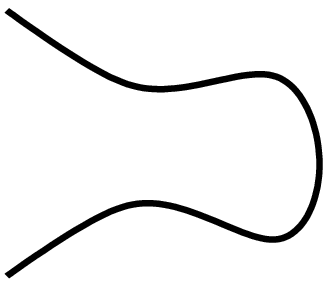,height=0.6cm}})=$
$(\parbox{1.0cm}{\psfig{figure=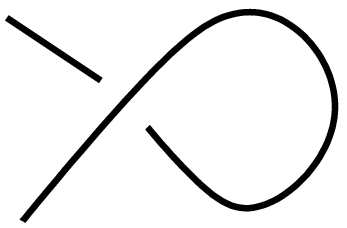,height=0.6cm}})$ and 
$r_{-1}(\parbox{0.8cm}{\psfig{figure=kinksmooth.eps,height=0.6cm}})=$
$(\parbox{1.0cm}{\psfig{figure=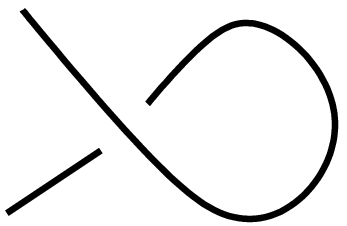,height=0.6cm}})$.

  With the notation we have introduced before, we can write 
the formula for the Kauffman bracket polynomial of a link 
diagram in the form\ 
$[D] = \sum_j A^j (\sum_i (-1)^{\frac{j-i}{2}}\sum_{S\in S_{i,j}} 1)=$
$\sum_j A^j (\sum_i (-1)^{\frac{j-i}{2}}dim C_{i,j})=$
$\sum_j A^j \chi_{i,j}(C_{*,j})$, where $\chi_{i,j}(C_{*,j}) = 
\sum_{i :\ j- i\equiv 0\ mod(2)}(-1)^{\frac{j-i}{2}}(dim C_{i,j})$ 
is (slightly adjusted) Euler characteristic of the chain 
complex $C_{*,j}$ ($j$ fixed).
This explains the phrase associated with Khovanov homology that it 
categorifies Kauffman bracket polynomial (or Jones polynomial)\footnote{In 
the narrow sense a categorification of a numerical or polynomial invariant 
is a homology theory whose Euler characteristic or polynomial Euler 
characteristic (the generating function of Euler characteristics) is 
the given invariant. We can quote after M.Khovanov \cite{Kh-1}:
 ``A speculative question now comes to mind: quantum invariants of knots 
and 3-manifolds tend to have good integrality properties. What if 
these invariants can be interpreted as Euler characteristics of some 
homology theories of 3-manifolds?".}.   

Below we list a few elementary properties of Khovanov homology
following from properties of Kauffman states used in the proof
of Tait conjectures \cite{K-5,M-4,This-3}; compare Chapter V.

The positive state $s_+=s_+D$ (respectively the negative
state $s_-=s_-D$) is the state with all positive markers
(resp. negative markers). The alternating diagrams without nugatory 
crossings (i.e. crossings in a diagram of the form
\parbox{2.2cm}{\psfig{figure=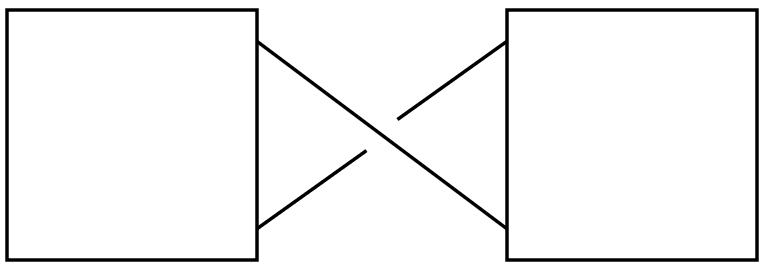,height=0.7cm}})
are generalized 
to adequate diagrams using properties of states $s_+$ and $s_-$. 
Namely, the diagram $D$ is $+$-adequate (resp. $-$-adequate) 
 if the state of positive (resp. negative) markers,
$s_+$ (resp. $s_-$), cuts the diagram to
the collection of circles, so that every crossing is
connecting different circles. $D$ is an adequate diagram if it is 
$+$- and $-$-adequate \cite{L-T}.

\begin{property}\label{1.4}\ \\
If $D$ is a diagram of $n$ crossings and its positive state $s_+$ has
$|s_+|$ circles then the highest term (in both grading indexes) 
 of Khovanov chain complex is
${\cal C}_{n,n+2|s_+|}(D)$; we have ${\cal C}_{n,n+2|s_+|}(D)= \Z$.
 Furthermore, if $D$ is a $+$-adequate diagram, 
 then the whole group
${\cal C}_{*,n+2|s_+|}(D)=\Z$ and $H_{n,n+2|s_+|}(D)= \Z$. 
Similarly the lowest term in the Khovanov
chain complex is ${\cal C}_{-n,-n-2|s_-|}(D)$. 
Furthermore, if $D$ is a $-$-adequate diagram,
 then the whole group
${\cal C}_{*,-n-2|s_-|}(D)$\\
$=\Z$ and $H_{-n,-n-2|s_-|}(D)= \Z$.
Assume that $D$ is
a non-split diagram then $|s_+| +|s_-| \leq n+2$ and the equality
holds if and only if $D$ is an alternating diagram or
a connected sum of such diagrams (Wu's dual state lemma \cite{Wu}; 
see Lemma V.3.13).
\end{property}
\begin{property} \label{1.5} \
Let $\sigma(L)$ be the classical (Trotter-Murasugi) 
signature\footnote{One should not mix the signature $\sigma(L)$ with
$\sigma(s)$ which is the signed sum of markers of the state $s$ of 
a link diagram.} of 
an oriented link $L$ and $\hat \sigma (L) = \sigma({L}) + lk({L})$, 
where $lk(L)$ is the global linking number of $L$, its Murasugi's 
version which does not depend on an orientation of $L$. Then  
(see Section V.4 for the proof of Traczyk results):
\begin{enumerate}
\item[(i)] [Traczyk's local property]\ 
If $D^v_0$ is a link diagram obtained from
an oriented alternating link diagram $D$ by smoothing its crossing 
$v$ and $D^v_0$ has the same number of (graph) components as $D$, then
$\sigma(D) = \sigma(D^v_0)- sgn(v)$. 
One defines the sign of a crossing $v$ as\
  $sgn(v) = \pm 1$ according to
the convention $sgn(\parbox{0.9cm}{\psfig{figure=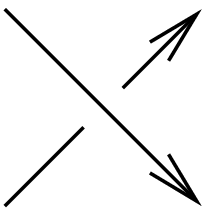,height=0.8cm}}
)= 1$ and $sgn(\parbox{0.9cm}{\psfig{figure=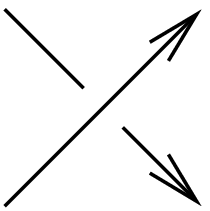,height=0.8cm}}
)=-1$.
\item[(ii)] [Traczyk Theorem \cite{T-2,AP}] \ \\
The signature, $\sigma (D)$, of the non-split alternating oriented 
link diagram $D$ is equal to 
$n^- - |s_-| +1 = -n^+ + |s_+|-1 = -\frac{1}{2}(n^+ - n^- 
-(|s_+|- |s_-|))=$ $n^- - n^+ + d^+ - d^-$, where $n^+(D)$ (resp. $n^-(D)$) 
is the number of positive (resp. negative) crossings
of $D$ and $d^+$ (resp. $d^-$) is the 
number of positive (resp.  negative) edges in a spanning forest of 
the Seifert graph\footnote{The Seifert graph, $GS(D)$, of 
an oriented link diagram $D$ is 
a signed graph whose vertices are in bijection with Seifert circles of $D$ 
and edges are in a natural bijection with crossings of $D$. For an alternating 
diagram the 2-connected components (blocks) of $GS(D)$ have edges of the 
same sign which makes $d^+$ and $d^-$ well defined.}  
of $D$. 
\item[(iii)]
[Murasugi's Theorem \cite{M-5,M-7}]\ \\
Let $\vec{D}$ be a non-split alternating oriented diagram without
nugatory crossings or a
connected sum of such diagrams. 
Let $V_{\vec{D}}(t)$ be its Jones polynomial\footnote{Recall that if 
$\vec{D}$ is an oriented diagram (any orientation
put on the unoriented diagram $D$),  and $w(\vec{D})$ is
its writhe or Tait number, $w(\vec{D}) = n^+-n^-$, 
then $V_{\vec{D}}(t) = A^{-3w(\vec{D})}<D>$ for $t=A^{-4}$.
P.G.Tait (1831-1901) was the first to consider the number $w(\vec{D})$
 and it is often called
the Tait number of the diagram $\vec{D}$ and denoted by
$Tait(\vec{D})$.},
 then the maximal degree $\max \ V_{\vec{D}}(t) = 
n^+(L) - \frac{\sigma(\vec{D})}{2}$
and the minimal degree $\min \ V_{\vec{D}}(t) = 
-n^-(\vec{D}) - \frac{\sigma(\vec{D})}{2}$.
\item[(iv)] [Murasugi's Theorem for unoriented link diagrams].
Let $D$ be a non-split alternating unoriented diagram without
nugatory crossings or a connected sum of such diagrams. \\
Then the maximal degree 
$\max  <D> = \max  [D] -2 = n+2|s_+|-2=2n+ sw(D)+2\hat\sigma(D)$\\
and the minimal degree 
$\min  <D> = \min  [D] +2 = -n - 2|s_-|+2= -2n+sw(D)+ 2\hat\sigma(D)$. 
The self-twist number of a diagram $sw(D)=
\sum_v sgn(v)$, where the sum is
taken over all self-crossings of $D$. A self-crossing involves arcs
from the same component of a link. $sw(D)$ does not depend on orientation 
of $D$.
\end{enumerate}
\end{property}
\begin{remark}\label{1.6} \ \\ 
In Section 5 we reprove the result of Lee \cite{Lee-1} 
that the Khovanov homology 
of non-split alternating links is supported on two adjacent 
diagonals of slope 2, that is $H_{i,j}(D)$ can be nontrivial only for 
two values of $j-2i$ which differ by $4$ (Corollary 5.5). 
One can combine Murasugi-Traczyk result with Viro's 
long exact sequence of Khovanov  
homology  and Theorem 7.3 to recover Lee's result (\cite{Lee-2}) 
that for alternating links
Khovanov homology has the same information as the Jones polynomial and
the classical signature\footnote{The beautiful
paper by  Jacob Rasmussen \cite{Ras} generalizes Lee's results  
and fulfills our dream
(with Pawe{\l} Traczyk) of constructing a ``supersignature" from Jones type
construction. Rasmussen signature $R(D)$ agrees with Trotter-Murasugi 
classical signature for alternating knots, satisfies Murasugi 
inequalities, $R(D_+) \leq R(D_-) \leq R(D_+)+2$ and allows to 
approximate the unknotting number, $u(K)\geq \frac{1}{2}R(K)$. 
Furthermore, it allows to find the unknotting number 
for positive knots, solving, in particular, 
Milnor's conjecture (Chapter V). For a positive diagram of a knot we 
have $R(D)=n(D) - Seif(D) +1$, where $Seif(D)$ is the number 
of Seifert circles of $D$.}. 
From properties 1.4 and 1.5 it follows that 
for non-split alternating diagram without nugatory crossings
$H_{n,2n+sw+2\hat\sigma +2}(D)=H_{-n,-2n+sw+2\hat\sigma -2}(D)= \Z$. 
Thus diagonals which support nontrivial $H_{i,j}(D)$ satisfy $j-2i= 
sw(D) +2\hat\sigma(D) \pm 2$. If we consider Khovanov cohomology 
$H^{i',j'}(D)$, as considered in \cite{Kh-1,Ba-2}, 
then $H^{i',j'}(D)=H_{i,j}(D)$ 
for $i'=\frac{w(D)-i}{2}$, $j'=\frac{3w(D)-j}{2}$ and thus $j'-2i' = 
\frac{-1}{2}(j-2i-w(D))=-\sigma(D)\mp 1$ as in Lee's Theorem. 
Note also that the self-linking number $sw(D)$ can be interpreted 
as the total framing number of an unoriented diagram $D$ (or the link 
defined by $D$ with blackboard framing). For a framed unoriented 
link we can define the (framed) signature $\sigma^f(D) = 
\hat\sigma(D) + \frac{1}{2}sw(D)= \sigma(\vec{D}) + \frac{1}{2}w(\vec{D})$. 
In this notation we get $j-2i = 2\sigma^f(D) \pm 2$, which is the framed 
version of Lee's theorem.
\end{remark}
\begin{exercise}\label{X.1.7}
 Compute the Khovanov homology for the diagram of the 
left handed trefoil knot presented in Fig. 1.2. Check in particular 
that the differential $d: C_{3,5}=\Z^3 \to Z^3=C_{1,5}$ is given 
by a matrix of determinant 2 (compare the first part of the proof of 
Theorem 2.2). 
Thus $H_{1,5}= Z_2$. Table 1.3 lists graded chain groups 
and graded homology of our diagram. Related calculation for reduced and 
co-reduced Khovanov homology is presented in detail in Example 7.6. 
\end{exercise}
\ \\
\centerline{\psfig{figure=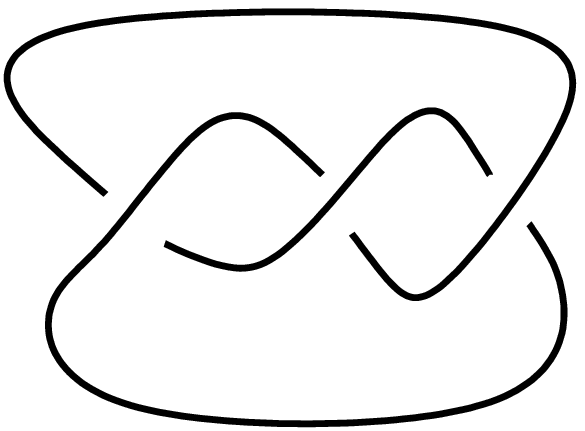,height=2.8cm}}
\begin{center}
Fig. 1.2
\end{center}
\ \\
\begin{center}
\begin{minipage}{11.5cm}
\[
\begin{array}{r||c|c|c|c|c|c|c|}
\cline {1-8}
9 &0 & 0 & 0 & 0 & 0 & 0 & \Z \\
\cline {1-8}
8 & 0 & 0 & 0 & 0 & 0 & 0 & 0 \\
\cline {1-8}
7&0 & 0 & 0 & 0 & 0 & 0 & 0 \\
\cline {1-8}
6&0 & 0 & 0 & 0 & 0 & 0 & 0 \\
\cline {1-8}
5&0 & 0 & 0 & 0 & \Z^3 & 0 & \Z^3 \\
\cline {1-8}
4&0 & 0 & 0 & 0 & 0 & 0 & 0 \\
\cline {1-8}
3& 0 & 0 & 0 & 0 & 0 & 0 & 0 \\
\cline {1-8}
2&0 & 0 & 0 & 0 & 0 & 0 & 0 \\
\cline {1-8}
1& \Z & 0 & \Z^3 & 0 & \Z^6 & 0 & \Z^3 \\
\cline {1-8}
0&0 & 0 & 0 & 0 & 0 & 0 & 0 \\
\cline {1-8}
-1&0 & 0 & 0 & 0 & 0 & 0 & 0 \\
\cline {1-8}
-2&0 & 0 & 0 & 0 & 0 & 0 & 0 \\
\cline {1-8}
-3& \Z^2 & 0 & \Z^3 & 0 & \Z^3 & 0 & \Z \\
\cline {1-8}
-4 &0 & 0 & 0 & 0 & 0 & 0 & 0 \\
\cline {1-8}
-5 & 0 & 0 & 0 & 0 & 0 & 0 & 0\\
\cline {1-8}
 -6&0 & 0 & 0 & 0 & 0 & 0 & 0 \\
\cline {1-8}
-7 &\Z & 0 & 0 & 0 & 0 & 0 & 0  \\
\cline {1-8}
\cline {1-8}
C_{i,j}&-3 & -2& -1& 0 & 1 & 2 & 3
\end{array} \ \ \ \ \ \ \ \ \
\begin{array}{r||c|c|c|c|c|c|c|}
\cline {1-8}
9 &0 & 0 & 0 & 0 & 0 & 0 & \Z \\
\cline {1-8}
8 & 0 & 0 & 0 & 0 & 0 & 0 & 0 \\
\cline {1-8}
7&0 & 0 & 0 & 0 & 0 & 0 & 0 \\
\cline {1-8}
6&0 & 0 & 0 & 0 & 0 & 0 & 0 \\
\cline {1-8}
5&0 & 0 & 0 & 0 & \Z_2 & 0 & 0 \\
\cline {1-8}
4&0 & 0 & 0 & 0 & 0 & 0 & 0 \\
\cline {1-8}
3&0 & 0 & 0 & 0 & 0 & 0 & 0 \\
\cline {1-8}
2&0 & 0 & 0 & 0 & 0 & 0 & 0 \\
\cline {1-8}
1&0 & 0 & 0 & 0 & \Z & 0 & 0 \\
\cline {1-8}
0&0 & 0 & 0 & 0 & 0 & 0 & 0 \\
\cline {1-8}
-1&0 & 0 & 0 & 0 & 0 & 0 & 0 \\
\cline {1-8}
-2&0 & 0 & 0 & 0 & 0 & 0 & 0 \\
\cline {1-8}
-3& \Z & 0 & 0 & 0 & 0 & 0 & 0 \\
\cline {1-8}
-4 &0 & 0 & 0 & 0 & 0 & 0 & 0 \\
\cline {1-8}
-5 & 0 & 0 & 0 & 0 & 0 & 0 & 0\\
\cline {1-8}
 -6&0 & 0 & 0 & 0 & 0 & 0 & 0 \\
\cline {1-8}
-7 &\Z & 0 & 0 & 0 & 0 & 0 & 0  \\
\cline {1-8}
\cline {1-8}
H_{i,j}&-3 & -2& -1& 0 & 1 & 2 &  3
\end{array}
\]
\end{minipage}
\end{center}
\begin{center}
Table 1.3
\end{center}

\section{Diagrams with odd cycle property}\label{2}
In the next few sections we use the concept of a graph, $G_s(D)$,
 associated to a link diagram $D$ and its state $s$. The graphs 
corresponding to states $s_+$ and $s_-$ are of particular interest.
If $D$ is an alternating diagram then $G_{s_+}(D)$ and $G_{s_-}(D)$ 
are the plane graphs first constructed by Tait.

\begin{definition}\label{2.1}\ \
\begin{enumerate}
\item[(i)]
Let $D$ be a diagram of a link and $s$ its Kauffman state. 
We form a graph, $G_s(D)$, associated to $D$ and $s$ as follows.
 Vertices of $G_s(D)$ correspond to circles of $D_s$.
Edges of $G_s(D)$ are in bijection with crossings of $D$ and
an edge connects given 
vertices if the corresponding crossing connects circles of $D_s$ 
corresponding to the vertices\footnote{If $S$ is an enhanced Kauffman 
state of $D$ then, in a similar manner, we associate to $D$ and $S$ the 
graph $G_S(D)$ with signed vertices. Furthermore, we can additionally 
equip $G_S(D)$ with a cyclic ordering of edges at every vertex 
following the ordering of crossings at any circle of $D_s$. 
The sign of each edge is the label of the corresponding crossing.
In short, 
we can assume that $G_S(D)$ is a ribbon (or framed) graph. 
We do not use this additional data in this chapter but it  
may be of great use in analysis of Khovanov homology.}.
\item[(ii)] 
In the language of associated graphs we can state the definition
of adequate diagrams as follows:\ 
the diagram $D$ is $+$-adequate (resp. $-$-adequate) if the graph 
$G_{s_+}(D)$ (resp. $G_{s_-}(D)$) has no loops. 
\end{enumerate}
\end{definition}
In this language we can formulate our first result about torsion in 
Khovanov homology. 
\begin{theorem}\label{2.2} \ \\
Consider a link diagram $D$ of $N$ crossings. Then 
\begin{enumerate}
\item[(+)]
If $D$ is $+$-adequate and $G_{s_+}(D)$ has a cycle 
of odd length, then  the Khovanov homology has $\Z_2$ torsion.
More precisely we show that \\ 
$H_{N-2,N+2|s_+|-4}(D)$ has $\Z_2$ torsion.
\item[(-)]
If $D$ is $-$-adequate and $G_{s_-}(D)$ has a cycle
of odd length, then \\ 
$H_{-N,-N-2|s_-|+4}(D)$ has $\Z_2$ torsion.
\end{enumerate}
\end{theorem}
\begin{proof} 
$(+)$ It suffices to show that the group\\ ${\cal C}_{N-2,N+2|s_+|-4}(D)/
d({\cal C}_{N,N+2|s_+|-4}(D))$ has $2$-torsion.

Consider first the diagram $D$ of the left handed torus knot $T_{-2,n}$ 
(Fig.2.1 illustrates the case of $n=5$). The associated graph 
$G_n=G_{s_+}(T_{-2,n})$ is an $n$-gon.
\\
\ \\
\centerline{\psfig{figure=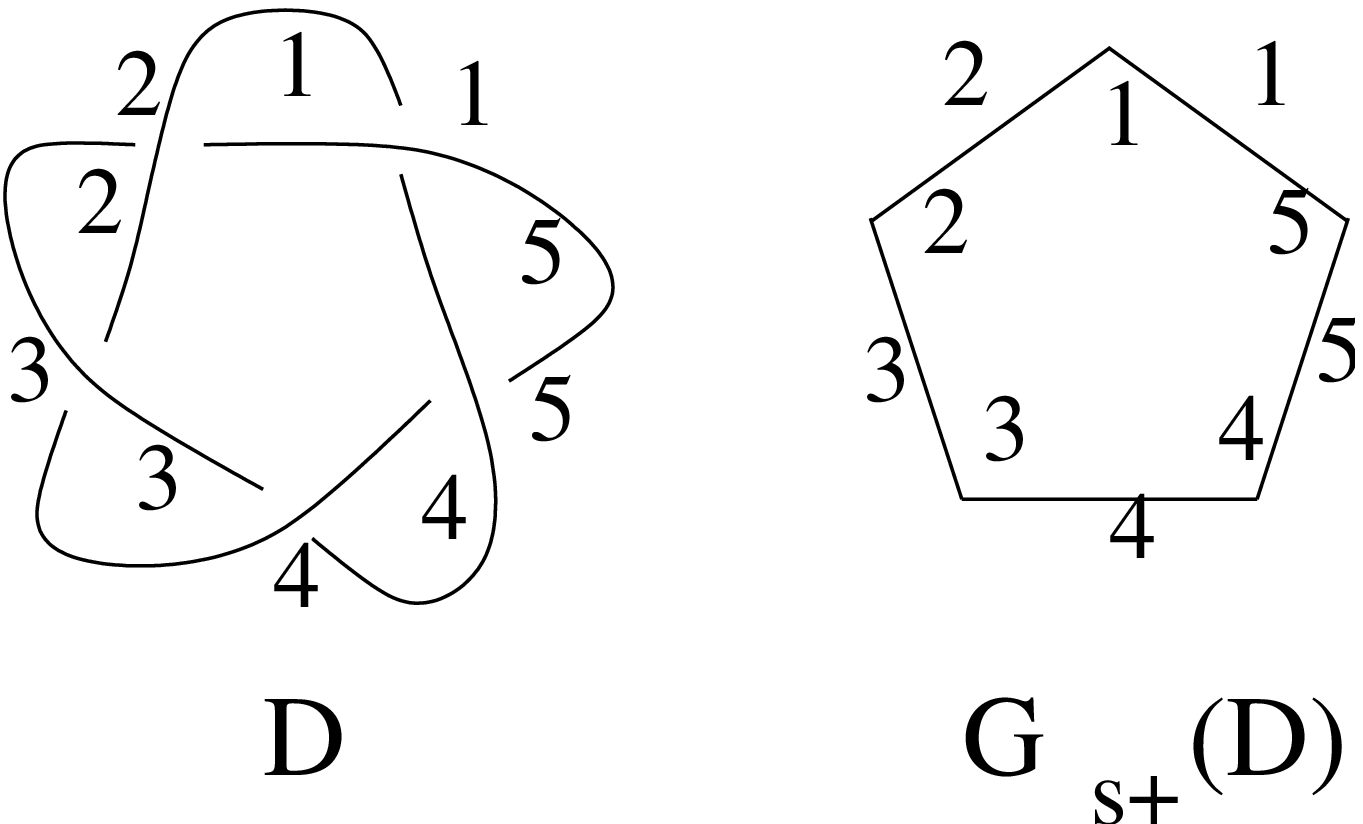,height=3.8cm}}
\begin{center}
Fig. 2.1
\end{center}

 For this diagram we have 
${\cal C}_{n,n+2|s_+|}(D)= \Z$, ${\cal C}_{n,n+2|s_+|-4}(D)=\Z^n$ and 
${\cal C}_{n-2,n+2|s_+|-4}(D) = \Z^n$, where enhanced states generating 
${\cal C}_{n,n+2|s_+|-4}(D)$ have all markers positive and 
exactly one circle (of $D_S$) negative\footnote{In this case $s_+=n$ but we 
keep the general notation so the generalization which follows is natural.}. 
Enhanced states generating
 ${\cal C}_{n-2,n+2|s_+|-4}(D)$ have exactly one negative marker and 
all positive circles of $D_S$. The differential 
$$d: {\cal C}_{n,n+2|s_+|-4}(D) \to  {\cal C}_{n-2,n+2|s_+|-4}(D)$$ can be 
described by an $n \times n$ circulant matrix (for the ordering of states 
corresponding to the ordering of crossings and regions as in Fig. 2.1)). 

$$\left( \begin{array}{cccccc} 
1 & 1 &  0 & \ldots & 0 & 0\\
0 & 1 & 1 &\ldots & 0 & 0 \\
\ldots & \ldots & \ldots & \ldots & \ldots & \ldots\\ 
0 & 0& \ldots &0 & 1 & 1 \\
1 & 0 & \ldots & 0 & 0 & 1
\end{array} \right)$$

Clearly the determinant of the matrix is equal to 2 (because $n$ is 
odd; for $n$ even the determinant is equal to $0$ because the alternating 
sum of columns gives the zero column). To see this one can 
consider for example the first row expansion\footnote{Because 
the matrix is a circulant one we know furthermore that its eigenvalues are 
equal to $1 + \omega$, where $\omega$ is any $n$'th root of unity 
($\omega^n = 1$), and that $\prod_{\omega^n = 1}(1+\omega)= 0$ for $n$ even 
and $2$ for $n$ odd.}. Therefore the group described by the matrix
is equal to $\Z_2$ (for an even $n$ one would get $\Z$). One more observation 
(which will be used later). The sum of rows of the matrix is equal to 
the row vector ($2,2,2,...,2,2$) but the row vector ($1,1,1,...,1,1$) 
is not an integral linear combination of rows of the matrix. 
In fact the element 
$(1,1,1,...,1,1)$ is the generator of $\Z_2$ group represented by the 
matrix. This can be easily checked because if $S_1,S_2,....S_n$ are 
states freely generating ${\cal C}_{n-2,n+2|s_+|-4}(D)$ then 
relations given by the image of ${\cal C}_{n,n+2|s_+|-4}(D)$ are
$S_2=-S_1, S_3=-S_2=S_1,..., S_1=-S_n=...=-S_1$ thus $S_1+S_2+...+S_n$ is 
the generator of the quotient group ${\cal C}_{n-2,n+2|s_+|-4}(D)/
d({\cal C}_{n,n+2|s_+|-4}(D))= \Z_2$. In fact we have proved that any 
sum of the odd number of states $S_i$ represents the generator of $\Z_2$. 

Now consider the general case in which $G_{s_+}(D)$ is a graph without 
a loop and with an odd polygon. Again, we build a matrix presenting 
the group 
${\cal C}_{N-2,N+2|s_+|-4}(D)/d({\cal C}_{N,N+2|s_+|-4}(D))$ with the
north-west block corresponding to the odd $n$-gon. This block is exactly 
the matrix described previously. Furthermore, the submatrix of the full matrix
 below this block is the zero matrix, as every column 
has exactly two nonzero entries (both equal 
to $1$). This is the case because each edge of the graph (generator) 
has two endpoints (belongs to exactly two relations). If we add  
all rows of the matrix 
 we get the row of all two's. On the other hand 
the row of one's cannot be created, even in the first block. Thus the row 
of all one's representing the sum of all enhanced states in 
${\cal C}_{N-2,N+2|s_+|-4}(D)$ is  $\Z_2$-torsion element in 
the quotient group (presented by the matrix) 
so also in $H_{N-2,N+2|s_+|-4}(D)$. \\
(-) This part follows from the fact that the mirror image of $D$, the diagram 
$\bar D$, satisfies the assumptions of the part (+) of the theorem. 
Therefore the quotient 
${\cal C}_{N-2,N+2|s_+|-4}(\bar D)/d({\cal C}_{N,N+2|s_+|-4}(\bar D))$ has 
$\Z_2$ torsion. Furthermore, the matrix describing the map 
$d: {\cal C}_{-N+2,-N-2|s_-|+4}(D) \to {\cal C}_{-N,-N-2|s_-|+4}(D)$ is
(up to sign of every row) equal to the transpose 
of the matrix describing the map 
$d: {\cal C}_{N,N+2|s_+|-4}(\bar D) \to {\cal C}_{N-2,N+2|s_+|-4}(\bar D)$.
Therefore the torsion of the group 
${\cal C}_{-N,-N-2|s_-|+4}(D)/d({\cal C}_{-N+2,-N-2|s_-|+4}(D)$ is the same 
as the torsion of the group 
${\cal C}_{N-2,N+2|s_+|-4}(\bar D)/d({\cal C}_{N,N+2|s_+|-4}(\bar D))$ 
and, in conclusion, $H_{-N,-N-2|s_-|+4}(D)$ has  $\Z_2$ 
torsion\footnote{Our reasoning reflects a more general fact observed  
 by Khovanov \cite{Kh-1} (see \cite{APS-2} for the case of $F \times I$)
 that Khovanov homology satisfies ``duality theorem", namely 
$H^{i,j}(D) = H_{-i,-j}({\bar D})$. This combined with the 
Universal Coefficients Theorem saying that 
$H^{i,j}(D) = H_{i,j}(D)/T_{i,j}(D) \oplus T_{i-2,j}(D)$, where 
$T_{i,j}(D)$ denote the torsion part of $H_{i,j}(D)$ gives:
$T_{-N,-N-2|s_-|+4}(D) = T_{N-2,N+2|s_+|-4}(\bar D)$ (notice that 
$|s_-|$ for $D$ equals to $|s_+|$ for $\bar D$).}.

\end{proof}

\begin{remark}\label{2.3}\ \  \\ 
Notice that the torsion part of the homology, $T_{N-2,N+ 2|s_+|-4}(D)$,
depends only on the graph $G_{s_+}(D)$.
Furthermore if $G_{s_+}(D)$ has no 2-gons then $H_{N-2,N+ 2|s_+|-4}(D)=
{\cal C}_{N-2,N+ 2|s_+|-4}(D)/d(C_{N,N+ 2|s_+|-4}(D))$ and depends only 
on the graph $G_{s_+}(D)$. See a generalization in Remark 3.6
\end{remark}
\section{Diagrams with an even cycle property}\label{3}

If every cycle of the graph $G_{s_+}(D)$ is even (i.e. the graph is 
a bipartite graph) we cannot expect that $H_{N-2,N+|s_+|-4}(D)$  
always has nontrivial torsion. The simplest link diagram 
without an odd cycle in
$G_{s_+}(D)$ is the left handed torus 
link diagram $T_{-2,n}$ for $n$ even. As mentioned before, 
in this case ${\cal C}_{n-2,n+2|s_+|-4}(D)/
d({\cal C}_{n,n+2|s_+|-4}(D))= \Z$, and, 
in fact $H_{n-2,n+2|s_+|-4}(D)=\Z$ except for $n=2$, i.e. the Hopf link, in 
which case $H_{0,2}(D)=0$.

To find torsion we have to look ``deeper" into the homology. We 
give a condition on a diagram $D$ with $N$ crossings which guarantees
that $H_{N-4,N+2|s_+|-8}(D)$  has  $\Z_2$ torsion, 
where $N$ is the number of crossings of $D$.

Analogously to the odd case, we will start from the left handed torus
link $T_{-2,n}$ and associated graph $G_{s_+}(D)$ being an $n$-gon with 
even $n \geq 4$; Fig.3.1. 
\\
\ \\
\centerline{\psfig{figure=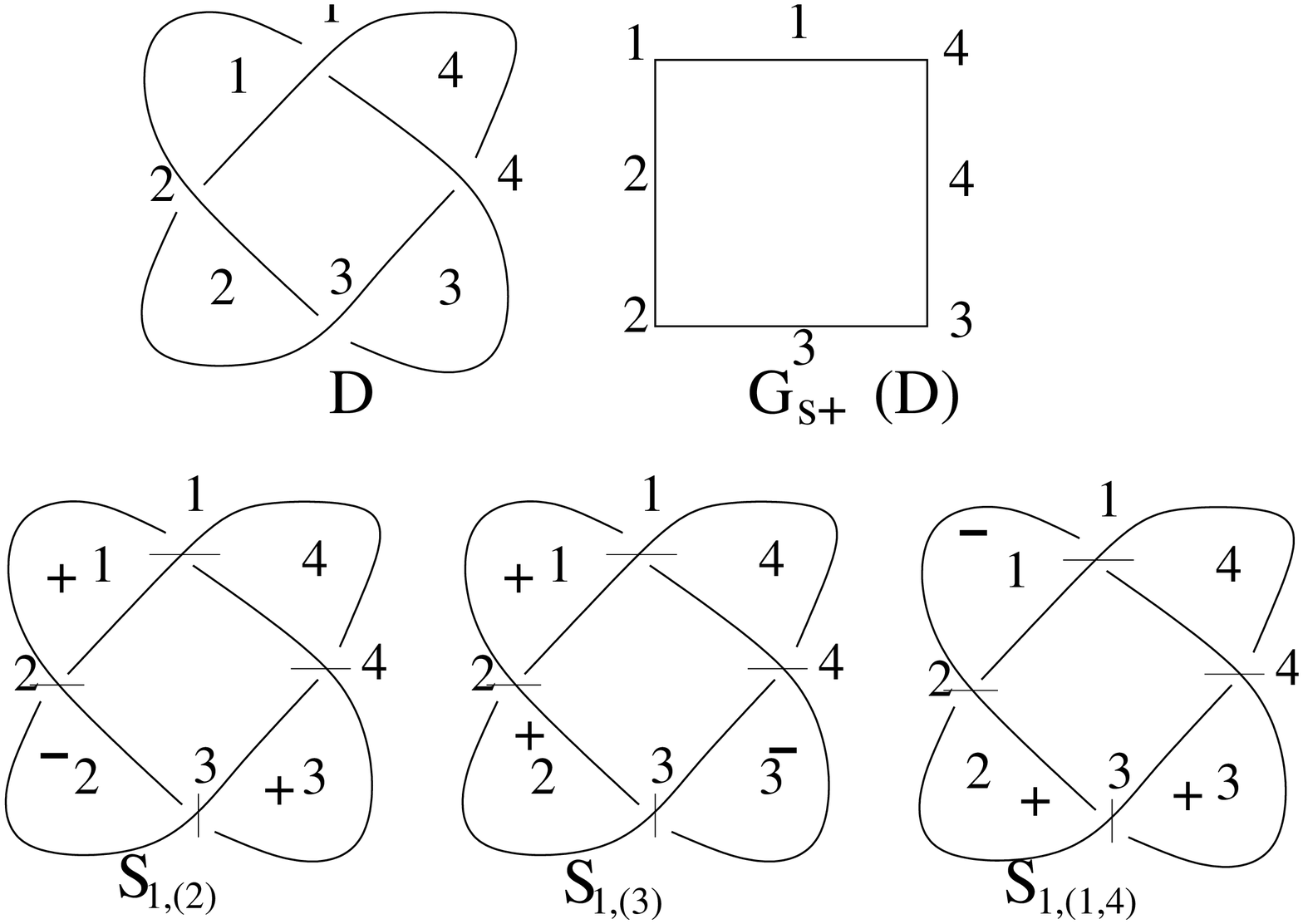,height=5.8cm}}
\begin{center}
Fig. 3.1
\end{center}

\begin{lemma}\label{3.1}\ \\
Let $D$ be the diagram of the left-handed torus link of type 
($-2,n$) with $n$ even, $n \geq 4$.\\
Then $H_{n-4,n+2|s_+|-8}(D)= H_{n-4,3n-8}(D) = 
{\cal C}_{n-4,3n-8}(D)/d({\cal C}_{n-2,3n-8}(D)) = \Z_2$.\\ 
Furthermore, every enhanced state from the basis of 
${\cal C}_{n-4,3n-8}(D)$ (or an odd sum of such states) is 
the generator of $\Z_2$.
\end{lemma}
\begin{proof}
We have  $n=|s_+|$. 
The chain group $C_{n-4, 3n-8}(D)= \Z^{\frac{n(n-1)}{2}}$
is freely generated by enhanced states $S_{i,j}$, where exactly 
$i$th and $j$th 
crossings have negative markers, and all the circles of $D_{S_{i,j}}$ are 
positive (crossings of $D$ and circles of $D_{s_+}$ are ordered in Fig. 3.1). 
We have to understand the differential
$d: C_{n-2, 3n-8}(D) \to C_{n-4, 3n-8}(D)$. The chain group 
$C_{n-2, 3n-8}(D) = \Z^{n(n-1)}$ is freely generated by enhanced 
states with $i$th negative marker and one negative circle of $D_S$. 
In our notation we will write $S_{i,(i-1,i)}$ if the negative circle is
obtained by connecting circles $i-1$ and $i$ in $D_{s_+}$ by a negative 
marker. Notation $S_{i,(j)}$ is used if we have $j$th negative circle, 
$j\neq i-1$, $j \neq i$. The states $S_{1,(2)}$, $S_{1,(3)}$ and 
$S_{1,(4,1)}$ are 
shown in Fig. 3.1 ($n=4$ in the figure). The quotient group 
$  {\cal C}_{n-4, 3n-8}(D)/d({\cal C}_{n-2, 3n-8}(D))$ can be presented 
by  a $n(n-1) \times \frac{n(n-1)}{2}$ matrix, $E_n$. One should just 
understand the images of enhanced states of $  {\cal C}_{n-2, 3n-8}(D)$.
In fact, for a fixed crossing $i$ the corresponding $n-1 \times n-1$ 
block is (up to sign of columns\footnote{In the $(n-1)\times (n-1)$ 
block corresponding to the $i$th crossing (i.e. we consider only states 
in which $i$th crossing has a negative marker), the column under the 
generator $S_{i,j}$ of $ {\cal C}_{n-4, 3n-8}$ has $+1$ entries 
if $i<j$ and $-1$ entries if $i>j$.}) 
the circulant matrix discussed in Section 2. 
Our goal is to understand  the matrix $E_n$, to show that 
it represents the group $\Z_2$ and to find natural representatives of 
the generator of the group. 
For $n=4$, $d: \Z^{12} \to \Z^6$ and it is given by: 
$d(S_{1,(2)}) = S_{1,2} + S_{1,3}$, $d(S_{1,(3)}) = S_{1,3} + S_{1,4}$, 
$d(S_{1,(1,4)}) = S_{1,2} + S_{1,4}$, $d(S_{2,(1,2)}) = -S_{2,1} + S_{2,3}$,
 $d(S_{2,(3)}) = S_{2,3} + S_{2,4}$, $d(S_{2,(4)}) = -S_{2,1} + S_{2,4}$, 
$d(S_{3,(1)}) = -S_{3,1} - S_{3,2}$, $d(S_{3,(2,3)}) = -S_{3,2} + S_{3,4}$,
$d(S_{3,(4)}) = -S_{3,1} + S_{3,4}$, $d(S_{4,(1)}) = -S_{4,1} - S_{4,2}$, 
$d(S_{4,(2)}) = -S_{4,2} - S_{4,3}$, $d(S_{4,(4,3)}) = -S_{4,1} - S_{4,3}$,

Therefore $d$ can be described by the $12 \times 6$ matrix. States are 
ordered lexicographically, e.g. $S_{i,j}$ ($i<j$) is before 
$S_{i',j'}$ ($i'<j'$) if $i<i'$ or $i=i'$ and $j<j'$.
$$\left( \begin{array}{cccccc}
1 & 1 & 0 & 0 & 0 & 0\\
0 & 1 & 1 & 0 & 0 & 0 \\
1 & 0 & 1 & 0 & 0 & 0 \\
-1& 0 & 0 & 1 & 0 & 0 \\
0 & 0 & 0 & 1 & 1 & 0 \\
-1& 0 & 0 & 0 & 1 & 0 \\
0 & -1& 0 & -1& 0 & 0 \\
0 & 0 & 0 & -1& 0 & 1 \\
0 & -1& 0 & 0 & 0 & 1 \\
0 & 0 & -1& 0 & -1& 0 \\
0 & 0 & 0 & 0 & -1& -1 \\
0 & 0 & -1& 0 & 0 & -1 
\end{array} \right),$$
In our example the rows correspond to 
$S_{1,(2)},S_{1,(3)},S_{1,(1,4)},S_{2,(1,2)},
S_{2,(3)}, S_{2,(4)}, S_{3,(1)},$ \\
$S_{3,(2,3)},S_{3,(4)},S_{4,(1)},S_{4,(2)},$
and $ S_{4,(4,3)}$, the columns correspond to\\  
$S_{1,2},S_{1,3},S_{1,4},S_{2,3},S_{2,4},S_{3,4}$ 
in this order. 
Notice that the sum of columns of the matrix gives the non-zero 
column of all $\pm 2$ or $0$.
Therefore over $\Z_2$ our matrix represents a nontrivial group. On the 
other hand, over $\Q$, the matrix represent the trivial group. Thus 
over $\Z$ the group represented by the matrix has $\Z_2$ torsion. 
More precisely, we can see that the group is $\Z_2$ as follows: 
The row relations can be expressed as: 
$S_{1,2}=-S_{1,3}=S_{1,4}=-S_{1,2}$, $S_{2,1}= S_{2,3}=-S_{2,4}= -S_{2,1}$, 
$S_{3,1} =- S_{3,2}= -S_{3,4}=-S_{3,1}$ and 
$S_{4,1}= -S_{4,2} = S_{4,3} =-S_{4,1}$. 
$S_{i,j}=S_{j,i}$ in our notation.
In particular, it follows from these 
equalities that the group 
given by the matrix is equal to $\Z_2$ and is generated by 
any basic enhanced state $S_{i,j}$ or the sum of odd number of $S_{i,j}$'s.  

Similar reasoning works for any even $n\geq 4$ (not only $n=4$). 

Furthermore, ${\cal C}_{n-6, 3n-8}=0$, therefore $H_{n-4,3n-8}= 
{\cal C}_{n-4, 3n-8}/d({\cal C}_{n-2, 3n-8}) = \Z_2$.
\end{proof}
We are ready now to use Lemma 3.1 in the general case of an even cycle.
\begin{theorem}\label{3.2}\ \\
Let $D$ be a connected diagram of a link of $N$ crossings such that 
the associated graph $G_{s_+}(D)$ has no loops (i.e. $D$ is $+$-adequate) 
and the graph has an even $n$-cycle with a singular edge 
(i.e. not a part of a 2-gon). 
Then $H_{N-4,N+2|s_+|-8}(D)$  has  $\Z_2$ torsion.
\end{theorem}
\begin{proof}
Consider an ordering of crossings of $D$ such that $e_1,e_2,...,e_n$ are
crossings (edges) of the $n$-cycle. The chain group 
${\cal C}_{N-2,N+2|s_+|-8}(D)$ is freely generated by 
 $N(V-1)$ enhanced
states,$ S_{i,(c)}$, where $N$ is the number of crossings of $D$ (edges
of $G_{s_+}(D)$) and $V=|s_+|$ is the number of circles 
of $D_{s_+}$ (vertices
of $G_{s_+}(D)$). $ S_{i,(c)}$ is the enhanced state in which the crossing
$e_i$ has the negative marker and the circle $c$ of $D_{s_i}$ is negative,
where $s_i$ is the state which has all positive markers except at $e_i$.
The chain group
${\cal C}_{N-4,N+2|s_+|-8}(D)$ is freely generated by enhanced states 
which we can partition into two groups.\\ 
(i) States $S_{i,j}$, where crossings 
$e_i$, $e_j$ have negative markers and corresponding edges of 
$G_{s_+}(D)$ do not form part of a multi-edge 
(i.e. $e_i$ and $e_j$ do not have the 
same endpoints). All circles of the state
$S_{i,j}$ are positive.\\ 
(ii) States $S'_{i,j}$ and $S''_{i,j}$, where crossings
$e_i$, $e_j$ have negative markers and corresponding edges of
$G_{s_+}(D)$ are parts of a multi-edge (i.e. $e_i$, $e_j$ have the same 
endpoints). All but one circle of $S'_{i,j}$ and $S''_{i,j}$ are positive 
and we have two choices for a negative circle leading to 
$S'_{i,j}$ and $S''_{i,j}$, i.e. the crossings $e_i$, $e_j$ touch two circles, 
and we give negative sign to one of them.\\
In our proof we will make the essential use of the assumption that the 
edge (crossing) $e_1$ is a singular edge.

We analyze the matrix presenting the group\\
 ${\cal C}_{N-4,N+2|s_+|-8}(D)/d({\cal C}_{N-2,N+2|s_+|-8}(D))$.

By Lemma 3.1, we understand already the $n(n-1) \times \frac{1}{2}n(n-1)$
block corresponding to the even $n$-cycle.  In this block every column has
$4$ non-zero entries (two $+1$ and two $-1$),
therefore columns of the full matrix corresponding
to states $S_{i,j}$, where $e_i$ and $e_j$ are in the $n$-gon, have zeros
outside our block. We use this property later.

We now analyze another block
represented by rows and columns associated to states having the first
crossing $e_1$ with the negative marker.
This  $(V-1)\times (N-1)$ block has entries equal to $0$ or $1$.
If we add rows in this block we obtain the vector row of two's ($2,2,...,2$),
following from the fact that every edge of $G_{s_+}(D)$  and of $G_{s_1}(D)$
has $2$ endpoints (we use the fact that $D$ is $+$ adequate and $e_1$ is a 
singular edge). Consider
now the bigger submatrix of the full matrix composed of the same rows as
our block but without restriction on columns.
All additional columns are $0$ columns as our row relations
involve only states with negative marker at $e_1$. Thus the sum of these rows
is equal to the row vector ($2,2,...,2,0,...,0$). We will argue now that
the half of this vector, ($1,1,...,1,0,...,0$), is not an integral linear
combination of rows of the full matrix and so represents  $\Z_2$-torsion
element of the group
${\cal C}_{N-4,N+2|s_+|-8}(D)/d({\cal C}_{N-2,N+2|s_+|-8}(D))$.
For simplicity assume that $n=4$ (but the argument holds for any
even $n \geq 4$). Consider the columns indexed by $S_{1,2},S_{1,3},S_{1,4},
S_{2,3},S_{2,4}$ and $S_{3,4}$. The integral linear combination of rows
restricted to this columns cannot give a row with odd number of one's, as
proven in Lemma 3.1. 
In particular we cannot get the row vector ($1,1,1,0,0,0$).
This excludes the row ($1,1,...,1,0,...,0$), as  an integral linear
combination of rows of the full matrix. Therefore the sum of enhanced states
with the marker of $e_1$ negative is  $2$-torsion element in
${\cal C}_{N-4,N+2|s_+|-8}(D)/d({\cal C}_{N-2,N+2|s_+|-8}(D))$ and therefore
in $H_{N-4,N+2|s_+|-8}(D)$.
\end{proof}
Similarly, using duality, we can deal with $-$-adequate diagrams.
\begin{corollary}\label{3.3}\ \\
Let $D$ be a connected, $-$-adequate diagram of a link and the
graph $G_{s_-}(D)$ has an even $n$-cycle, $n \geq 4$, with a singular edge.
Then $H_{-N+2,-N-2|s_-|+8}(D)$  has  $\Z_2$ torsion.
\end{corollary}

\begin{remark}\label{3.4}\ \\
The restriction on $D$ to be a connected diagram is not essential 
(it just simplifies the proof) as for a non-connected diagram, 
$D=D_1\sqcup D_2$ we have ``K\"unneth formula" 
$H_*(D)= H_*(D_1) \otimes H_*(D_2)$ so if any 
of $H_*(D_i)$ has  torsion then $H_*(D)$ has torsion as well.
\end{remark}

We say that a link diagram is doubly $+$-adequate if its graph 
$G_{s_+}(D)$ has no loops and 2-gons. In other words, if a state 
$s$ differs from the state $s_+$ by two markers then $|s| = |s_+| -2 $.
We say that a link diagram is doubly $-$-adequate if its mirror image 
is doubly $+$-adequate.

\begin{corollary}\label{3.5}\ \\
Let $D$ be a connected doubly $+$-adequate diagram of a link of 
$N$ crossings, then either $D$ represents the trivial knot 
or one of the groups
 $H_{N-2,N+2|s_+|-4}(D)$ and $H_{N-4,N+2|s_+|-8}(D)$ has  $\Z_2$ torsion. 
\end{corollary}
\begin{proof}
The associated graph $G_{s_+}(D)$ has no loops and 2-gons.
If $G_{s_+}(D)$ has an odd cycle then 
by Theorem 2.2 $H_{N-2,N+2|s_+|-4}(D)$ has $\Z_2$ torsion. If $G_{s_+}(D)$ 
has an even $n$-cycle, $n \geq 4$ then $H_{N-4,N+2|s_+|-8}(D)$ has 
$\Z_2$ torsion by Theorem 3.2 (every edge of $G_{s_+}(D)$ is a singular edge 
as $G_{s_+}(D)$ has no $2$-gons). Otherwise $G_{s_+}(D)$ is a tree, 
each crossing of $D$ is a nugatory crossing and $D$ 
represents the trivial knot.
\end{proof}

We can generalize and interpret Remark 2.3 as follows.
\begin{remark}\label{3.6} \ \\
 Assume that the associated graph $G_{s_+}(D)$ 
has no $k$-gons, for every $k\leq m$. Then
the torsion part of Khovanov homology, $T_{N-2m,N+ 2|s_+|-4m}(D)$
depends only on the graph $G_{s_+}(D)$.
Furthermore, $H_{N-2m+2,N+ 2|s_+|-4m+4}(D)=
{\cal C}_{N-2m+2,N+ 2|s_+|-4m+4}(D)/d(_{N-2m+4,N+ 2|s_+|-4m+4}(D))$ 
and it depends only on the graph $G_{s_+}(D)$. On a more philosophical 
level
\footnote{In order to be able to recover the full Khovanov homology 
from the graph $G_{s_+}$ we would have to equip the graph with additional 
data: \ ordering of signed edges adjacent to every vertex. This allows us 
to construct a closed surface and the link diagram $D$ on it so that 
$G_{s_+} =G_{s_+}(D)$. The construction imitates the $2$-cell embedding of 
Heffter-Edmonds (but every vertex corresponds to a circle 
and signs of edges regulate whether 
an edge is added inside or outside of the circle). 
If the surface we obtain is equal to $S^2$ we get the classical Khovanov 
homology. If we get a higher genus surface we have to use \cite{APS-2} 
theory. This can also be utilized  to construct  Khovanov homology of 
virtual links (via Kuperberg minimal genus embedding theory \cite{Ku}). 
For example, if the graph $G_{s_+}$ is a loop with adjacent edge(s) 
ordered $e,-e$ then the diagram is composed of a meridian 
and a longitude on the torus.} 
 our observation is related to the fact that if the edge $e_c$ in $G_{s_+}(D)$  
corresponding to a crossing $c$ in $D$ is not a loop then for 
the crossing $c$ the graphs $G_{s_+}(D_0)$ and 
$G_{s_+}(D_{\infty})$ are the graphs obtained from $G_{s_+}(D)$ by 
deleting ($G_{s_+}(D)-e_c$) and contracting ($G_{s_+}(D)/e_c$), 
respectively, the edge $e_c$ (compare Fig.3.2).
\end{remark}
\centerline{\psfig{figure=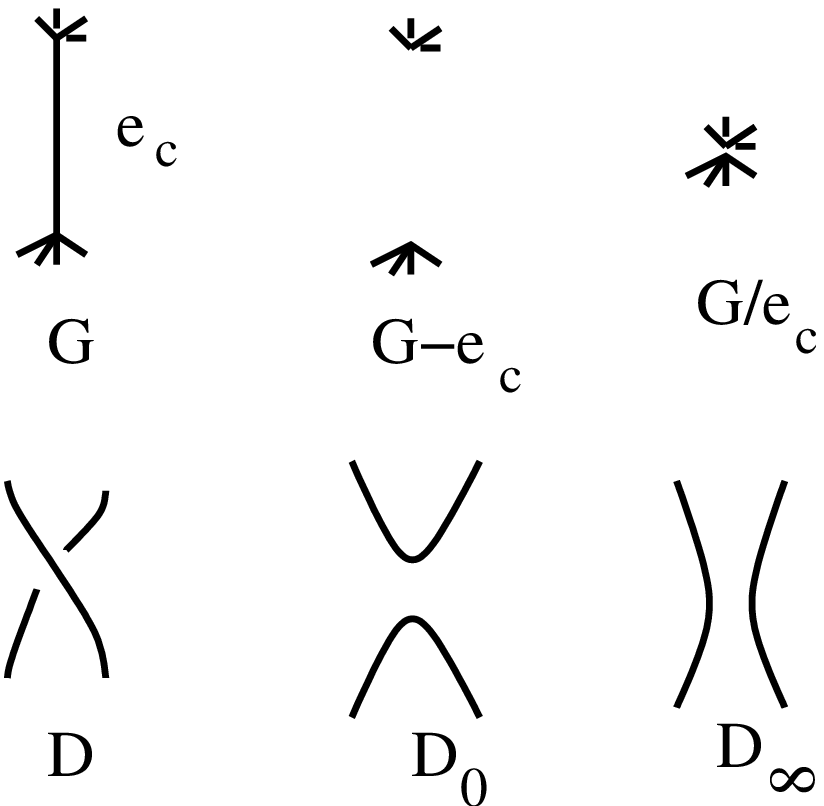,height=4.4cm}}
\begin{center}
Fig. 3.2
\end{center}

\begin{example}\label{3.7}\ \\
 Consider the 2-component alternating link
$6^2_2$  ($\frac{10}{3}$ rational link), with $G_{s_+}(D) = 
G_{s_-}(D)$ being a square with one edge tripled (this is 
a self-dual graph); see Fig 3.3. 
Corollary 3.5 does not apply to this case but Theorem 3.2
guarantees $\Z_2$ torsion at $H_{2,6}(D)$ and $H_{-4,-6}(D)$. \\
In fact, the KhoHo \cite{Sh-2} computation gives the following 
Khovanov homology\footnote{Tables and programs by Bar-Natan 
and Shumakovitch \cite{Ba-4,Sh-2} use 
the version of Khovanov homology for oriented diagrams, 
and the variable 
$q=A^{-2}$, therefore their monomial $q^at^b$ corresponds to the 
free part of the group $H_{i,j}(D;\Z)$ for 
$j=-2b+3w(D)$, $i= -2a +w(D)$ and the monomial $Q^at^b$ corresponds to 
the $\Z_2$ part of the group again with $j=-2b+3w(D)$, $i= -2a +w(D)$. 
KhoHo gives the torsion part of the polynomial for the oriented 
link $6^2_2$, with $w(D)=-6$, as 
$Q^{-6}t^{-1} + Q^{-8}t^{-2} + Q^{-10}t^{-3} + Q^{-12}t^{-4}$.}: 
$H_{6,14}=H_{6,10}=H_{4,10}=\Z$, $H_{2,6}= \Z \oplus \Z_2$, $H_{2,2}=\Z$,
$H_{0,2}= \Z \oplus \Z_2$, $H_{0,-2}= \Z$, $H_{-2,-2} = \Z \oplus \Z_2$,
$H_{-2,-6} = \Z$, $H_{-4,-6}=\Z_2$, $H_{-4,-10}=\Z$, $H_{-6,-10}= H_{-6,-14}=\Z$.
\end{example}

\centerline{\psfig{figure=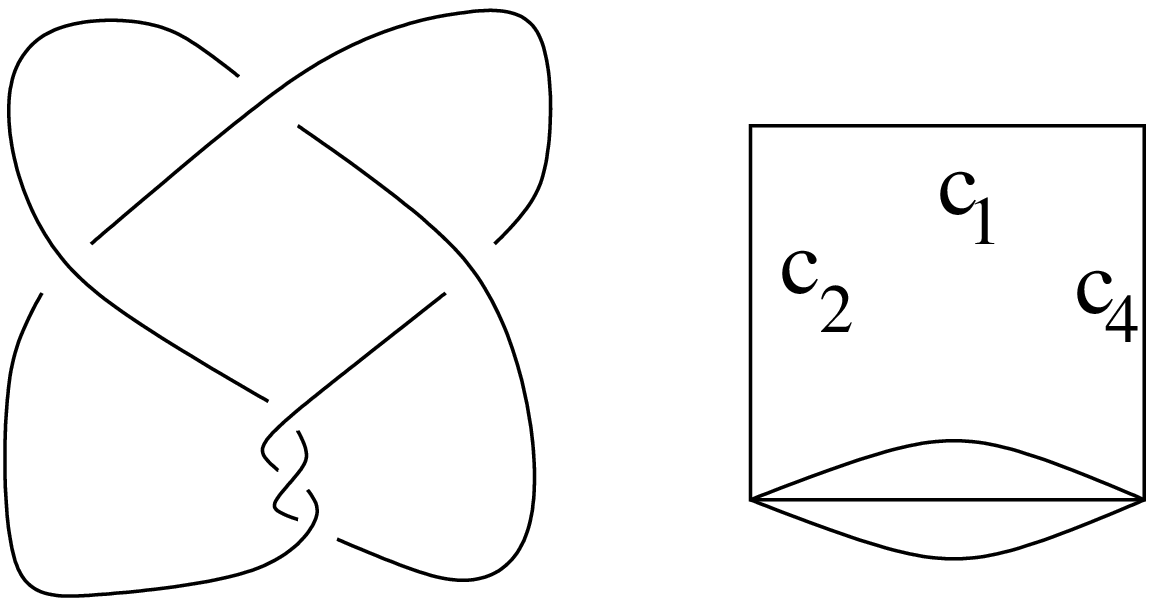,height=4.3cm}}
\begin{center}
Fig. 3.3
\end{center}

\section{Torsion in the Khovanov homology of alternating 
and adequate links}\label{4}
\markboth{\hfil{\sc Khovanov Homology }\hfil}
{\hfil{\sc Torsion for alternating and adequate links}\hfil}

We show in this section how to use technical results from the previous 
sections to prove Shumakovitch's result on torsion in the 
Khovanov homology of alternating links 
and the analogous result for a class of adequate diagrams.

\begin{theorem}[Shumakovitch]\label{4.1} 
The alternating link has torsion free 
Khovanov homology if and only if it is the trivial knot, the Hopf
link or the connected or split sum of copies of them. The nontrivial 
torsion always contains the $\Z_2$ subgroup.
\end{theorem}
The fact that the Khovanov homology of the connected sum of Hopf links 
is a free group, is 
discussed in Section 6 (Corollary 6.6).   

We start with the ``only if" part of the proof by showing 
the following geometric fact.
\begin{lemma}\label{4.2} Assume that $D$ is a link diagram which 
contains a clasp: either $T_{[-2]} = 
\parbox{0.9cm}{\psfig{figure=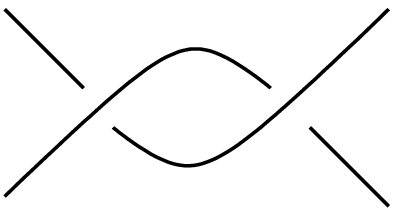,height=0.8cm}}$ \ \ \ \ or 
$T_{[2]} = \parbox{0.9cm}{\psfig{figure=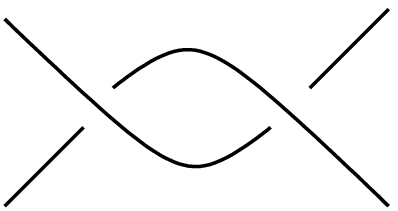,height=0.8cm}}$ 
 \ \ \ \ . 
Assume additionally that the clasp is not a part of the Hopf link summand 
of $D$. Then if the clasp is of 
$T_{[-2]}$ type then the associated graph $G_{s_+}(D)$ has a singular 
edge. If the clasp is of $T_{[2]}$ type then 
the associated graph $G_{s_-}(D)$ has a singular edge. Furthermore the 
singular edge is not a loop.
\end{lemma}
\begin{proof}
Consider the case of the clasp  $T_{[-2]}$, the case of $T_{[2]}$ 
being similar. 
The region bounded by the clasp 
corresponds to the vertex of degree $2$ in $G_{s_+}(D)$. The two edges 
adjacent to this vertex are not loops and they are not singular edges 
only if they share the second 
endpoint as well. In that case our diagram looks like on the 
Fig. 4.1 so it clearly has a Hopf link summand (possibly it is just a 
Hopf link) as the north part is separated by a clasp from the south part 
of the diagram.
\end{proof}
\centerline{\psfig{figure=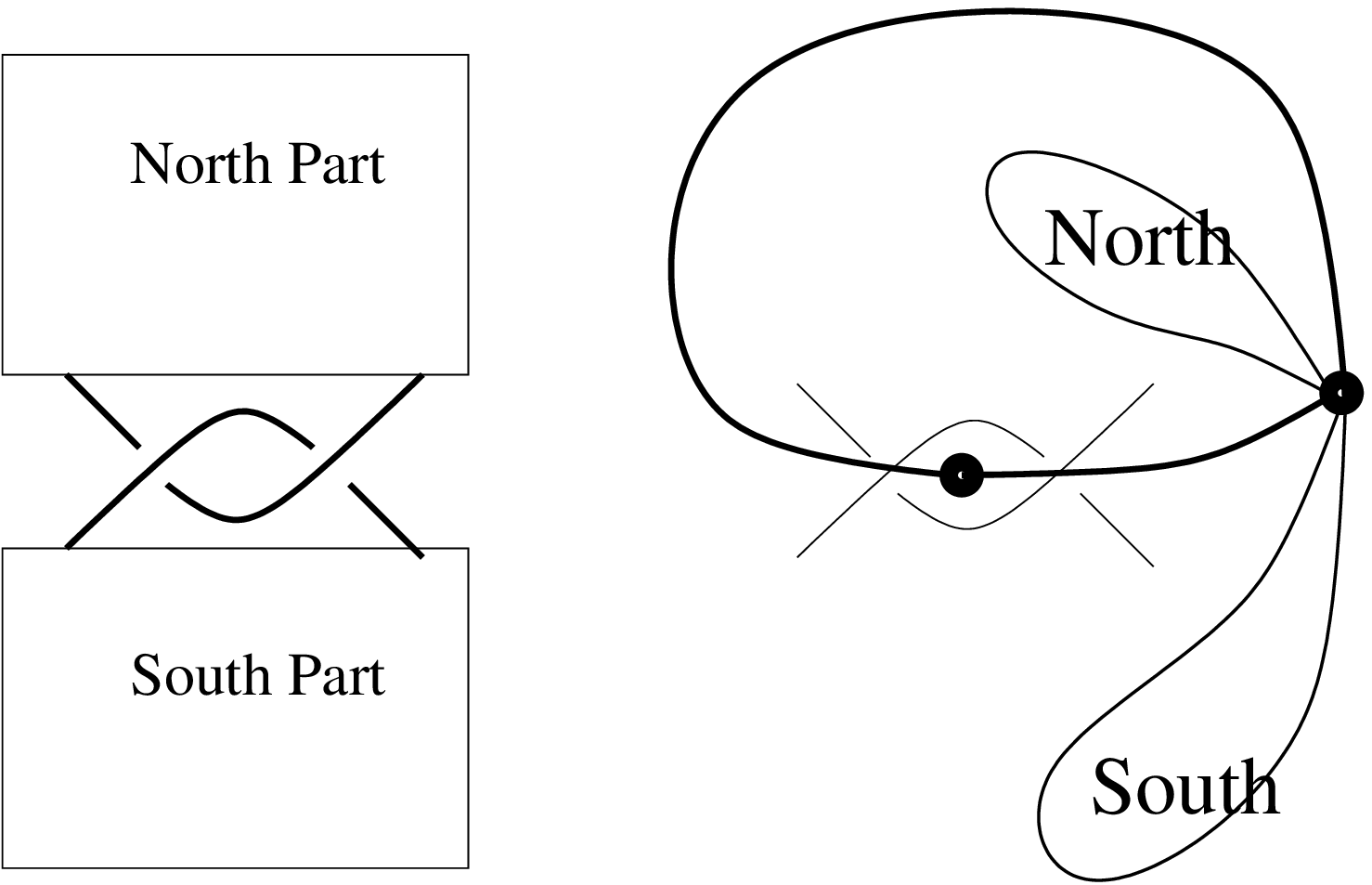,height=3.2cm}}
\begin{center}
Fig. 4.1
\end{center}
\begin{corollary}\label{4.3}
If $D$ is a $+$-adequate diagram (resp. $-$-adequate diagram) with a 
clasp of type $T_{[-2]}$ (resp. $T_{[2]}$), then Khovanov homology 
contains  $\Z_2$-torsion or $T_{[-2]}$ (resp. $T_{[2]}$)) is a part of a 
Hopf link summand of $D$.
\end{corollary}
\begin{proof}
Assume that $T_{[-2]}$ is not a part of Hopf link summand of $D$. By Lemma 
4.2 the graph $G_{s_+}(D)$ has a singular edge. Furthermore,
the graph $G_{s_+}(D)$ has no loops as $D$ is $+$-adequate. If the graph 
has an odd cycle then $H_{N-2,N+2|s_+|-4}(D)$ has $\Z_2$ torsion by Theorem 2.2.
If $G_{s_+}(D)$ is bipartite (i.e. it has only even cycles), then consider 
the cycle containing the singular edge. It is an even cycle of length 
at least 4, so by Theorem 3.2 $H_{N-4,N+2|s_+|-8}(D)$ has $\Z_2$ torsion. 
A similar proof works in $-$-adequate case.
\end{proof}
With this preliminary result we can complete our proof of Theorem 4.1.
 
\begin{proof} First we prove the theorem for non-split, prime 
alternating links. Let $D$ be a diagram of such a link without a 
nugatory crossing. $D$ is an adequate diagram (i.e. it is $+$ and $-$ 
adequate diagram), so it is enough to show that if $G_{s_+}(D)$ (or 
$G_{s_-}(D)$) has a double edge then $D$ can be modified by Tait flypes 
into a diagram with $T_{[-2]}$ (resp. $T_{[2]}$) clasp. This is a standard fact,
justification of which is illustrated in Fig.4.2\footnote{For alternating 
diagrams, $G_{s_+}(D)$ and $G_{s_-}(D)$ are Tait graphs of $D$. These graphs 
are plane graphs and the only possibilities when multiple edges are not 
``parallel" is if our graphs are not 3-connected (as $D$ is not a split link,
graphs are connected, and because $D$ is a prime link, the graphs are 
2-connected). Tait flype corresponds to the special case of change of the 
graph in its 2-isomorphic class as illustrated in Fig.4.2.}. 
\\ \ \\
\centerline{\psfig{figure=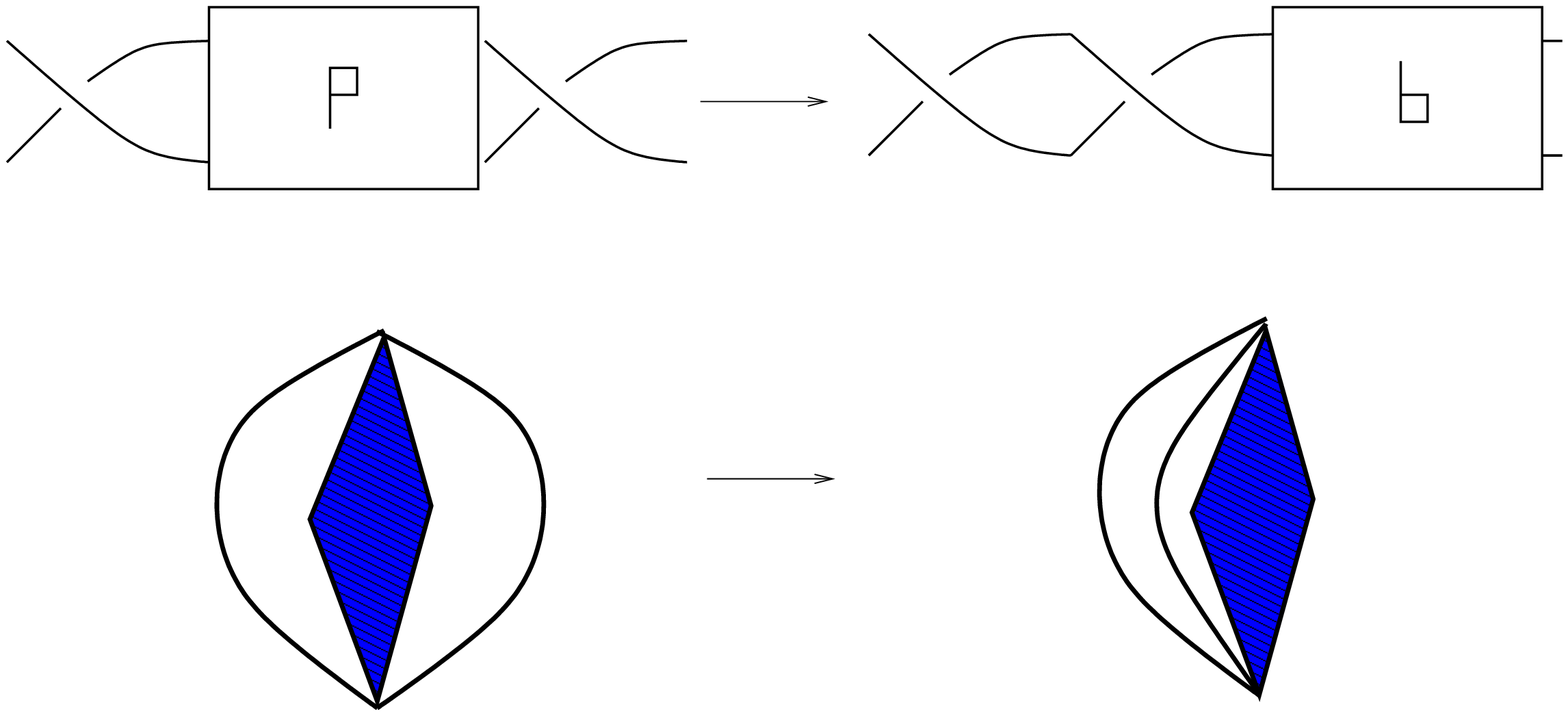,height=4.4cm}}
\begin{center}
Fig. 4.2
\end{center}
If we do not assume that $D$ is a prime link then we use the theorem by 
Menasco \cite{Me} that prime decomposition is visible on the level of 
a diagram. In particular the Tait graphs $G_{s_+}(D)$ and $G_{s_-}(D)$ 
have block structure, where each block (2-connected component) 
corresponds to prime factor of a link. Using the previous results we 
see that the only situation when we didn't find  torsion is if every 
block represents a Hopf link so $D$ represents the sum of Hopf links 
(including the possibility that the graph is just one vertex 
representing the trivial knot).

If we relax condition that $D$ is a non-split link then we use the fact, 
mentioned before, that for $D=D_1\sqcup D_2$, Khovanov homology satisfies 
K\"unneth's formula, $H(D)=H(D_1)\otimes H(D_2)$.
\end{proof}
\begin{example}[The $8_{19}$ knot]\label{4.4} \ \\
The first entry in the knot tables which is not alternating 
is the $(3,4)$ torus knot, $8_{19}$. It is $+$-adequate as it is a positive 
$3$-braid, the closure of $(\sigma_1\sigma_2)^4$. Every positive braid is 
$+$-adequate but its associated graph $G_{s_+}(D)$ 
is composed of 2-gons. Furthermore 
the diagram $D$ of $8_{19}$ is not 
$-$-adequate, Fig.4.3. 
KhoHo shows that the Khovanov homology of $8_{19}$ has 
 torsion, namely $H_{2,2}=\Z_2$.
This torsion is hidden deeply inside the homology spectrum\footnote{The 
full graded homology group is: $H_{8,14}(D)=H_{8,10}(D)= H_{4,6}(D)=\Z$, 
$H_{2,2}=\Z_2$, $H_{0,2}(D)=H_{2,-2}(D)=H_{0,-2}(D)=H_{-2,-4}(D)= 
H_{-2,-10}(D)=\Z$.}, which 
starts from maximum $H_{8,14}(D)=\Z$ and ends on the minimum 
$H_{-2,-10}(D)=\Z$.\\
\end{example}
\centerline{\psfig{figure=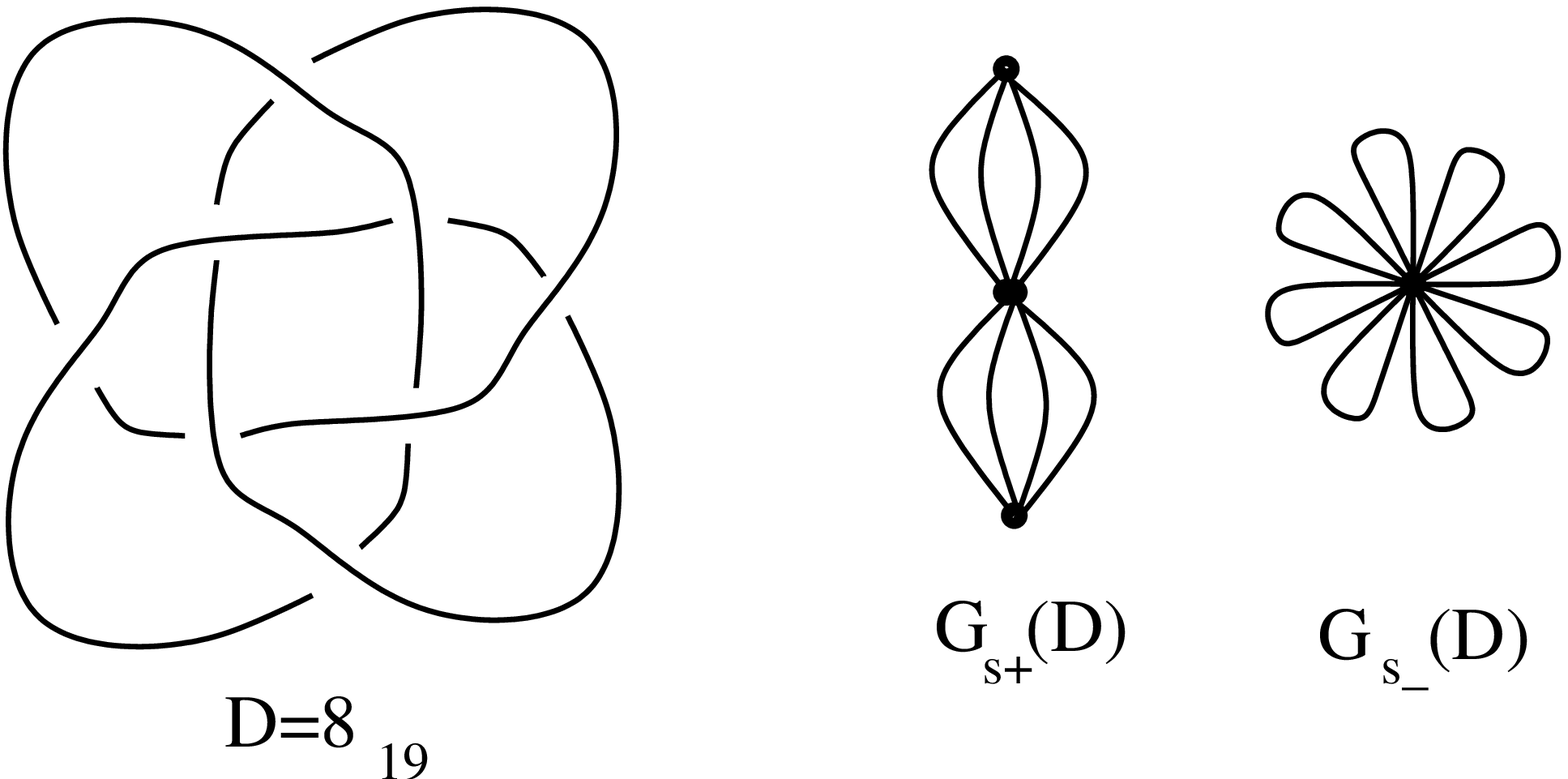,height=5.4cm}}
\begin{center}
Fig. 4.3
\end{center}

\  \\ 

The simplest alternating link which satisfies all conditions of 
Theorem 3.2 except for the existence of a singular edge, is 
the four component alternating link of 8 crossings
$8^4_1$ \cite{Rol}; Fig.4.4. We know that $H_{**}(8^4_1)$ has torsion 
(by using duality) but Theorem 3.2 does not guarantee torsion in 
$H_{N-4,N+ 2|s_+|-8}(D)= H_{4,8}(D)$, the graph $G_{s_+}(D)$
is a square with every edge doubled; Fig.4.3. We checked, however 
using KhoHo the torsion part and in fact $T_{4,8}(D)=\Z_2 $. This 
suggests that Theorem 3.2 can be improved\footnote{In \cite{Ba-3} the 
figure describes, by mistake, the mirror image of $8^4_1$. 
The full homology is as follows:\ 
$H_{8,16}=\Z$, $H_{8,12}=\Z= H_{6,12}$,
$H_{4,8}=\Z_2 \oplus \Z^4$, $H_{4,4}=\Z$, $H_{2,4}=\Z_2^4$, 
$H_{2,0}=\Z^4$, $ H_{0,0}=\Z^7$, $H_{0,-4}=\Z^6$, 
$H_{-2,-4}=\Z_2^3 \oplus \Z^3$, $H_{-2,-8}=\Z$, $H_{-4,-8}=\Z^3$,
$H_{-4,-12}=\Z^3$, $H_{-6,-12}=\Z_2^3$, 
$H_{-6,-16}=\Z^3$, $H_{-8,-16}=\Z$, $H_{-8,-20}=\Z$. In KhoHo the 
generating polynomials, assuming $w(8^4_1)=-8$, are:
KhPol("8a",21)= $[((q^{18} + q^{16})*t^8 + 3*q^{16}*t^7 + (3*q^{14} + 
3*q^{12})*t^6 + (q^{12} + 3*q^{10})*t^5 + (6*q^{10} + 7*q^8)*t^4 + 
4*q^8*t^3 + (q^6 + 4*q^4)*t^2 + q^2*t + (q^2 + 1))/(q^{20}*t^8)$,\\
$(3*Q^{10}*t^5 + 3*Q^8*t^4 + Q^6*t^3 + 4*Q^2*t + 1)/(Q^{16}*t^6)].$}.
\ \\
\ \\

\centerline{\psfig{figure=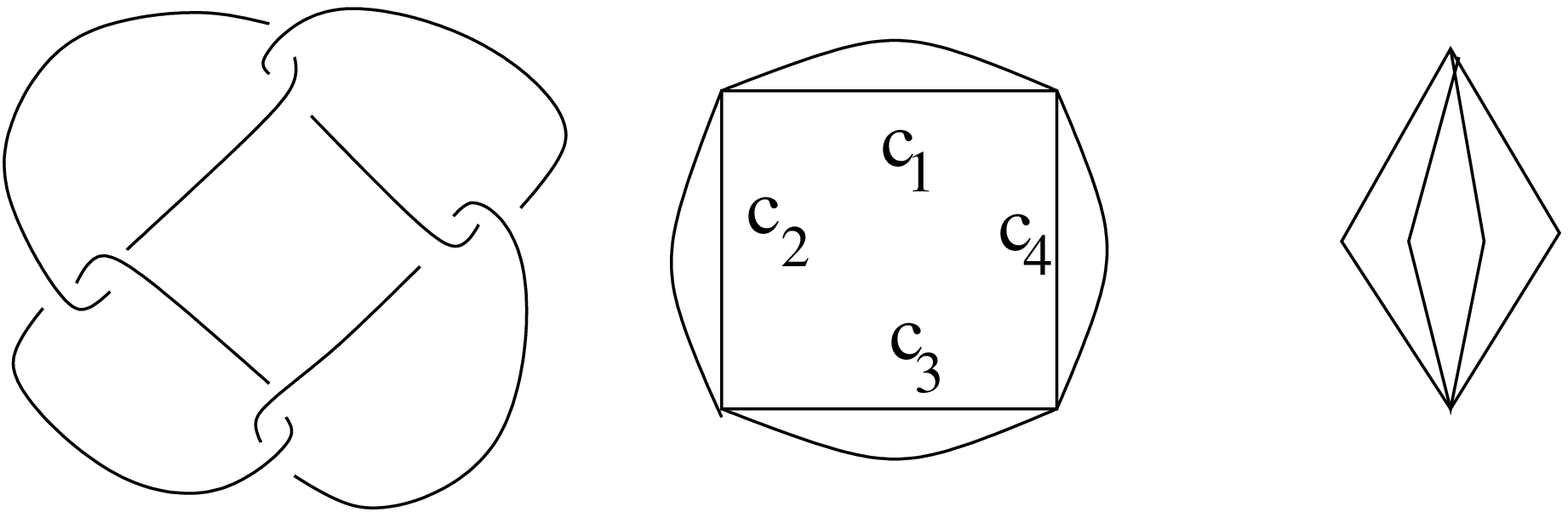,height=4.4cm}}
\begin{center}
Fig. 4.4
\end{center}


\section{Thickness of Khovanov homology 
and almost alternating links}\label{5}
\markboth{\hfil{\sc Khovanov Homology }\hfil}
{\hfil{\sc Thickness of Khovanov homology links}\hfil}

We define, in this section, the notion of an $H$-$k$-thick link diagram and
relate it to $(k-1)$-almost alternating diagrams. In particular we 
give a short proof of Lee's theorem \cite{Lee-1} (conjectured by Khovanov, 
Bar-Natan, and 
Garoufalidis) that alternating non-split links are $H$-$1$-thick 
($H$-thin in Khovanov terminology). 

\begin{definition}\label{5.1}
We say that a link is $k$-almost alternating if it has a diagram 
which becomes alternating after changing $k$ of its crossings.
\end{definition}
As noted in Property 1.4 the ``extreme" terms of Khovanov chain comples 
are $C_{N,N+2|s_+|}(D)= C_{-N,-N-2|s_-|}(D)= \Z$. In the 
following definition 
of a $H$-$(k_1,k_2)$-thick diagram we compare indices of actual 
Khovanov homology of $D$ with lines of slope $2$ going through the 
points $(N,N+2|s_+|)$ and $(-N,-N-2|s_-|)$.
\begin{definition}\label{5.2}
\begin{enumerate}
\item[(i)] 
We say that a link diagram, $D$ of $N$ crossings is $H$-$(k_1,k_2)$-thick 
if $H_{i,j}(D)=0$ 
with a possible exception of $i$ and $j$ satisfying:
$$ N- 2|s_-| - 4k_2 \leq j-2i \leq 2|s_+| -N + 4k_1.$$

\item[(ii)]
We say that a link diagram of $N$ crossings is 
$H$-$k$-thick\footnote{Possibly, the more appropriate name would be 
$H$-$k$-thin diagram, as the width of Khovanov homology is bounded from 
above by $k$. Khovanov (\cite{Kh-2}, page 7) suggests the term 
homological width; $hw(D)=k$ if homology of $D$ lies on $k$ adjacent 
diagonals (in our terminology, $D$ is $k-1$ thick).}
 if,
it is $H$-$(k_1,k_2)$-thick where $k_1$ and $k_2$ satisfy:
$$k\geq k_1 + k_2 + \frac{1}{2} (|s_+| + |s_-| - N).$$
\item[(iii)] 
We define also $(k_1,k_2)$-thickness (resp. $k$-thickness) 
of Khovanov homology separately for the torsion part 
(we use the notation $TH$-$(k_1,k_2)$-thick 
diagram), and for the free part (we use the notation 
$FH$-$(k_1,k_2)$-thick diagram).
\end{enumerate} 
\end{definition}

Our $FH$-$1$-thick diagram is a $H$-thin diagram 
in \cite{Kh-2,Lee-1,Ba-2,Sh-1}.

With the above notation we are able to formulate our main result of 
this section.
\begin{theorem}\label{5.3}\ \\
If the diagram $D_{\infty}= 
D(\parbox{0.9cm}{\psfig{figure=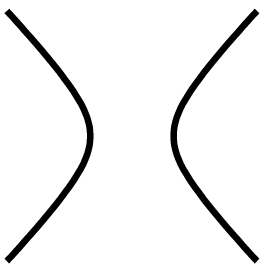,height=0.7cm}})$ 
is $H$-$(k_1(D_{\infty}),
k_2(D_{\infty}))$-thick 
and the diagram 
$D_0 = D(\parbox{0.9cm}{\psfig{figure=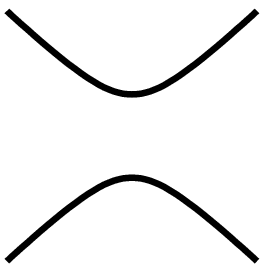,height=0.7cm}})$ 
is $H$-$(k_1(D_0),
k_2(D_0))$-thick, 
then the diagram 
$D_+ = D(\parbox{0.9cm}{\psfig{figure=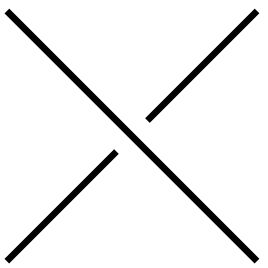,height=0.7cm}})$ 
is $H$-$(k_1(D_+),
k_2(D_+))$-thick 
where \\
$k_1(D_+)= 
\max(k_1(D_{\infty}) + 
\frac{1}{2}(|s_+(D_{\infty})| - 
|s_+(D_+)| +1), 
k_1(D_0))$ and .\\
$k_2(D_+)=
\max(k_2(D_{\infty}),
k_2(D_0)) + \frac{1}{2}(|s_-(D_0)| - |s_-(D_+)| +1))$.\\ 
In particular 
\begin{enumerate}
\item[(i)] if $|s_+(D_+)| - |s_+(D_{\infty})| =1$, as is always 
the case for a $+$-adequate diagram, then 
$k_1(D_+)= \max(k_1(D_{\infty}), k_1(D_0))$,
\item[(ii)] if $|s_-(D_+)| - |s_-(D_0| =1$, as is always
the case for a $-$-adequate diagram, then
$k_2(D_+)= \max(k_2(D_{\infty}), k_2(D_0))$.
\end{enumerate}
\end{theorem}
\begin{proof}
We formulated our definitions so that our proof follows almost 
immediately via the Viro's long exact sequence of Khovanov homology:
$$...\to H_{i+1,j-1}(D_0) \stackrel{\partial}{\to}  
 H_{i+1,j+1}(D_{\infty}) \stackrel{\alpha}{\to} 
H_{i,j}(D_+)
 \stackrel{\beta}{\to}$$
$$H_{i-1,j-1}(D_0) \stackrel{\partial}{\to} H_{i-1,j+1}(D_{\infty})
 \to ...$$
If $0\neq h \in H_{i,j}(D_+)$ then 
either $h=\alpha (h')$ for $0\neq h' \in H_{i+1,j+1}(
D_{\infty})$ or 
 $0\neq \beta(h) \in H_{i-1,j-1}(D_0)$. Thus if
$H_{i,j}(D_+)\neq 0$ 
then either $H_{i+1,j+1}(D_{\infty})\neq 0$ or 
$H_{i-1,j-1}(D_0)\neq 0$.
The first possibility gives the inequalities involving $(j+1)- 2(i+1)$:
$$N(D_{\infty}) - 2|s_-(D_{\infty})| -4k_2(D_{\infty})   
\leq j-2i-1 
\leq 2|s_+(D_{\infty})| - N(D_{\infty}) +4k_1(D_{\infty})$$
which, after observing that $|s_-(D_+)|=|s_-(D_{\infty})|$, leads to:
$$ N(D_+) - 2|s_-(D_+)| - 4k_2(D_{\infty})   \leq j-2i \leq $$ $$ 
2|s_+(D_+)| - N(D_+) + 4k_1(D_{\infty}) + 2(|s_+(D_{\infty})| - 
|s_+(D_+)| +1).$$
The second possibility gives the inequalities involving $(j-1)- 2(i-1)$:
$$N(D_{0}) - 2|s_-(D_{0})| -4k_2(D_{0})  \leq j-2i+1
\leq 2|s_+(D_{0})| - N(D_{0}) +4k_1(D_{0})$$
which, after observing that $|s_+(D_+)|=|s_+(D_0)|$, leads to:
$$ N(D_+) - 2|s_-(D_+)| - 4k_2(D_{0}) 
-2(|s_-(D_{0})| -|s_-(D_+)| +1) \leq  $$\ $ j-2i \leq \ 
2|s_+(D_+)| - N(D_+) + 4k_1(D_{0}).$

Combining these two cases we obtain the conclusion of Theorem 5.3.
\end{proof}
\begin{corollary}\label{5.4}
If $D$  is an adequate diagram such that, for some crossing of $D$, 
the diagrams $D_0$ and $D_{\infty}$ are $H$-$(k_1,k_2)$-thick (resp. 
$H$-$k$-thick) then $D$ is $H$-$(k_1,k_2)$-thick (resp.  $H$-$k$-thick).
\end{corollary}
\begin{corollary}\label{5.5}
Every alternating non-split diagram without a nugatory 
crossing is $H$-$(0,0)$-thick and $H$-$1$-thick.
\end{corollary}
\begin{proof}
The $H$-$(k_1,k_2)$-thickness in 
Corollary 5.4 follows immediately from Theorem 5.3. To show 
$H$-$k$-thickness we observe additionally that for an adequate 
diagram $D_+$ one has $|s_+(D_0)| + |s_-(D_0)| - N(D_0) =
|s_+(D_+)| + |s_-(D_+)| - N(D_+) = 
|s_+(D_{\infty})| + |s_-(D_{\infty})| - N(D_{\infty})$.

 We prove 
Corollary 5.5 using induction on the number of crossings a slightly
more general statement allowing nugatory crossings.\\ 
(+) If $D$ is an alternating non-split $+$-adequate diagram 
then $H_{i,j}(D)\neq 0$ 
can happen only for $j-2i \leq 2|s_+(D)| - N(D)$.\\
(--) If $D$ is an alternating non-split $-$-adequate diagram 
then $H_{i,j}(D)\neq 0$
can happen only for $N(D) - 2|s_-(D)| \leq j-2i$.\\
If the diagram $D$ from (+) has only nugatory crossings then it represents 
the trivial knot and its nontrivial Khovanov homology are 
$H_{N,3N-2}(D)=H_{N,3N+2}(D)=\Z$. Because $|s_+(D)|=N(D) + 1$ in this case, 
the inequality (+) holds.  In the inductive step we use the property of 
a non-nugatory crossing of a non-split $+$-adequate diagram, namely 
$D_0$ is also an alternating non-split $+$-adequate diagram and inductive 
step follows from Theorem 5.3.\\
Similarly one proves the condition (--). Because the non-split 
alternating diagram without nugatory crossings is an adequate diagram, 
therefore Corollary 5.5 follows from Conditions (+) and (--). 
\end{proof}
 
The conclusion of the theorem is the same if we are interested only in 
the free part of Khovanov homology (or work over a field). In the case 
of the torsion part of the homology we should take into account the 
possibility 
that torsion ``comes" from the free part of the homology, that is 
$H_{i+1,j+1}(D_{\infty})$ may be torsion free but its image under $\alpha$ 
may have torsion element.

\begin{theorem}\label{5.6}
If $T_{i,j}(D_+) \neq 0$ then either \\
(1)\  $T_{i+1,j+1}(D_{\infty}) \neq 0$ or 
$T_{i-1,j-1}(D_{0}) \neq 0$, \\
or\\
(2)\ $FH_{i+1,j+1}(D_{\infty}) \neq 0$ and 
$FH_{i+1,j-1}(D_0) \neq 0$.
\end{theorem}
\begin{proof}
From the long exact sequence of Khovanov homology it follows that the only 
way the torsion is not related to the torsion of $H_{i+1,j+1}(D_{\infty})$ or 
$H_{i-1,j-1}(D_{0})$ is the possibility of torsion created by taking 
the quotient \\ 
$FH_{i+1,j+1}(D_{\infty})/\partial(FH_{i+1,j-1}(D_0))$ and 
in this case both groups $FH_{i+1,j+1}(D_{\infty})$ and 
$FH_{i+1,j-1}(D_0)$ have to be nontrivial.
\end{proof}

\begin{corollary}\label{5.7}
If $D$ is an alternating non-split diagram without a nugatory crossing 
then $D$ is $TH$-$(0,-1)$-thick and $TH$-$0$-thick. In other words if 
$T_{i,j}(D)\neq 0$ then $j-2i = 2|s_+(D)| -N(D)= N(D) -2|s_-(D)|+4$.
\end{corollary}
\begin{proof} We proceed in the same (inductive) manner as in the proof 
of Corollary 5.5, using Theorem 5.7 and Corollary 5.5. In the first step 
of the induction we use the fact that the trivial knot has no torsion 
in Khovanov homology.
\end{proof}

The interest in $H$-thin diagrams was 
motivated by the observation (proved by Lee) that diagrams of non-split 
alternating links are $H$-thin (see Corollary 5.5). Our approach allows 
the straightforward generalization to $k$-almost alternating diagrams.

\begin{corollary}\label{5.8} 
Let $D$ be a non-split $k$-almost alternating diagram without a nugatory 
crossing. Then $D$ is $H$-$(k,k)$-thick and $TH$-$(k,k-1)$-thick.
\end{corollary}
\begin{proof}
The corollary holds for $k=0$ (alternating diagrams) and we use an induction 
on the number of crossings needed to change the diagram $D$ to an alternating 
digram, using Theorem 5.3 in each step.
\end{proof}

We were assuming throughout the section that our diagrams are non-split. 
This assumption was not always necessary. In particular even the split 
alternating diagram without nugatory crossings is $H$-$(0,0)$-thick as 
follows from the following observation.
\begin{lemma}\label{5.9}
If the diagrams $D'$ and $D^{\prime\prime}$ are $H$-$(k_1',k_2')$-thick and 
$H$-$(k_1^{\prime\prime},k_2^{\prime\prime})$-thick, respectively, 
then the diagram $D=D' \sqcup D^{\prime\prime}$ 
is $H$-$(k_1'+k_1^{\prime\prime},k_2'+k_2^{\prime\prime})$-thick.
\end{lemma}
\begin{proof}
Lemma 5.9 follows from the obvious fact that in the split sum 
$D=D' \sqcup D''$ we have 
$N(D) = N(D') + N(D^{\prime\prime})$, $|s_+(D)| = |s_+(D')| + 
|s_+(D^{\prime\prime})|$ and 
$|s_-(D)| = |s_-(D')| + |s_-(D^{\prime\prime})|$.
\end{proof}

Khovanov observed (\cite{Kh-2}, Proposition 7) that adequate 
non-alternating knots are not $H$-$1$-thick. We are able to proof 
the similar result about torsion of adequate non-alternating links.

\begin{theorem}\label{5.10}\ \\
Let $D$ be a connected adequate diagram which  does not represent an 
alternating link and  such that $G_{s_+}(D)$ and 
$G_{s_-}(D)$ have either an odd cycle or an even cycle with a 
singular edge, then $D$ is not $TH$-$0$-thick diagram. More generally,
$D$ is at best $TH$-$\frac{1}{2}(N+2-(|s_+(D)| + |s_-(D)|)$-thick.
\end{theorem}
\begin{proof}
The first part of Theorem 5.10 follows from the second part because by 
Proposition 1.4 (Wu's Lemma), $\frac{1}{2}(N+2-(|s_+(D)| + |s_+(D)|) > 0$ 
for a diagram which is not a connected sum of alternating diagrams. 
By Theorems 2.2, 3.2 and Corollary 3.3, $TH_{i,j}(D)$ is nontrivial 
on slope $2$ diagonals $j-2i = 2|s_+|-N$ and $N- 2|s_-| +4$. The $j$ 
distance between these diagonals is $N- 2|s_-| +4 - (2|s_+|-N) = 
2(N+2-(|s_+(D)| + |s_+(D)|)$, so the theorem follows. 
\end{proof}
\begin{example}\label{5.11}\ \\
Consider the knot $10_{153}$ (in the notation of \cite{Rol}). 
It is an adequate 
non-alternating knot. Its associated graphs $G_{s_+}(10_{153})$ and 
$G_{s_-}(10_{153})$ have triangles (Fig.5.1) so Theorem 5.10 applies. 
Here $|s_+|=6$, $|s_-|=4$ and by Theorem 2.2, $H_{8,18}(10_{153})$ and
$H_{-10,-14}(10_{153})$ have $\Z_2$ torsion. Thus support of torsion 
requires at least $2$ adjacent diagonals\footnote{Checking \cite{Sh-2},
gives the full torsion of the Khovanov homology of $10_{153}$ as: \ 
$T_{8,18}=T_{4,10}=T_{2,6}=T_{0,6}=T_{-2,-2}=T_{-4,-2}=T_{-6,-6}= 
T_{-10,-14}= \Z_2$.}
\end{example}

\centerline{\psfig{figure=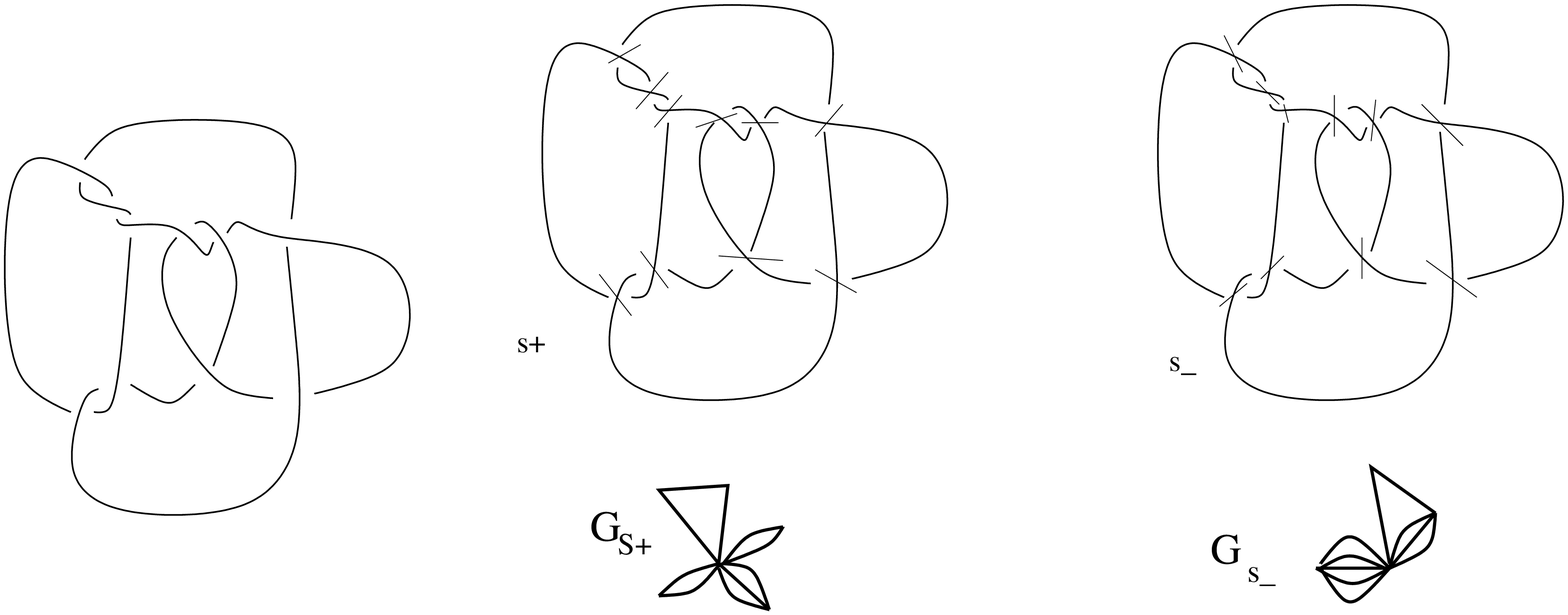,height=6.3cm}}
\begin{center}
Fig. 5.1\end{center}
\begin{corollary}\label{5.12}
Any doubly adequate link which is not an alternating link is not 
$TH$-$0$-thick.
\end{corollary}
 
\section{Hopf link addition}\label{6}
We find, in this section, the structure of the Khovanov homology of 
connected sum of $n$ copies of the Hopf link, as promised in Section 5. 
As a byproduct of our method, we are able to compute Khovanov homology 
of a connected sum of a diagram $D$ and the Hopf link $D_h$, Fig 6.1,  
confirming a conjecture by A.Shumakovitch that the 
Khovanov homology of the connected sum of D with the Hopf link,
is the double of the 
Khovanov homology of $D$.
\\
\ \\
\centerline{\psfig{figure=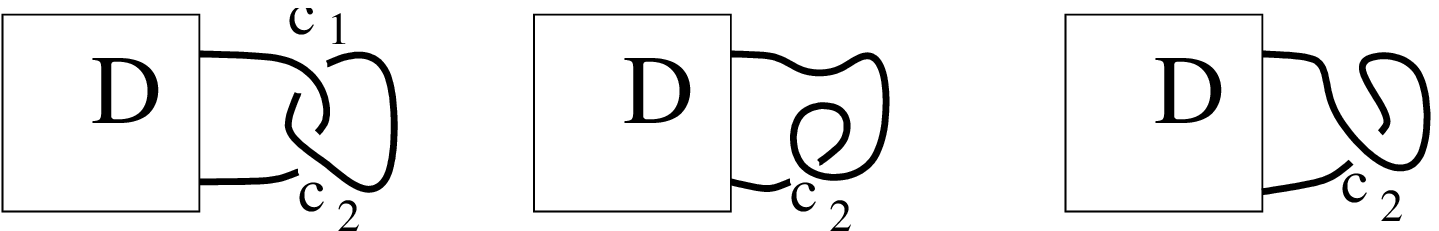,height=2.2cm}}
\centerline{ \ \ \  $D\#D_h$ \ \ \ \ \ \ \ \ \ \ \ \ \ \ \ \ \ \ \ \ 
 $(D\#D_h)_0$ \ \ \ \ \ \ \ \ \ \  \ \ \ \ \ \ \ 
 $(D\#D_h)_{\infty}$  \ \ \ \ } 
\begin{center}
Fig. 6.1
\end{center}

\begin{theorem}\label{6.1}
For every diagram $D$ we have the short exact sequence of Khovanov 
homology\footnote{Theorems 6.1 and 6.2 hold for a diagram $D$ on any surface 
$F$ and for any ring of coefficients $\cal R$ with the restriction 
that for $F=RP^2$ we need $2{\cal R}=0$. In this more general case of 
a manifold being $I$-bundle over a surface, we use definitions and 
setting of Section 8.}
$$0\to H_{i+2,j+4}(D) \stackrel{\alpha_h}{\to} 
H_{i,j}(D\#D_h) \stackrel{\beta_h}{\to} H_{i-2,j-4}(D) \to 0$$
where $\alpha_h$ is given on a state $S$ by Fig.6.2(a)  and $\beta_h$ is a 
projection given by Fig.6.2(b) (and $0$ on other states). The theorem 
holds for any ring of coefficients, $\cal R$, not just ${\cal R}=\Z$.
\end{theorem}
  \ \\
\centerline{\psfig{figure=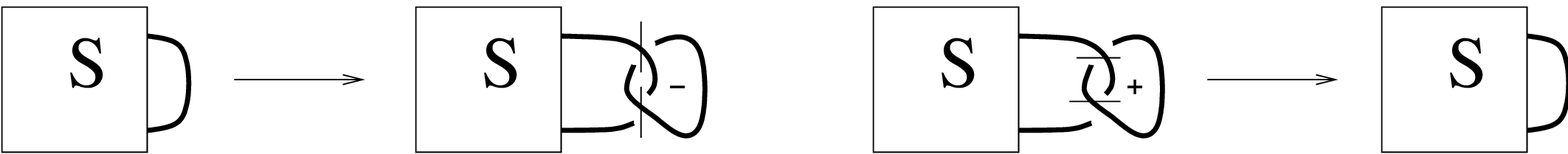,height=1.5cm}}
\ \\
\centerline{ (a) \ \ \ \ \ \ \  \ \ \ \ \ \ \ \ 
$\alpha_h$ \ \ \ \ \ \ \ \ \ \ \ \ \ \ \ \ \ \ \ \ \ \ \ \ \ \ \ \ \ \
 (b) \ \ \ \ \ \ \ \ \ \  \ \ \ \ \ \ \ \ \ \ \ \ \ \ \ \ \ \ \ \
 $\beta_h$  \ \ \ \ \ \ \ \ \ \ \ \ \ } 
\begin{center}
Fig. 6.2
\end{center}
\begin{theorem}\label{6.2}\ \\
The short exact sequence of homology from Theorem 6.1 splits, so we have
$$H_{i,j}(D\#D_h) = H_{i+2,j+4}(D) \oplus H_{i-2,j-4}(D).$$
\end{theorem}
\begin{proof} To prove Theorem 6.1 we 
consider the long exact sequence of the Khovanov homology of the diagram 
$D\#D_h$ with respect to the first crossing of the diagram, $e_1$ (Fig.6.1). 
To simplify the notation we assume that ${\cal R}=\Z$ but our 
proof works for any ring of coefficients.
$$ ... \to H_{i+1,j-1}((D\#D_h)_0) \stackrel{\partial}{\to} 
H_{i+1,j+1}((D\#D_h)_{\infty}) \stackrel{\alpha}{\to} H_{i,j}(D\#D_h) 
\stackrel{\beta}{\to} $$ 
$$H_{i-1,j-1}((D\#D_h)_0) \stackrel{\partial}{\to} 
H_{i-1,j+1}((D\#D_h)_{\infty}) \to ...$$

We show that the homomorphism $\partial$ is the zero map. We 
use the fact that $(D\#D_h)_0$ differs from $D$ by a positive first 
Reidemeister move $R_{+1}$ and that $(D\#D_h)_{\infty}$ 
differs from $D$ by a negative first
Reidemeister move $R_{-1}$; Fig.6.1. 
 We know, see 
\cite{APS-2} for example, that the chain map \\
$r_{-1}: {\cal C}(D) \to {\cal C}(R_{-1}(D))$ 
given by 
$r_{-1}(\parbox{0.8cm}{\psfig{figure=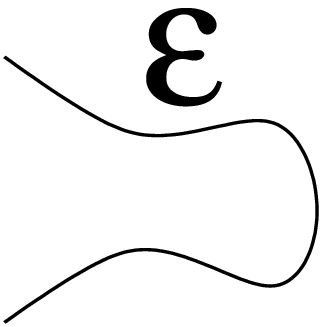,height=0.6cm}}) = 
\parbox{1.2cm}{\psfig{figure=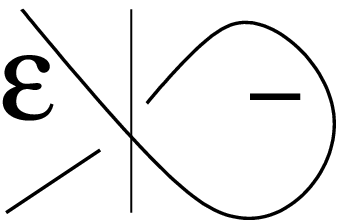,height=0.6cm}}$ 
yields the isomorphism of homology:
$$r_{-1*}: H_{i,j}(D) \to H_{i-1,j-3}(R_{-1}(D))$$ and the chain map 
$\bar{r}_{+1}({\cal C}(R_{+1}(D))={\cal C}((D)$ 
given by the projection with $\bar{r}_{+1}(
\parbox{1.2cm}{\psfig{figure=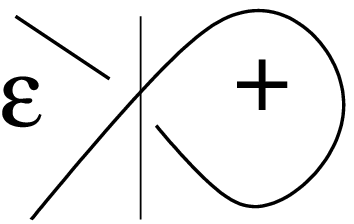,height=0.6cm}})= 
(\ \parbox{1.1cm}{\psfig{figure=kinksmoothe.eps,height=0.6cm}})$ and $0$ 
otherwise,  
induces the isomorphism of homology:
$$\bar r_{+1*}: H_{i+1,j+3}(R_{+1}(D)) \to H_{i,j}(D).$$ 
From these we get immediately that the composition  homomorphism:
$$r_{-1*}^{-1}\partial {\bar r}_{+1*}^{-1}: H_{i,j-4}(D) \to H_{i+2,j+4}(D)$$ 
is the zero map by considering the composition of homomorphisms 
$$ H_{i,j-4}((D)) \stackrel{{\bar r}_{+1*}^{-1}}{\to}
H_{i+1.j-1}((D\#D_h)_{0}) \stackrel{\partial}{\to}
H_{i+1.j+1}((D\#D_h)_{\infty})) \stackrel{r_{-1*}^{-1}}{\to}
H_{i+2,j+4}(D).$$
\end{proof}
Let $h(a,b)(D)$ (resp. $h_{{\cal F}}(a,b)(D)$ for a field ${\cal F}$) 
be the generating polynomial of the free part of $H_{**}(D)$ (resp. 
$H_{**}(D;{\cal F})$),
where $kb^ia^j$  (resp. $k_{{\cal F}}b^ia^j$) 
represents the fact that the free 
part of $H_{i,j}(D)$, $FH_{i,j}(D)= \Z^k$ 
(resp. $H_{i,j}(D;{\cal F})= {\cal F}^k$).

Theorem 6.2 will be proved in several steps.
\begin{lemma}\label{6.3}\ \\
If the module $H_{i-2,j-4}(D;{\cal R})$ is free (e.g. ${\cal R}$ is 
a field) then the sequence from Theorem 6.1 splits and 
 $H_{i,j}((D\#D_h);{\cal R})=H_{i-2,j-4}(D;{\cal R})\oplus 
H_{i+2,j+4}(D;{\cal R})$ or shortly
 $H_{**}(D\#D_h;{\cal R})= 
H_{**}(D;{\cal R})(b^2a^4+ b^{-2}a^{-4})$. 

For the free part we have always $FH_{i,j}((D\#D_h) = 
FH_{i+2,j+4}(D) \oplus FH_{i-2,j-4}(D))$ or in the 
language of generating functions: 
$h(a,b)(D\#D_h)=(b^2a^4+ b^{-2}a^{-4})h(a,b)(D)$.
\end{lemma}
\begin{proof}
The first part of the lemma follows immediately from Theorem 6.1 which holds 
for any ring of coefficients, in particular 
$rank(FH_{i,j}(D\#D_h)=rank(FH_{i+2,j+4}(D)) + rank(FH_{i-2,j-4}(D))$.
\end{proof}
\begin{lemma}\label{6.4}\ \\
There is the exact sequence of $\Z_p$ linear spaces:
$$ 0 \to H_{i+2,j+4}(D)\otimes \Z_p \to H_{i,j}(D\# D_h)\otimes \Z_p
\to H_{i-2,j-4}(D)\otimes \Z_p \to 0 .$$
\end{lemma}
\begin{proof}
Our main tool is the 
universal coefficients theorem
(see, for example, \cite{Ha}; Theorem 3A.3) combined with Lemma 6.3. 
By the second part of Lemma 6.3 it suffices to prove that:
$$T_{i,j}(D\# D_h)\otimes \Z_p = T_{i+2,j+4}(D)\otimes \Z_p \oplus 
T_{i-2,j-4}(D)\otimes \Z_p.$$
From the universal coefficients theorem we have:
$$H_{i,j}(D\#D_h);\Z_p) = 
H_{i,j}(D\#D_h)\otimes \Z_p \oplus  Tor(H_{i-2,j}(D\#D_h), \Z_p)$$ and 
$Tor(H_{i-2,j}(D\#D_h), \Z_p) = T_{i-2,j}(D\#D_h)\otimes \Z_p$ and 
the analogous formulas for the Khovanov homology of $D$. 
Combining this with both parts of 
Lemma 6.3, we obtain:
$$ T_{i,j}(D\# D_h)\otimes \Z_p \oplus T_{i-2,j}(D\# D_h)\otimes \Z_p = $$
$$ 
(T_{i+2,j+4}(D)\otimes \Z_p \oplus T_{i,j+4}(D)\otimes \Z_p) \oplus 
(T_{i-2,j-4}(D)\otimes \Z_p \oplus T_{i-4,j-4}(D)\otimes \Z_p).$$ 
We can express this 
in the language of generating functions assuming that $t(b,a)(D)$ is the 
generating function of dimensions of $T_{i,j}(D)\otimes \Z_p$:
$$(1+b^{-2})t(b,a)(D\# D_h) = (1+b^{-2})(b^2a^4 + b^{-2}a^{-4}t(b,a)(D).$$
Therefore $t(b,a)(D\# D_h) = (b^2a^4 + b^{-2}a^{-4})t(b,a)(D)$ and 
 $dim(T_{i,j}(D\# D_h)\otimes \Z_p) = 
dim(T_{i+2,j+4}(D)\otimes \Z_p + dim(T_{i,j+4}(D)\otimes \Z_p)$.
The lemma follows by observing that 
the short exact sequence with $\Z$ coefficients
leads to the sequence
$$0 \to ker(\alpha_p) \to H_{i+2,j+4}(D)\otimes 
\Z_p \stackrel{\alpha_p}{\to} $$ 
$$H_{i,j}(D\# D_h)\otimes \Z_p
\to H_{i-2,j-4}(D)\otimes \Z_p \to 0 .$$ 
By the previous computation $dim(ker(\alpha_p))=0$ and the 
proof is completed.
\end{proof}
To finish our proof of Theorem 6.2 we only need the following lemma.
\begin{lemma}\label{6.5}
Consider a short exact sequence of finitely generated abelian groups:
$$0\to A \to B \to C \to 0. $$ If for every prime number $p$ we 
have also the exact sequence: 
$$0\to A\otimes \Z_p \to B\otimes \Z_p \to C\otimes \Z_p \to 0 $$ then 
the exact sequence $0\to A \to B \to C \to 0 $ splits and $B= A\oplus C$.
\end{lemma}
\begin{proof} Assume, for contradiction, that the sequence 
$0\to A \stackrel{\alpha}{\to} B \to C \to 0 $ does not split. Then there 
is an element $a \in A$ such that ${\alpha}(a)$ is not $p$-primitive in $B$, 
that is ${\alpha}(a)=pb$ for $b\in B$ and $p$ a prime number and 
$b$ does not lies in the subgroup of $B$ span by ${\alpha}(a)$ 
(to see that such an $a$ exists one 
can use the maximal decomposition of $A$ and $B$ into cyclic subgroups 
(e.g. $A= \Z^{k}\oplus_{p,i} \Z_{p^i}^{k_{p,i}}$)). Now comparing dimensions 
of linear spaces $A\otimes \Z_p, B\otimes \Z_p, C\otimes \Z_p$ 
(e.g. $dim(A\otimes \Z_p)= k+k_{p,1}+k_{p,2}+...$ we see that the sequence 
$0\to A\otimes \Z_p 
\to B\otimes \Z_p \to C\otimes \Z_p \to 0 $ is not exact, a contradiction.
\end{proof}

\begin{corollary}\label{6.6}
For the connected sum of $n$ copies of the Hopf link 
we get\footnote{In the oriented version
(with the linking number equal to $n$, so the writhe number $w=2n$) and
with Bar-Natan notation one gets:
$q^{3n}t^n(q+q^{-1})(q^2t + q^{-2}t^{-1})^n$, as computed first by
Shumakovitch.}\\ 
$H_{*,*}(D_h \# ...\# D_h) = h(a,b)(D)=(a^2+a^{-2})(a^4b^2 + a^{-4}b^{-2})^n$
\end{corollary}

\begin{remark}\label{6.7}
Notice that $h(a,b)(D_h)- h(a,b)(OO) = (a^2+a^{-2})(a^4b^2 + a^{-4}b^{-2}) - 
(a^2+a^{-2})^2= b^{-2}a^{-4}(a^2+a^{-2})(1+ba)(1-ba)(1+ba^3)(1-ba^3)$.
This equality may serve as a starting point to formulate a conjecture 
for links, analogous to Bar-Natan-Garoufalidis-Khovanov conjecture 
\cite{Kh-2,Ga},\cite{Ba-2} (Conjecture 1), formulated for knots and 
proved for alternating knots by Lee \cite{Lee-1}.
\end{remark}

\section{Reduced Khovanov homology}\label{7}
Most of the results of Sections 5 and 6 can be adjusted to the case of reduced 
Khovanov homology\footnote{Introduced by Khovanov; we follow here 
Shumakovitch's approach adjusted to the framed link version.}. 
We introduce the concept of $H^r$-($k_1,k_2$)-thick diagram and formulate the 
result analogous to Theorem 5.3. The highlight of this section is the 
exact sequence connecting reduced and unreduced Khovanov homology.

Choose a base point, $b$, on a link diagram $D$. Enhanced states, $S(D)$ can 
be decomposed into disjoint union of enhanced states $S_+(D)$ and $S_-(D)$,
where the circle containing the base point is positive, respectively negative.
The Khovanov abelian group 
${\cal C}(D)= {\cal C}_+(D)\oplus {\cal C}_-(D)$ 
where
${\cal C}_+(D)$ is spanned by $S_+(D)$ and ${\cal C}_-(D)$ 
is spanned by $S_-(D)$. ${\cal C}_+(D)$ is a chain subcomplex of 
${\cal C}(D)$. Its homology, $H^{r}(D)$, is called the reduced 
Khovanov homology of $D$, or more precisely, of $(D,b)$ (it may 
depends on the component on which the base point lies).
Using the long exact sequence of reduced Khovanov homology (Theorem 7.1)
 we can reformulate most of the results of Sections 5 and 6.
\begin{theorem}\label{7.1}\ \\
For any skein triple $D_{\infty},D_p,D_0$ (at a crossing $p$), Fig. 7.1, 
of a link diagram $D_p$, consider the map 
$\alpha_0: C_{i,j}(D_{\infty}) \to
C_{i-1,j-1}(D_p)$ given by the embedding shown in Fig. 7.1(a) 
and the map $\beta: C_{ijk}(D_p) \to
C_{i-1,j-1,k}(D_0)$ which is the projection shown in
Fig. 7.1(b).
\\
\ \\
\centerline{\psfig{figure=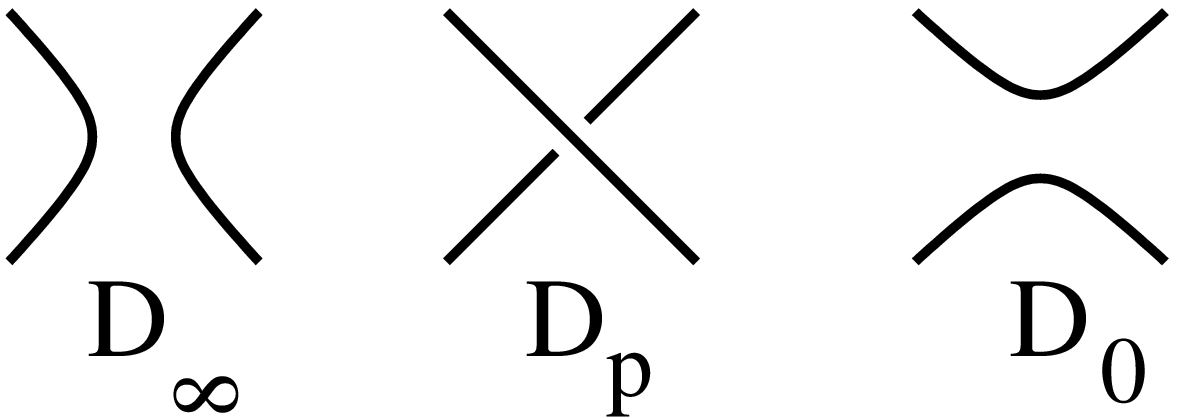,height=1.4cm}}
\ \\
\centerline{\psfig{figure=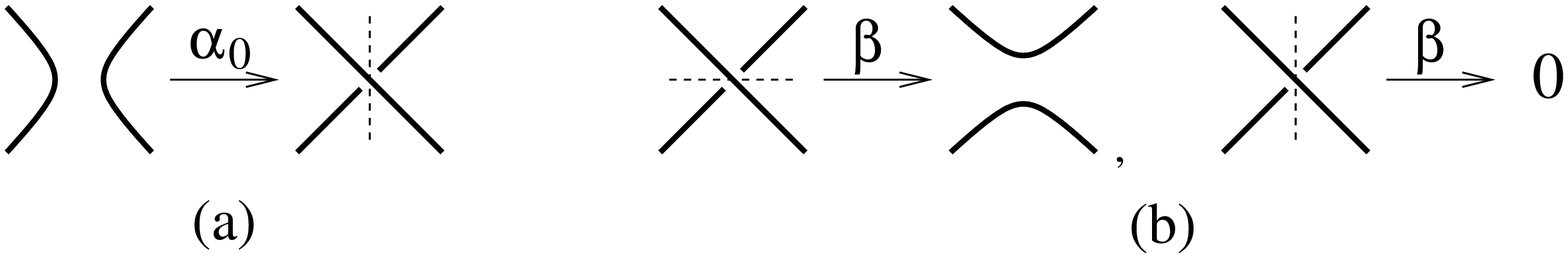,height=1.5cm}}
\begin{center}
Fig. 7.1
\end{center}

Let the orderings of the crossings in $D_\infty$ and in $D_0$ be
inherited from the ordering of crossings in $D_p$, and let
$\alpha(S)=(-1)^{t'(S)}\alpha_0(S)$ where
$t'(S)$ is the number of negatively labeled
crossings in $S$ before $p.$ Under the above conditions
\begin{enumerate}
\item[(i)] The maps $\alpha: C^r_{i,j}(D_{\infty}) \to
C^r_{i-1,j-1}(D_p),$ $\beta: C^r_{i,j}(D_p) \to
C^r_{i-1,j-1}(D_0)$ are chain maps, and
\item[(ii)] The sequence 
$$0 \to C^r(D_{\infty}) \stackrel{\alpha}{\to}
C^r(D_p) \stackrel{\beta}{\to} C^r(D_{0}) \to  0
$$
is exact.
\item[(iii)]
The short exact sequence (ii) leads to the
following long exact sequence of homology groups:
$$... \to H^r_{i,j}(D_{\infty}) \stackrel{\alpha_*}{\to} H^r_{i-1,j-1}(D_p)
 \stackrel{\beta_*}{\to}
H^r_{i-2,j-2}(D_0) \stackrel{\partial}{\to} H^r_{i-2,j}(D_{\infty})
 \to ...$$
\end{enumerate}
\end{theorem}

\begin{definition}\label{7.2}\ \\
We say that a link diagram, $D$ of $N$ crossings is $H^r$-$(k_1,k_2)$-thick
if $H^r_{i,j}(D)=0$
with a possible exception of $i$ and $j$ satisfying:
$$ N- 2|s_-| - 4k_2 +4 \leq j-2i \leq 2|s_+| -N + 4k_1.$$
\end{definition}

With this definition we have
\begin{theorem}\label{7.3}\ \
\begin{enumerate}
\item[(i)]
If the diagram $D_{\infty}$ is $H^r$-$(k_1(D_{\infty}),
k_2(D_{\infty}))$-thick
and the diagram
$D_0$ is $H^r$-$(k_1(D_0),
k_2(D_0))$-thick,
then the diagram
$D_+$ is $H^r$-$(k_1(D_+),
k_2(D_+))$-thick
where \\
$k_1(D_+)=
\max(k_1(D_{\infty}) +
\frac{1}{2}(|s_+(D_{\infty})| -
|s_+(D_+)| +1),
k_1(D_0))$ and $k_2(D_+)=
k_2(D_+)=
\max(k_2(D_{\infty}),
k_2(D_0)) + \frac{1}{2}(|s_-(D_0)| - |s_-(D_+)| +1))$.
\item[(ii)]
Every alternating non-split diagram $D$ without a nugatory
crossing is $H^r$-$(0,0)$-thick (that is $H^r_{i,j}(D)=0$ except, possibly, 
for $j-2i=2|s_+|-N$), and $H^r_{**}(D)$ is torsion free 
\cite{Lee-1,Sh-1}. Furthermore if
the Kauffman bracket $<D>= \sum_{i=0}^N a_i(-1)^{|s_+|-i}A^{N+2|s_+|-4i-2}$ 
then $H^r_{N-2i,N+2|s_+|-4i} = \Z^{a_i}$. 
By \cite{This-3} we know that $a_i \geq 0$, 
$a_0=a_N=1$ and if $D$ is not a torus link then $a_i > 0$.
\end{enumerate}
\end{theorem}

The graded abelian group ${\cal C}_-(D)= 
\bigoplus_{i,j}{\cal C}_{i,j;-}(D)$ is not a sub-chain complex of ${\cal C}(D)$,
as $d(S)$ is not necessary in ${\cal C}_-(D)$, for $S \in S_-(D)$. 
However the quotient ${\cal C}^-(D)= {\cal C}(D)/{\cal C}_+(D)$ is 
a graded chain complex and as a graded abelian group it can 
be identified with ${\cal C}_-(D)$. We call this chain complex 
a co-reduced chain complex of $D$ and associated homology are 
called {\it co-reduced Khovanov homology} of a link diagram $D$.
\begin{theorem}\label{7.4}
\begin{enumerate}
\item[(i)] We have the following short exact sequence of chain complexes:
$$ 0 \to {\cal C}_+(D) \stackrel{\phi}{\to}
 {\cal C}(D) \stackrel{\psi}{\to} {\cal C}^-(D)\to 0.$$
\item[(ii)] We have the following long exact sequence of homology:
$$... \to H^r_{i,j}(D) \stackrel{\phi_*}{\to} H_{i,j}(D) 
\stackrel{\psi_*}{\to} H^{\bar r}_{i,j}(D) 
 \stackrel{\partial}{\to} H^r_{i-2,j}(D) \to ...$$ where
$H^{\bar r}_{i,j}(D)$ is the homology of $C^-(D)$. The boundary map 
can be roughly interpreted for a state $S\in S_-(D)$ as $d(S)$ restricted 
to ${\cal C}_+(D)$.
\end{enumerate}
\end{theorem}

\begin{corollary}\label{7.5}
If $D$ is an alternating non-split diagram $D$ without a nugatory
crossing then (in notation of Theorem 7.3):\\ 
$H^{\bar r}_{N-2i,N+2|s_+|-4i-4} = \Z^{a_i}= H^{r}_{N-2i,N+2|s_+|-4i}$.
\end{corollary}

We conclude this section by careful (computer free) calculation 
of reduced and co-reduced homology of the left-handed trefoil knot.

To have a consistent, easy to use, notation we put $C^r(D)= C_+(D)$ and 
$C^{\bar r}(D)= C^-(D)$. In this notation we have the isomorphism of 
groups $C_{i,j}^r = C_{i,j-4}^{\bar r}$ which
follows immediately for any link diagram from the definition
of reduced and co-reduced chains.

\begin{example}\label{7.6} 
We compute reduced and co-reduced homology of the left-handed trefoil knot 
diagram with the base point $b$ shown in Fig. 7.2 and 
the connecting boundary map 
$\partial : H^{\bar r}_{i,j} \to H^r_{i-2,j}$.\\
As the result of calculations we will get:
$$ H^r_{3,9}= H^r_{1,5}=H^r_{-3,-3}=\Z$$
$$ H^{\bar r}_{3,5}=H^{\bar r}_{1,1}=H^{\bar r}_{-3,-7}= \Z,$$
$ \partial: H^{\bar r}_{3,5} \to H^r_{1,5}$ is the multiplication by $2$. Otherwise 
$\partial$ is the zero map.\\ 
Notice that we obtained $H_{i,j}^r = H_{i,j-4}^{\bar r}$
and reduced and co-reduced homology are torsion free 
(as predicted by Theorem 7.3 and Corollary 7.5). The Kauffman bracket 
polynomial (for the standard diagram) is $<3_1> = A^7 - A^3 - A^{-5}$ 
in agreement with Theorem 7.3(ii). 
\\
\ \\
\centerline{\psfig{figure=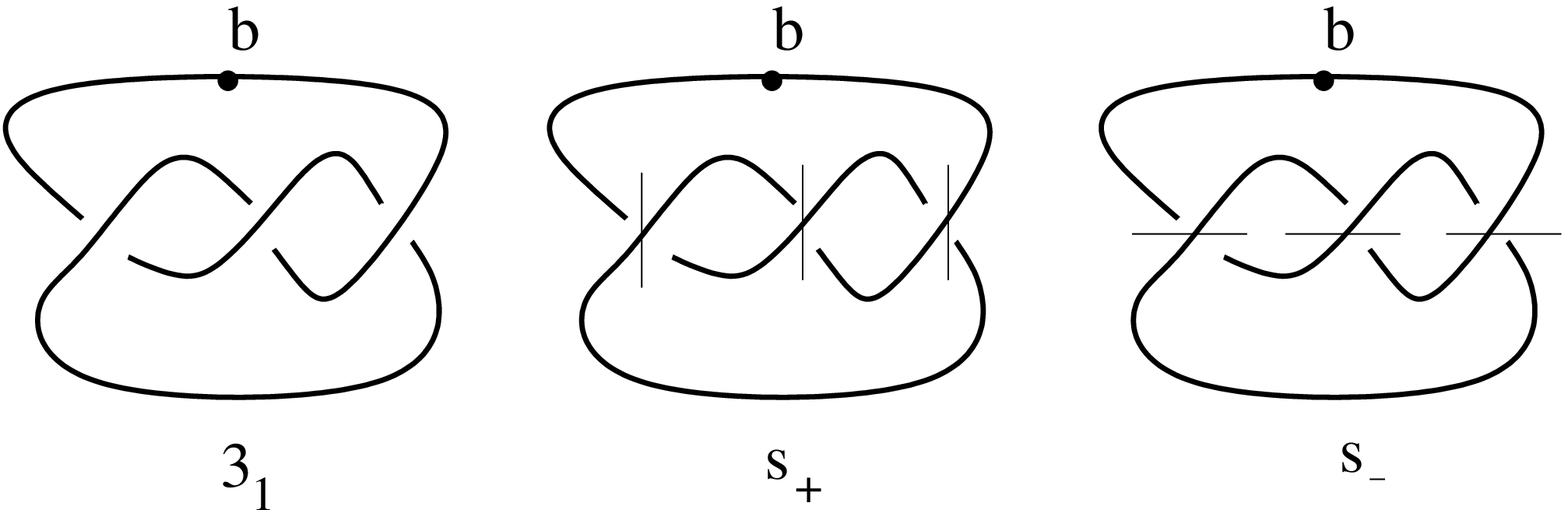,height=2.8cm}}
\begin{center}
Fig. 7.2
\end{center}
Now let us present our detail calculation. In the tables below, Fig. 7.3,
 we list all chain groups of $C_{*,*}^r$ and  $C_{*,*}^{\bar r}$.
\\
\ \\
\begin{center}
\begin{minipage}{11.5cm}
\[
\begin{array}{r||c|c|c|c|c|c|c|}
\cline {1-8}
9 &0 & 0 & 0 & 0 & 0 & 0 & \Z \\
\cline {1-8}
8 & 0 & 0 & 0 & 0 & 0 & 0 & 0 \\
\cline {1-8}
7&0 & 0 & 0 & 0 & 0 & 0 & 0 \\
\cline {1-8}
6&0 & 0 & 0 & 0 & 0 & 0 & 0 \\
\cline {1-8}
5&0 & 0 & 0 & 0 & \Z^3 & 0 & \Z^2 \\
\cline {1-8}
4&0 & 0 & 0 & 0 & 0 & 0 & 0 \\
\cline {1-8}
3& 0 & 0 & 0 & 0 & 0 & 0 & 0 \\
\cline {1-8}
2&0 & 0 & 0 & 0 & 0 & 0 & 0 \\
\cline {1-8}
1& \Z & 0 & \Z^3 & 0 & \Z^3 & 0 & \Z \\
\cline {1-8}
0&0 & 0 & 0 & 0 & 0 & 0 & 0 \\
\cline {1-8}
-1&0 & 0 & 0 & 0 & 0 & 0 & 0 \\
\cline {1-8}
-2&0 & 0 & 0 & 0 & 0 & 0 & 0 \\
\cline {1-8}
-3& \Z & 0 & 0 & 0 & 0 & 0 & 0 \\
\cline {1-8}
-4 &0 & 0 & 0 & 0 & 0 & 0 & 0 \\
\cline {1-8}
-5 & 0 & 0 & 0 & 0 & 0 & 0 & 0\\
\cline {1-8}
 -6&0 & 0 & 0 & 0 & 0 & 0 & 0 \\
\cline {1-8}
-7 &0 & 0 & 0 & 0 & 0 & 0 & 0  \\
\cline {1-8}
\cline {1-8}
C_{i,j}^{r}&-3 & -2& -1& 0 & 1 & 2 & 3 
\end{array} \ \ \ \ \ \ \ \ \ 
\begin{array}{r||c|c|c|c|c|c|c|}
\cline {1-8}
9 &0 & 0 & 0 & 0 & 0 & 0 & 0 \\
\cline {1-8}
8 & 0 & 0 & 0 & 0 & 0 & 0 & 0 \\
\cline {1-8}
7&0 & 0 & 0 & 0 & 0 & 0 & 0 \\
\cline {1-8}
6&0 & 0 & 0 & 0 & 0 & 0 & 0 \\
\cline {1-8}
5&0 & 0 & 0 & 0 & 0 & 0 & \Z \\
\cline {1-8}
4&0 & 0 & 0 & 0 & 0 & 0 & 0 \\
\cline {1-8}
3&0 & 0 & 0 & 0 & 0 & 0 & 0 \\
\cline {1-8}
2&0 & 0 & 0 & 0 & 0 & 0 & 0 \\
\cline {1-8}
1&0 & 0 & 0 & 0 & \Z^3 & 0 & \Z^2 \\
\cline {1-8}
0&0 & 0 & 0 & 0 & 0 & 0 & 0 \\
\cline {1-8}
-1&0 & 0 & 0 & 0 & 0 & 0 & 0 \\
\cline {1-8}
-2&0 & 0 & 0 & 0 & 0 & 0 & 0 \\
\cline {1-8}
-3& \Z & 0 & \Z^3 & 0 & \Z^3 & 0 & \Z \\
\cline {1-8}
-4 &0 & 0 & 0 & 0 & 0 & 0 & 0 \\
\cline {1-8}
-5 & 0 & 0 & 0 & 0 & 0 & 0 & 0\\
\cline {1-8}
 -6&0 & 0 & 0 & 0 & 0 & 0 & 0 \\
\cline {1-8}
-7 &\Z & 0 & 0 & 0 & 0 & 0 & 0  \\
\cline {1-8}
\cline {1-8}
C_{i,j}^{\bar r}&-3 & -2& -1& 0 & 1 & 2 &  3
\end{array}
\]
\end{minipage}
\end{center}
\begin{center}
Table 7.3
\end{center}

We analyze the Khovanov chain complex and reduced and co-reduced 
homology for every $j$ separately; we can say that we look 
into the tables row after row    
($j=9,5,1,-1,-3$ and $-7$ are of interest).\\
$j=9:$ $C_{*,9}= C_{3,9}= H_{3,9}= H^r_{3,9} = span(e^+_{3,9})= \Z$ 
as the only state with $j=9$ has all positive markers 
and all positive circles of $D_{s_+}$ (Fig. 7.4). The upper index is 
$+$ reflecting the fact that the circle containing $b$ in $D_{s_+}$ 
has the positive label.\\
$j=5:$. $C_{*,5}= C_{3,5} \oplus C_{1,5} = \Z^3 \oplus \Z^3$. 
All six enhanced states with $j=5$ are illustrated 
in Fig. 7.4 and in particular we have:\\
$C^r_{3,5} = \Z^2 = span(e^{+1}_{3,5},e^{+2}_{3,5})$ and 
$C^{\bar r}_{3,5} = \Z = span(e^{-}_{3,5})$. 
$C_{1,5}= C^r_{1,5}= \Z^3 = span(e^{+1}_{1,5},e^{+2}_{1,5},e^{+3}_{1,5})$\\
Furthermore $d(e^{+1}_{3,5}) = e^{+2}_{1,5} +  e^{+3}_{1,5}$, \ \
$d(e^{+2}_{3,5}) = e^{+1}_{1,5} +  e^{+2}_{1,5}$, \ \
$d(e^{-}_{3,5}) = \partial (e^{-}_{3,5}) = e^{+1}_{1,5} +  e^{+3}_{1,5}$.\\
From this we get: $H^r_{3,5}=0= H^{\bar r}_{1,5}$, $H^r_{1,5}= \Z = 
H^{\bar r}_{3,5}$ and $H_{1,5} =\Z_2$
\newpage
\centerline{\psfig{figure=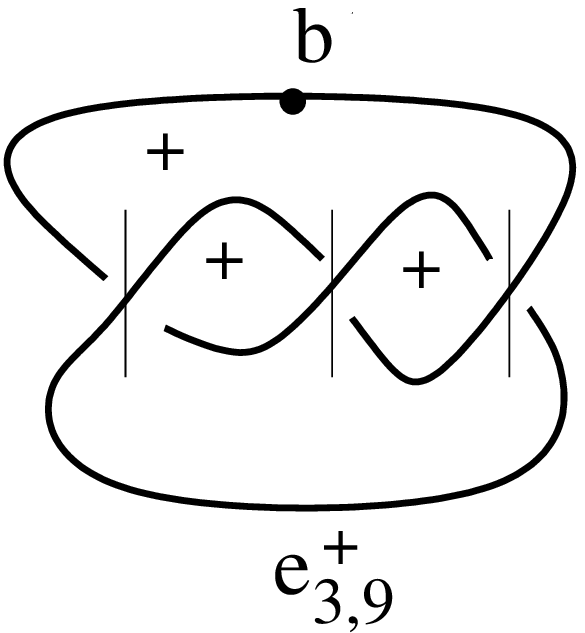,height=2.9cm}}\ \\
\centerline{\psfig{figure=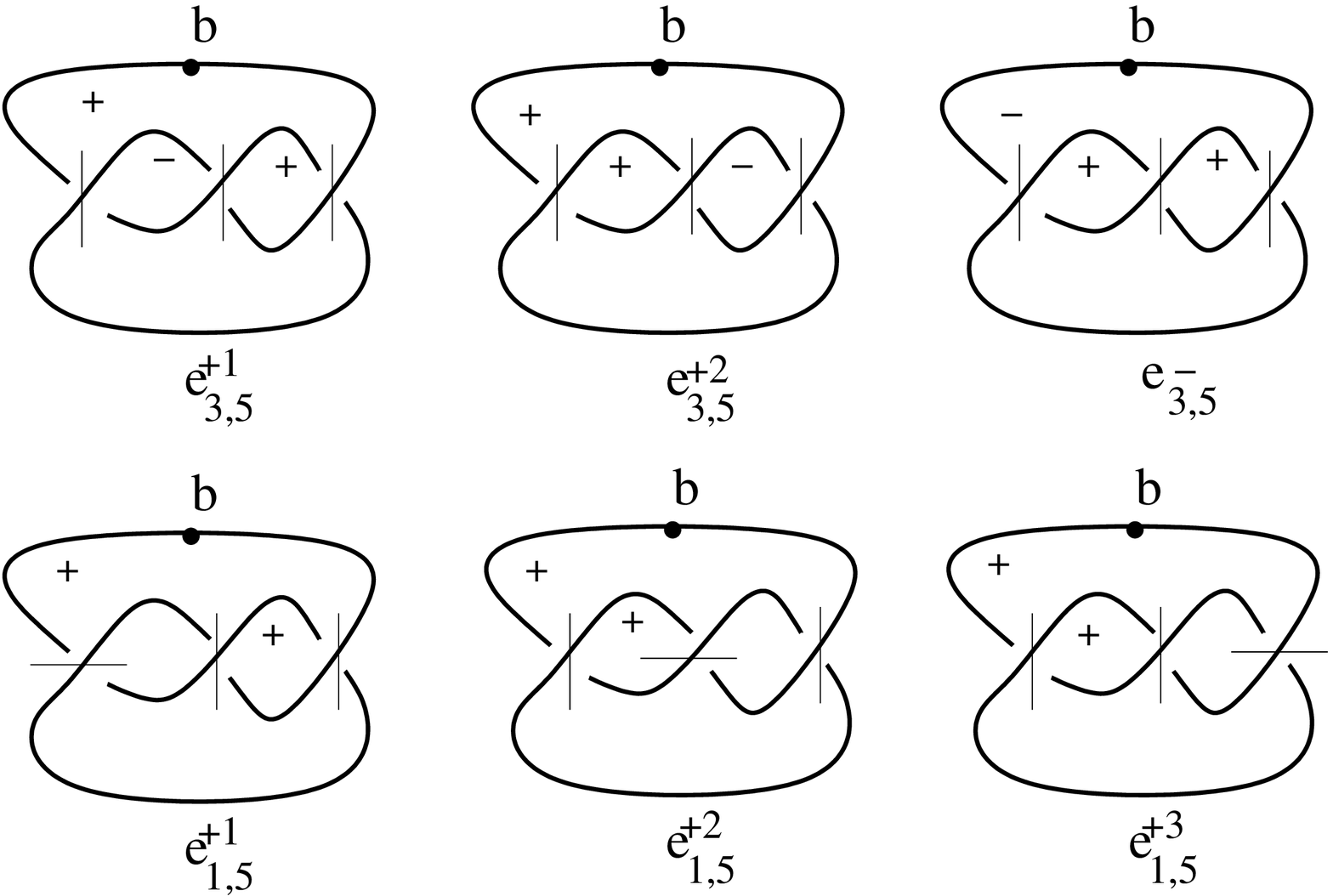,height=6.5cm}}
\begin{center}
Fig. 7.4
\end{center}

$j=1:$. $C_{*,1} = C_{3,1}\oplus C_{1,1} \oplus C_{-1,1} \oplus C_{-3,1}
= \Z^3 \oplus 
\Z^6 \oplus \Z^3 \oplus \Z$.\ $C^r_{3,1}= \Z = span(e^+_{3,1})$, 
$C^{\bar r}_{3,1}= \Z^2 = span(e^{-1}_{3,1},e^{-2}_{3,1})$ (Fig. 7.5).\\
\ \\
\centerline{\psfig{figure=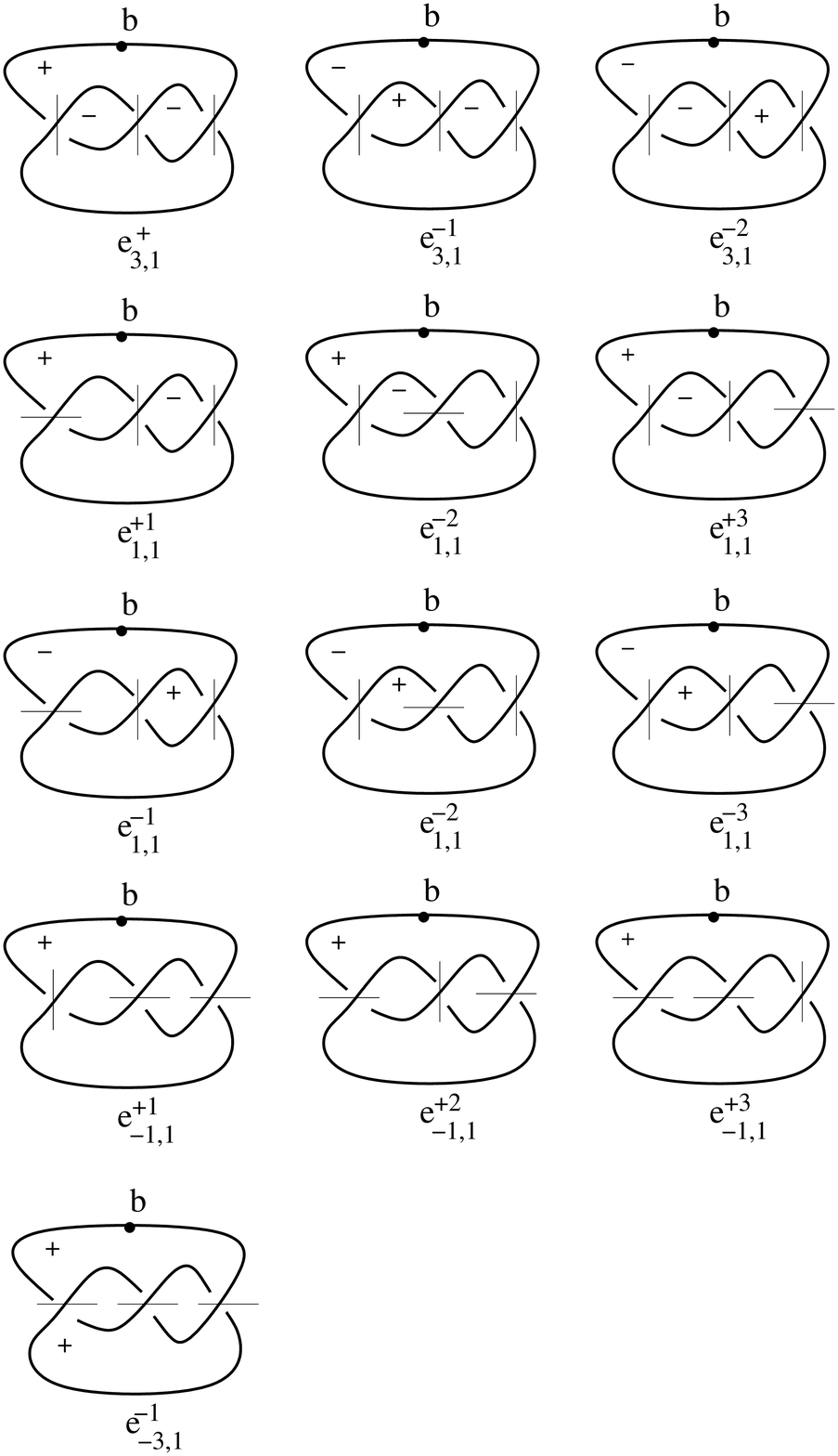,height=12.5cm}}
\begin{center}
Fig. 7.5
\end{center}

$C^r_{1,1}= \Z^3 =span(e_{1,1}^{+1},e_{1,1}^{+2},e_{1,1}^{+3})$,  
$C^{\bar r}_{1,1}= \Z^3 =span(e_{1,1}^{-1},e_{1,1}^{-2},e_{1,1}^{-3})$, and 
$C_{-1,1}= C^r_{-1,1}= \Z^3 = span(e_{-1,1}^{+1},e_{-1,1}^{+2},
e_{-1,1}^{+3})$, Fig. 7.5.  $C_{-3,1}= C^r_{-3,1} = span(e^{+}_{-3,1})$ 
Further we have: $d(e^+_{3,1})= e_{1,1}^{+1}+e_{1,1}^{+2}+e_{1,1}^{+3}$, 
$d(e^{-1}_{3,1}) = e_{1,1}^{+1}+e_{1,1}^{-2}+e_{1,1}^{-3}$, 
$d(e^{-2}_{3,1}) = e_{1,1}^{-1}+e_{1,1}^{-2}+e_{1,1}^{+3}$, 
$\partial (e^{-1}_{3,1}) = e_{1,1}^{+1}$ and 
$\partial (e^{-2}_{3,1}) =e_{1,1}^{+3}$.\\
$d(e^{+1}_{1,1}) = e_{-1,1}^{+1} + e_{-1,1}^{+2}$, 
$d(e^{+2}_{1,1}) = -e_{-1,1}^{+1} + e_{-1,1}^{+3}$, 
$d(e^{+3}_{1,1}) = -e_{-1,1}^{+2} - e_{-1,1}^{+3}$,
$d(e^{-1}_{1,1}) = e_{-1,1}^{+1} + e_{-1,1}^{+2}$, 
$d(e^{-2}_{1,1}) = -e_{-1,1}^{+1} + e_{-1,1}^{+3}$, 
$d(e^{-3}_{1,1}) = -e_{-1,1}^{+2} - e_{-1,1}^{+3}$, 
$d(e^{+1}_{-1,1}) = e^{+}_{-3,1}$, 
$d(e^{+2}_{-1,1}) =- e^{+}_{-3,1}$, 
$d(e^{+3}_{-1,1}) = e^{+}_{-3,1}$, \\
$\partial : H_{1,1}^{\bar r} \to H_{-1,1}^r $ is the $0$ map. 
On the level of chain maps $\partial : C_{1,1}^{\bar r} \to C_{-1,1}^r $
 is given by the matrix:
\[\left[ 
\begin{array}{ccc}
1 & 1  & 0 \\
-1 & 0 & 1 \\
0 &-1 & -1 
\end{array}
\right]
\]
Furthermore, the rows of the matrix are also images of basic elements 
under $d: C_{1,1}^{r} \to C_{-1,1}^r $ (as we observed 
$d(e_{1,1}^{+i}) = d(e_{1,1}^{-i})$, $i=1,2,3$). Thus 
$\partial(C_{1,1}^{\bar r})= d(C_{1,1}^{r})$ and therefore 
$\partial (H_{1,1}^{\bar r})= 0$.  \\ 
$j=-3$. $C_{*,-3}= C_{1,-3} \oplus C_{-1,-3} \oplus  C_{-3,-3} =
\Z^3 \oplus \Z^3 \oplus \Z^2$; Fig. 7.5.\\
\ \\
\centerline{\psfig{figure=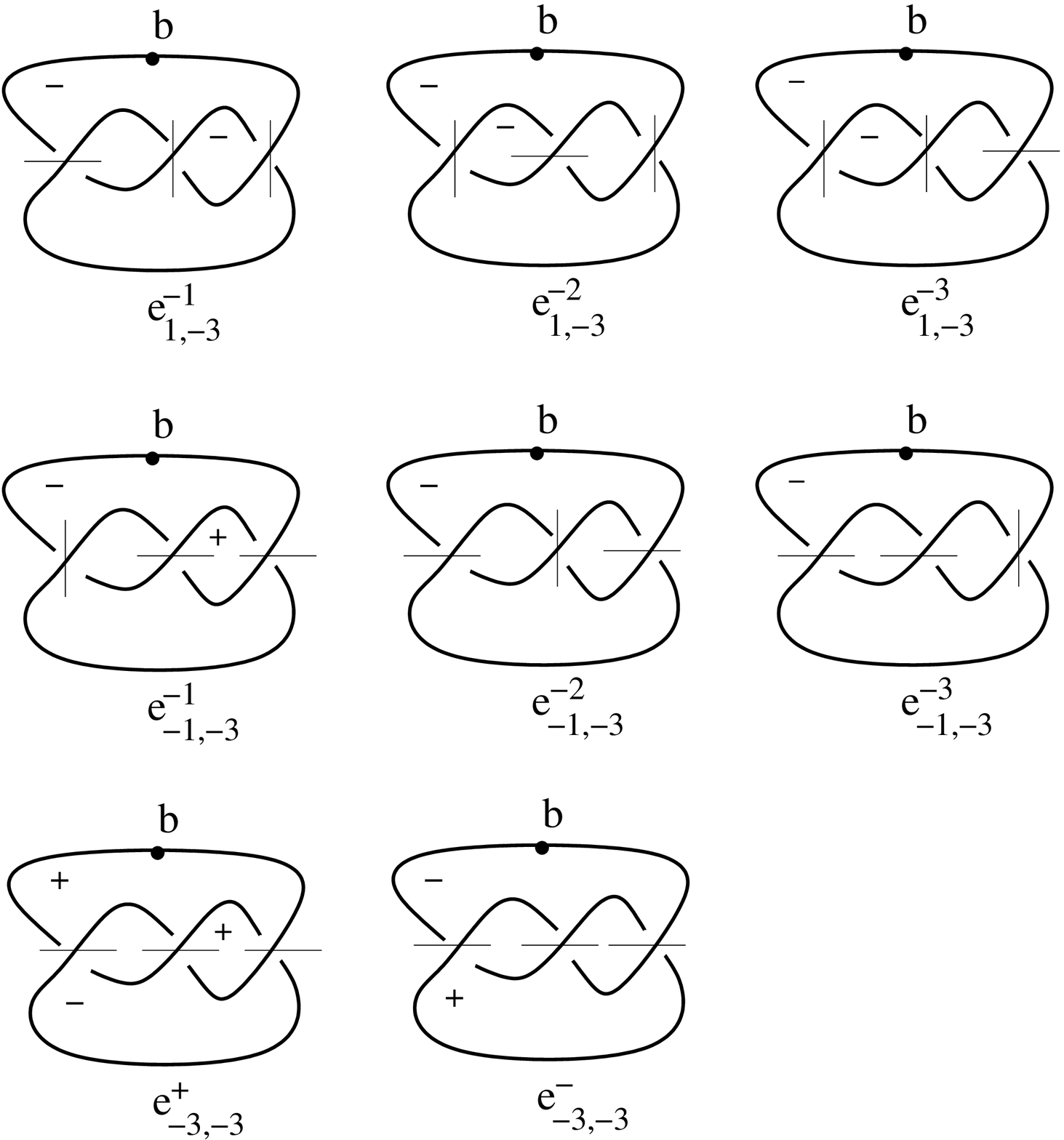,height=8.5cm}}
\begin{center}
Fig. 7.6
\end{center}

More precisely:\\
$C_{1,-3} = C^{\bar r}_{1,-3}= 
span(e^{-1}_{1,-3},e^{-2}_{1,-3},e^{-3}_{1,-3})$ ,\\
$C_{-1,-3} = C^{\bar r}_{-1,-3}= 
span(e^{-1}_{-1,-3},e^{-2}_{-1,-3},e^{-3}_{-1,-3})$ ,\\
$d(e^{-1}_{1,-3}) = e^{-1}_{-1,-3} + e^{-2}_{-1,-3}$, \\
$d(e^{-2}_{1,-3}) = e^{-1}_{-1,-3} - e^{-3}_{-1,-3}$, \\
$d(e^{-3}_{1,-3}) = -e^{-2}_{-1,-3} - e^{-3}_{-1,-3}$, \\
$C_{-3,-3} = C^r_{-3,-3}\oplus C^{\bar r}_{-3,-3}= \Z \oplus \Z $, \\
$C^r_{-3,-3} = span (e^{+}_{-3,-3})$, $C^{\bar r}_{-3,-3}= 
span(e^{-}_{-3,-3})$,\\
$d(e^{-1}_{-1,-3}) = e^{+}_{-3,-3} + e^{-}_{-3,-3}$, 
$d(e^{-2}_{-1,-3}) = -e^{+}_{-3,-3}  - e^{-}_{-3,-3}$,
$d(e^{-3}_{-1,-3}) = e^{+}_{-3,-3} + e^{-}_{-3,-3}$.\ For the differential 
$d^-_{i,j}: C_{i,j}^{\bar r} \to C_{i-2,j}^{\bar r}$\ 
$ker (d_{-1,-3}^-)$ is generated by $e^{-1}_{-1,-3} + e^{-2}_{-1,-3}$ and 
$e^{-2}_{-1,-3} + e^{-3}_{-1,-3}$, which is exactly the same as 
$d_{1,-3}^-(C^{\bar r}_{1,-3})$.\
Thus $H^{\bar r}_{-1,-3} = 0$, and therefore the boundary operation 
$\partial : H^{\bar r}_{-1,-3}=0 \to H^r_{-3,-3}$ is the zero map.\\
$j=-7\ :$\
$C_{*,-7} = C_{-3,-7} =C^{\bar r}_{-3,-7} = span(e^{-}_{-3,-7})$, where 
$e^{-}_{-3,-7}$ is the enhanced state with all markers and circles 
negative, Fig. 7.7.  Thus $H_{-3,-7}=H^{\bar r}_{-3,-7} = \Z$. 
\end{example}
\centerline{\psfig{figure=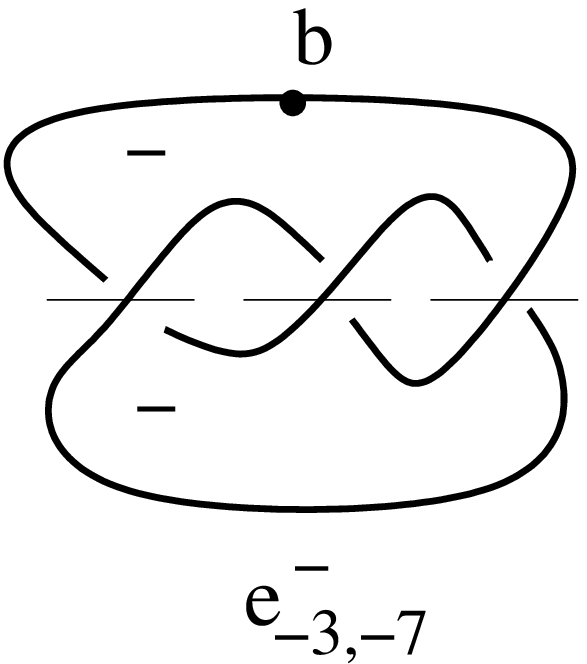,height=2.9cm}}\ \\
\begin{center}
Fig. 7.7
\end{center}
 

\section{Categorification of the Kauffman bracket skein module 
of $F \times I$}\label{8}
\markboth{\hfil{\sc Khovanov Homology }\hfil}
{\hfil{\sc Categorification of skein module}\hfil}

In this section we define Khovanov homology for a link in 
an $I$-bundle over the surface, $F\tilde\times I$,
 and show, in the case of a product $F \times I$, that if
we stratify the homology properly than 
we can recover the original 
skein element of the link (in the standard basis $B(F)$) from 
stratified Euler characteristic of homology. In other words, the 
 stratified Khovanov homology ``categorifies"
the coefficients of a link expressed in the natural basis of the
Kauffman bracket skein module, ${\cal S}_{2,\infty}(F \times I)$.

In the first sub-section we give the definition of
Khovanov homology of links in an $I$-bundle
over a surface $F,$ $F\tilde\times I$. We introduce stratification
of enhanced states of a diagram which lead to stratification of
Khovanov chain complex and eventually stratified Khovanov homology
of links in $F \tilde\times I$. 

In the second sub-section, we prove our main result that for a link $L$ in
$F \times I$ and stratified Khovanov chain
complex (so also Khovanov homology)
we can recover coefficients of $L$ in the standard basis $B(F)$
of the Kauffman bracket skein module of $F \times I$.
In other words if $L = \sum_b a_b(A) b$ where
the sum is taken over all basic elements, $b\in B(F)$, then each coefficient
$a_b(A)$ can be recovered from polynomial Euler characteristics of
the properly stratified Khovanov homology.

\subsection{Stratified Khovanov homology of links 
in $I$-bundles over surfaces}\label{8.1}

Our main result holds only for a link in the product $F \times I$; in
the case of a twisted $I$-bundle ($F \hat\times I$, $F$ unorientable) 
we are able to recover
coefficients of a link (in the standard basis $B(F)$ of the 
Kauffman bracket skein module) only partially. 
First recall (after Chapter 9) that 
for an oriented manifold $M$, being an $I$-bundle over a surface $F$
$(M= F \tilde\times I)$, 
the Kauffman bracket skein module ${\cal S}_{2,\infty}(M)$
is a free $Z[A^{\pm 1}]$-module with 
the standard basis $B(F)$ composed of collections of pairwise disjoint 
nontrivial simple closed curves in $F$, including the empty family 
$\emptyset$. The subset of $B(M)$ consisting of nontrivial simply 
closed curves in $F$ not bounding a M\"obius band is denoted by $B_0(F)$ 
and elements of $B_0(F)$ are called non-bounding curves. If $\gamma$ is 
a non-bounding curve we associate to it the variable $x_{\gamma}$ and 
$x_{\gamma}^k$ denotes the $k$ parallel copies of $\gamma$. Let 
us order elements (curves) of $B_o(F)$ as $\gamma_1,\gamma_2,\gamma_3...$.  
In this notation we write $D= \sum_b a_b(A)b = 
\sum_{k_1,k_2,k_3,...}(x_{\gamma_1}^{k_1}x_{\gamma_2}^{k_2}\cdots) 
\sum_{b'\in B'(F)}a_{b',k_1,k_2,k_3,...}(A)
b'$, where $B'(F)$ is the subset of 
$B(F)$ composed of families of bounding curves. 
Of course the product in the formula is finite as every $b$ is composed 
of a finite number of curves. 
Finally we will often 
substitute $x_{\gamma}= a_{\gamma}+ a^{-1}_{\gamma}$ which reflects 
decorating of circles in a Kauffman state by $+$ or $-$.

Below we give the general definition of stratified Khovanov homology
 stressing when the case of oriented surfaces differ from that of 
an unorientable one.

As in the classical case ($R^3$) we start by defining 
enhanced Kauffman states for link diagram in $F$;
 we follow Definitions 1.2 and 1.3. closely.

\begin{definition}\label{8.1}
\begin{enumerate}
\item[(i)]
An enhanced Kauffman state $S$ of
an unoriented framed link diagram
$D$ in $F$ is a Kauffman state $s$ (i.e. every crossing has $+$ or $-$ 
marking on it, Fig. 1.1) with an additional assignment of
$+$ or $-$ sign to each circle of $D_s.$ 
\item[(ii)]
Enhanced states are graded 
by indexes $i,j,k$ and they maybe more delicately additionally 
stratified by indexes $k_{\gamma}$. The set $S_{i,j,k}$  (resp. 
$S_{i,j,\{k_{\gamma}\}}$) in the set of all enhanced states ${\cal S}(D)$ 
is defined as follows: 
${\cal S}_{i,j,k}(D)=\{S\in {\cal S}(D): I(S)=i,\, J(S)=j\ \
{\rm and}\ \kappa(S)=k\}$ (resp. $S_{i,j,\{k_{\gamma}\}})=
\{S\in {\cal S}(D): I(S)=i,\, J(S)=j\ \ 
{\rm and}\ \kappa_{\gamma}(S) =k_{\gamma}\}$.\\
We define $I$, $J$, $\kappa$, $\kappa_{\gamma}$, and $\tau$ as follows:
\begin{eqnarray*}
I(S) &=& \sharp \{{\mbox positive \ markers}\} - \sharp \{
{\mbox negative \ markers} \}.\\
J(S) &=&I(S)+2 \tau(S), where\\
\tau(S) &=& \sharp \{{\mbox positive \ trivial \ circles} \} -\sharp \{
{\mbox negative \ trivial \ circles}\}.\\
\kappa(S)&=&\sharp \{\mbox{positive\ non{\rm -}bounding\  circles}\}
- \sharp \{\mbox{negative\ non{\rm -}bounding\ circles}\},\\
\kappa_{\gamma}(S) & = & \sharp \{\mbox{positive\ circles \ parallel \ 
to \ $\gamma$}\} - \sharp \{\mbox{negative\ circles \ parallel \
to \ $\gamma$}\},
\end{eqnarray*}
A circle is trivial if it bounds a disk in $F$. 
A circle is bounding if it bounds either a disc or a M\"obius band in $F.$
We have $\kappa(S) = \sum_{\gamma \in B_0(F)} \kappa_{\gamma}(S)$; 
we consider $\gamma$ to be a non-bounding circle.
\end{enumerate}
\end{definition}
As in the classical case we can express a link diagram $D$ as a state sum 
$D=\sum_{b\in B(F)}a_b(A)b = 
\sum_{s\in \ Kauffman \ states}A^{I(s)}(-A^2-A^{-2})^{|s|_t} b_s$ where 
$|s|_t$ is the number of trivial circles in $D_s$ and $b_s$ is 
the element of $B(F)$ obtained from $D_s$ by removing trivial 
components.
Furthermore using enhanced states we have 
$\sum_{b\in B(F)}2^{|b|}a_b(A) = \sum_S (-1)^{I(S)}A^{I(S)+2\tau(S)}$
where $|b|$ is the number of circles in the basic element $b$.
Proceeding similarly but taking stratification along the index $k$ 
into account, one gets:\\
$\sum_{b\in B(F)}2^{|b|_m}(a+a^{-1})^{|b|_n}a_b(A) = 
\sum_jA^j\sum_ka^k\sum_i(-1)^{\frac{j-i}{2}}(\sum_{S_{i,j,k}}1)$, where 
$|b|_m$ is the number of nontrivial bounding circles in $b$.\\
Finally taking full stratification into account we get:\\
$\sum_{b\in B(F)}2^{|b|_m}
\prod_{\gamma \in B_0(F)}(a_{\gamma}+a_{\gamma}^{-1})^{|b|_{\gamma}}a_b(A) =$ 
$\sum_S (-1)^{\tau(S)}A^{I(S)+2\tau(S)} 
a_{\gamma_1}^{\kappa_1}(S) a_{\gamma_2}^{\kappa_2(S)}\cdots
=$\\
$\sum_jA^j\sum_{k_1,k_2,...}a_{\gamma_1}^{k_1}a_{\gamma_2}^{k_2}... 
a_{\gamma_M}^{k_M}\sum_i(-1)^{\frac{j-i}{2}}(\sum_{S_{i,j,\{k_{\gamma}\}}}1) $.
Here $|b|_n$ is the number of non-bounding circles in $b$ and 
$|b|_{\gamma}$ is the number of circles in $b$ parallel to $\gamma$. 
Thus $|b|= |b|_n + |b|_m$ 
and $|b|_n = \sum_{\gamma \in B_0(F)}|b|_{\gamma}$.
The above state sum formulas will have more meaning when we define 
stratified Khovanov 
homology and $S_{i,j,\{k_{\gamma}\}}$ will be a basis of a stratum 
(Definition 8.2 and Proposition 8.3).  

For an orientable surface the factor $2^{|b|_m}$ disappears and we are
able to categorify the last formula (Section 8.2).

\begin{definition}[Khovanov chain complex]\label{8.2}\ \
\begin{enumerate}
\item[(i]) 
The group ${\cal C}(D)$
(resp. ${\cal C}_{i,j,k}(D)$ and ${\cal C}_{i,j,\{k_{\gamma}\}}$) 
is defined to be the free abelian group
spanned by ${\cal S}(D)$ (resp. ${\cal S}_{i,j,k}(D)$ and 
${\cal S}_{i,j,\{k_{\gamma}\}}$). ${\cal C}(D) =
\bigoplus_{i,j, k\in \Z} {\cal C}_{i,j,k}(D)$ is a free abelian
 group with (tri)-gradation.
\item[(ii]) For a link diagram $D$ with ordered crossings, we define the
chain complex $({\cal C}(D),d)$ where $d=\{d_{i,j,k}\}$ and
the differential $d_{i,j,k}: {\cal C}_{i,j}(D) \to
{\cal C}_{i-2,j,k}(D)$ satisfies 
$d(S) = \sum_{S'} (-1)^{t(S:S')}[S:S'] S'$
with $S\in {\cal S}_{i,j,k}(D)$, $S'\in {\cal S}_{i-2,j,k}(D)$, and
 $[S:S']$
equal to $0$ or $1$. $[S:S']=1$ if and only if markers of $S$ and $S'$
differ exactly at one crossing, call it $c$, and all the circles
of $D_S$ and $D_{S'}$
not touching $c$ have the same sign\footnote{From our conditions
it follows that at the crossing $c$ the marker of $S$ is positive,
 the marker of $S'$ is negative, and
that $\tau(S') = \tau(S)+1$.}. Furthermore, $t(S:S')$ is the number of
negative markers assigned to crossings in $S$ bigger than $c$ in
the chosen ordering.
\item[(iii)] The Khovanov homology of the diagram $D$ is defined to be
the homology of the chain complex $({\cal C}(D),d)$;
$H_{i,j,k}(D) = ker(d_{i,j,k})/d_{i+2,j,k}({\cal C}_{i+2,j,k}(D))$. The
Khovanov cohomology of the diagram $D$ are defined to be the cohomology
of the chain complex $({\cal C}(D),d)$.
\item[(iv)] The differential $d$ preserves filtration 
${\cal C}_{i,j,\{k_{\gamma}\}}$, that is, 
$d({\cal C}_{i,j,\{k_{\gamma}\}}) \subset {\cal C}_{i-2,j,\{k_{\gamma}\}}$.
Thus we can define refined Khovanov homology based on this filtration,
$H_{i,j,\{k_{\gamma}\}}$.
\end{enumerate}
\end{definition}

If $[S:S']=1$ then $\tau(S) = \tau(S')+1,$
and, therefore, $S$ and $S'$ in a neighborhood of the crossing
$v$ where $S$ and $S'$ have different markers 
have one of the following forms, where $\varepsilon=+$ or $-$.
\begin{center}
$\begin{array}{||ccc||}
\hline S & \Rightarrow & S' \\
\hline {\psfig{figure=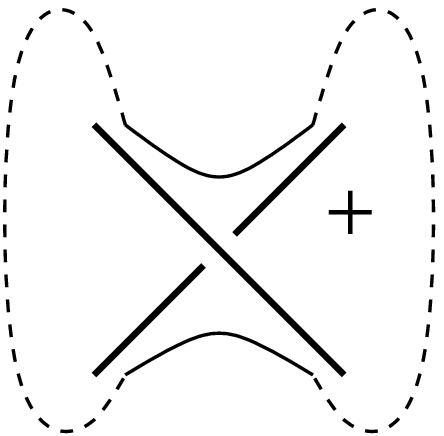,height=1.5cm}}
& &
{\psfig{figure=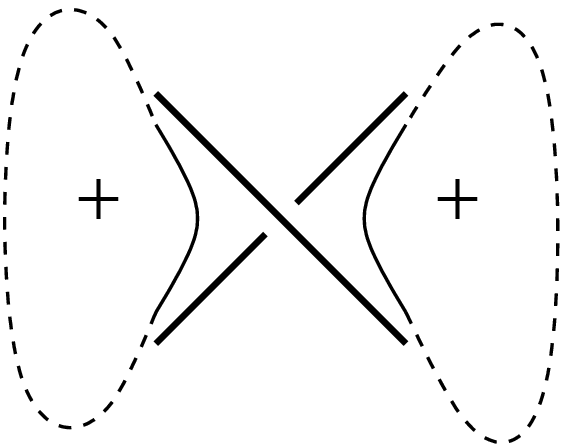,height=1.5cm}}
\\
\hline
{\psfig{figure=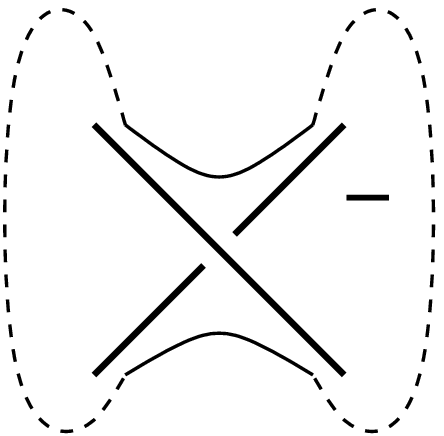,height=1.5cm}}
& &
{\psfig{figure=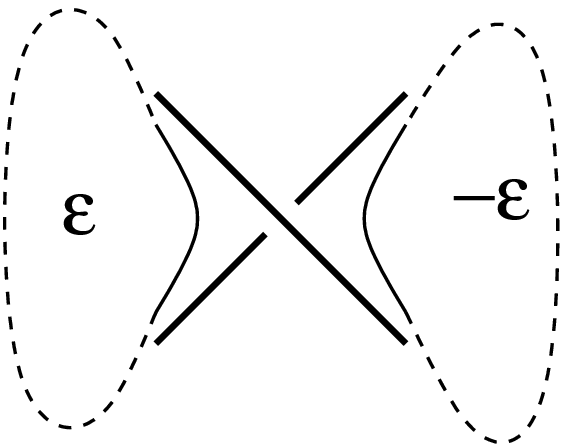,height=1.5cm}}
\\
\hline
{\psfig{figure=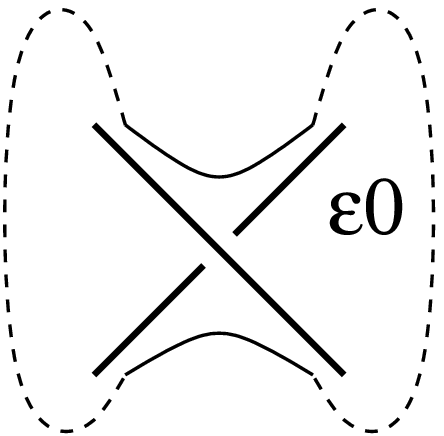,height=1.5cm}}
&  &
{\psfig{figure=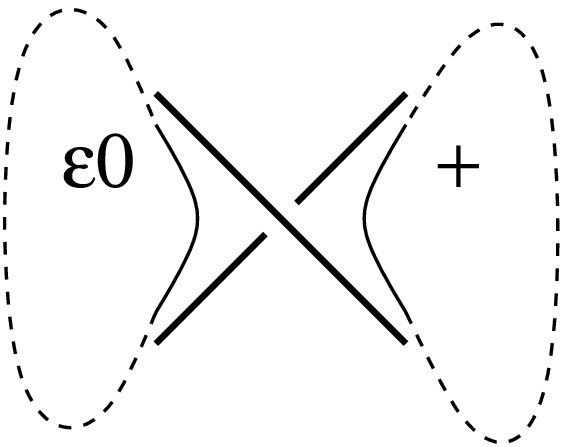,height=1.5cm}}
\\
\hline {\psfig{figure=+one-.eps,height=1.5cm}}& &
{\psfig{figure=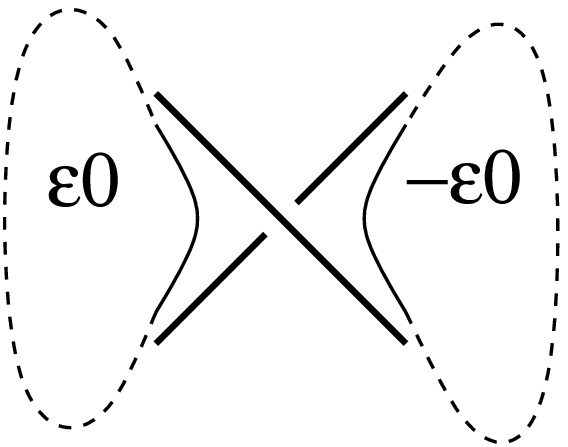,height=1.5cm}}\\
\hline
{\psfig{figure=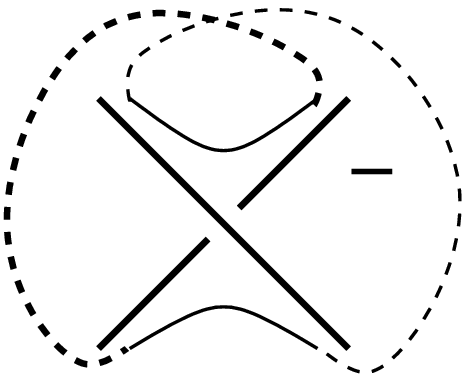,height=1.5cm}}
& &
{\psfig{figure=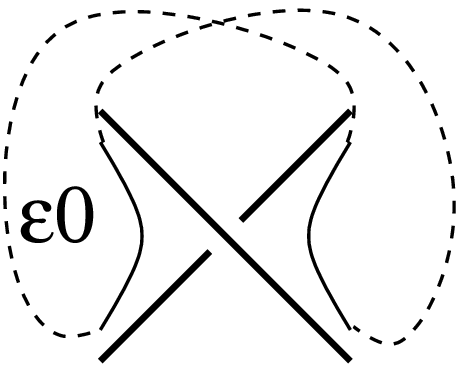,height=1.5cm}}
\\
\hline
\end{array}$
$\begin{array}{||ccc||}
\hline S & \Rightarrow & S'\\ \hline
{\psfig{figure=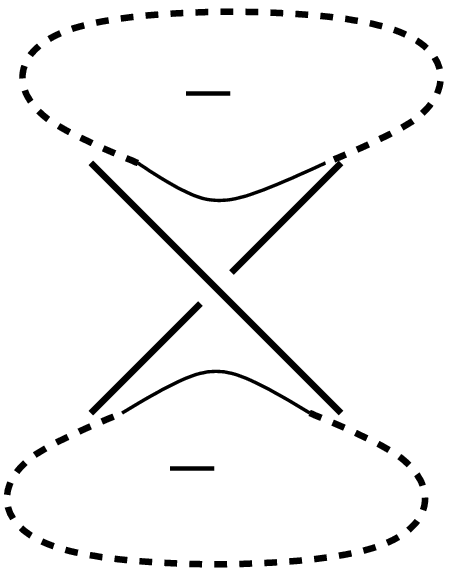,height=1.5cm}}& &
{\psfig{figure=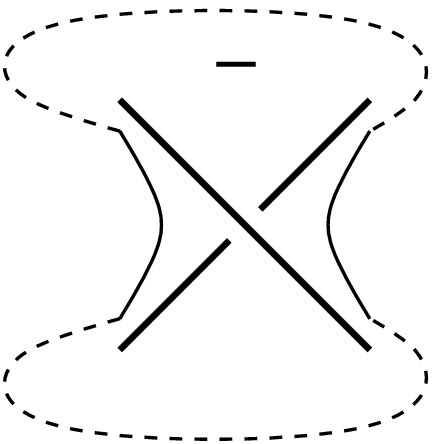,height=1.5cm}}\\ \hline
{\psfig{figure=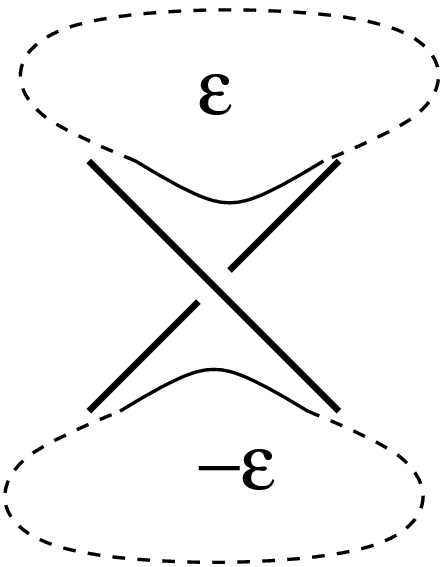,height=1.5cm}}& &
{\psfig{figure=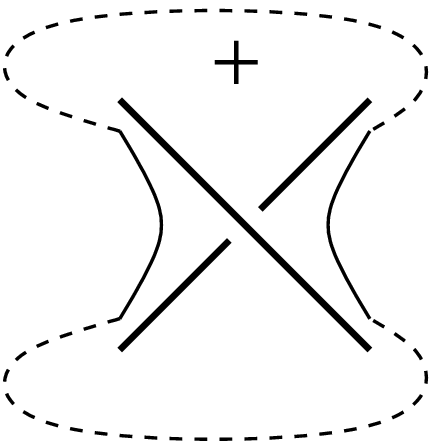,height=1.5cm}}\\ \hline
{\psfig{figure=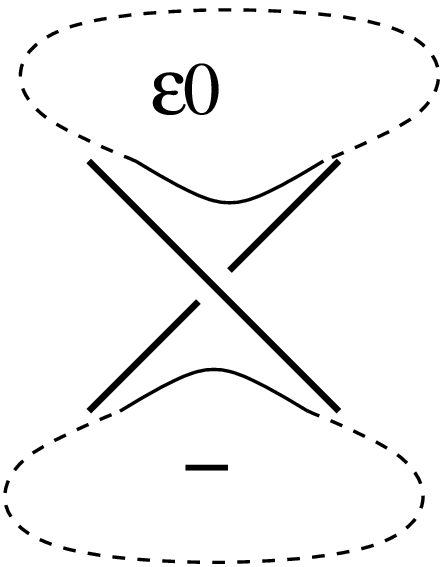,height=1.5cm}}& &
{\psfig{figure=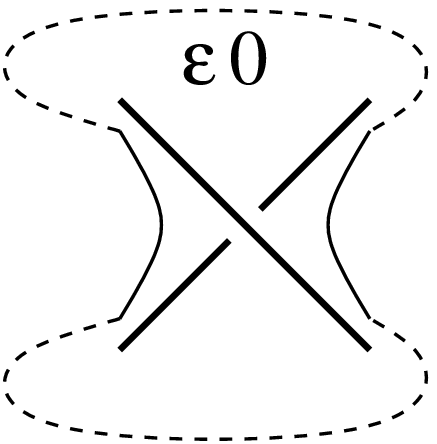,height=1.5cm}}\\ \hline
{\psfig{figure=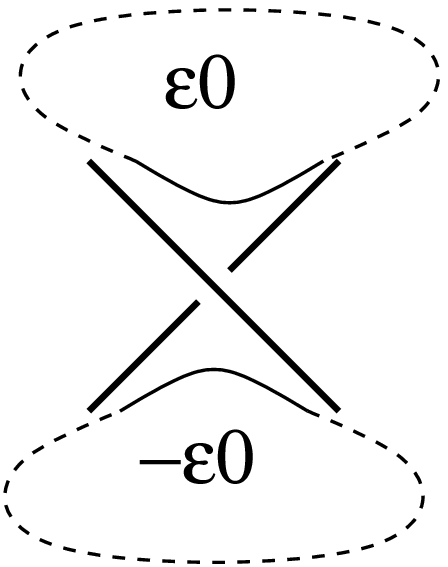,height=1.5cm}} &  &
{\psfig{figure=monep.eps,height=1.5cm}}\\ \hline
{\psfig{figure=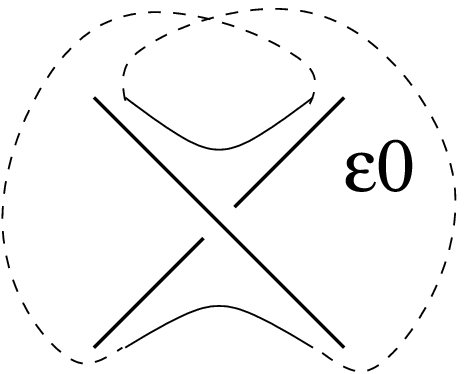,height=1.5cm}}& &
{\psfig{figure=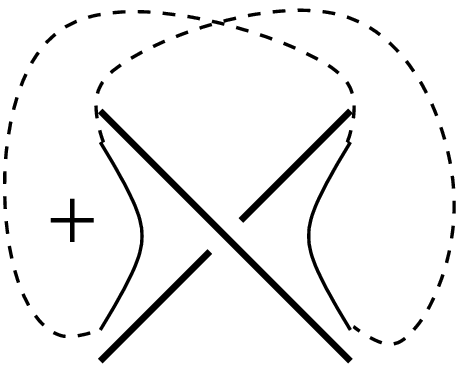,height=1.5cm}}\\
\hline
\end{array}$  \end{center}

\centerline{Table 8.1: Neighborhood of $v$}
Our convention is to label the trivial circles by $+$ or $-$ and the
non-trivial circles by $+0$ or $-0$ -- to emphasize that the signs of
the non-trivial circles do not count towards the $j$-grading.
Notice that the states in last row labeled by $\varepsilon 0$ bound
a M\"obius band and, therefore, they may exist in unorientable surfaces only.
The table has the following symmetry: it does not change
by reversing all signs and markers together with exchanging $S$ and $S'.$

Note also that the degree of $d$ is $-2$ and that $C_{N,N+2|s_+|_t,s_+|_k}=\Z$,
where $|s_+|_t$ is the number of trivial circles in $D_{s_+}$ and 
$|s_+|_k$ is the number of not bounding circles in $D_{s_+}$. 
Furthermore,  $C_{i,j,*}(D)=0$ unless $i\equiv N \ (mod \ 2)$, 
$j \equiv N+2|s_+|_t \ (mod\ 2)$ and $k \equiv |s_+|_k \ (mod\ 4)$
(if $F$ is a planar surface then $j \equiv N+2|s_+|_t \ (mod\ 4)$). 

To visualize our algebraic definition of stratification we have equivalently
that two enhanced states $S_a$
and $S_b$ are in the same stratum, $S_{*,j,\{k_{\gamma}\}}$,
if they are in the same class of equivalence
relation of enhanced states which is generated by adding or deleting trivial
components or a component bounding a M\"obius band, or
 adding and deleting pairs of parallel curves of opposite signs, and
additionally $J(S_a)=J(S_b)$. Thus Table 8.1 explains that the chain complex
$C_{*,j,\{k_{\gamma}\}}$ is a subcomplex of $(C(D),d)$,
$d: C_{i,j,\{k_{\gamma}\}} \to C_{i-2,j,\{k_{\gamma}\}}$ as in (iv) of
Definition 8.1.

%
%

\subsection{Categorification by stratified Khovanov homology}

We prove in this subsection the theorem,
that for a link $L$ in
$F \times I$ and the proper stratification of Khovanov chain
complex (so also Khovanov homology)
we can recover coefficients of $L$ in the standard basis $B(F)$
of the Kauffman bracket skein module of $F \times I$.
In other words if $L = \sum_b a_b(A) b$ where
the sum is taken over all basic elements, $b\in B(F)$, then each coefficient
$a_b(A)$ can be recovered from polynomial Euler characteristics of
the properly stratified Khovanov homology. 
The formula we stratify is a slight modification of the formula given 
before Definition 8.2.
\begin{proposition}\label{8.3}\ 
\begin{enumerate}
\item[(i)]
$\sum_{b\in B(F)}2^{|b|_m}(a+a^{-1})^{|b|_n}a_b(A) =
\sum_jA^j\sum_ka^k\sum_i(-1)^{\frac{j-i}{2}}(\sum_{S_{i,j,k}}1)=$\\
$\sum_{j,k}a^kA^j\sum_i(-1)^{\frac{j-i}{2}}dim(C_{i,j,k}) =
\sum_{j,k}a^kA^j\chi_{i,j,k}(C_{*,j,k}).$\
\item[(ii)]
$\sum_{b\in B(F)}2^{|b|_m}
\prod_{\gamma \in B_0(F)}(a_{\gamma}+a_{\gamma}^{-1})^{|b|_{\gamma}}a_b(A) =$\\
$\sum_S (-1)^{\tau(S)}A^{I(S)+2\tau(S)}
a_{\gamma_1}^{\kappa_1 (S)} a_{\gamma_2}^{\kappa_2(S)}\cdots
=$\\
$\sum_jA^j\sum_{k_1,k_2,...}a_{\gamma_1}^{j_1}a_{\gamma_2}^{j_2}...
a_{\gamma_M}^{j_M}\sum_i(-1)^{\frac{j-i}{2}}(\sum_{S_{i,j,\{k_{\gamma}\}}}1)= $.
\\ $\sum_jA^j\sum_{k_1,k_2,...}a_{\gamma_1}^{k_1}a_{\gamma_2}^{k_2}...
a_{\gamma_M}^{k_M}\sum_i(-1)^{\frac{j-i}{2}}(dim(C_{i,*,\{*\}})) = $ \\ 
$\sum_jA^j\sum_{k_1,k_2,...}a_{\gamma_1}^{k_1}a_{\gamma_2}^{k_2}...
a_{\gamma_M}^{k_M}\chi_{i,j,\{k_{\gamma}\}}(C_{i,*,\{*\}})$.

\end{enumerate}
\end{proposition}


We can formulate now our categorification result:
\begin{theorem}\label{8.4}
Let $F$ be a surface and $L$ a link in $F \tilde\times I$ and
 $B(F)$ the natural basis of the Kauffman bracket skein module
of $F \tilde\times I$.
If $L = \sum_b a_b(A) b$ then: 
\begin{enumerate}
\item[(i)]
If $F$ is oriented then the coefficients $a_b(A)$ of the sum can
be recovered from Euler characteristics of stratified Khovanov homology.
\item[(ii)] 
If $F$ is unorientable then for $b\in B(F)$ which does not contain  
 a component bounding a M\"obius band, one can recover 
the sums of coefficients $\sum_{b'}2^{|b'|_m}a_{b'}(A)$ over all elements 
of $B(F)$ which may differ from 
$b$ only by curves bounding M\"obius bands. 
\end{enumerate}
\end{theorem}

\begin{proof} All is prepared for a short proof of Theorem 8.4. 
Assume for a moment that $F$ is oriented so $|b|_{m}=0$ for any  
$b \in B(F)$ (part (i) of the theorem). 
The careful look at Proposition 8.3(ii) tell us that we can recover 
coefficients with respect to variables $A$ and $a_{\gamma}$ from 
stratified Euler characteristics. Furthermore, 
 $x_{\gamma}= a_{\gamma}+a^{-1}_{\gamma}$
so we can recover coefficients of $x_{\gamma}$ as well (basic property 
of symmetric functions). In the case of nonorientable $F$ we loose 
track of curves bounding M\"obius bands 
so we get only the weak statement of part (ii) of Theorem 8.4.
\end{proof}

\begin{example}\label{8.4}
 We illustrate our method by computing stratified Khovanov homology 
of the knot diagram $D$ with one positive crossing in the 
annulus, Fig. 8.2.
We draw all, six, enhanced states generating $C_{i,j,k} (D)$ 
From the figure we see that $C_{1,1,2}=\Z$, $C_{1,1,0}= \Z^2$,
$C_{1,1,-2}=\Z$, $C_{-1,1,0}=\Z = C_{-1,-,0}.$
Furthermore the only nontrivial differential is 
$d: C_{1,1,0} \to C_{-1,1,0}$
which is onto (precisely $d(e^{+,-}_{1,1,0})=d(e^{-,+}_{1,1,0})= 
e^{+}_{-1,1,0})$. Therefore the nonzero homology groups are:
 $H_{1,1,+2} = H_{1,1,0}=H_{1,1,-2} = H_{-1,-3,0}= \Z$.
The annulus has only one nontrivial scc, $x$ (i.e. $B_0(F)$ 
has one element and any 
element in $B(F)$ is of the form $x^i$). Thus 
$C_{i,j,k}$ and $H_{i,j,k}$ are already stratified Khovanov chain complexes 
and homology. 

From this we conclude that the polynomial Euler characteristics 
corresponding to 
$k=2,-2,0$ are $A$, $A$ and $-A- A^{-3}$ respectively. 
Using the method of the proof of the theorem we consider
the sum $Aa^2 + Aa^{-2} + (-A- A^{-3})\emptyset = A(a+ a^{-1})^2
-A-A^{-3} = Ax^2 +A^{-1}(-A^2-A^{-2})\emptyset $ which describes $D$ in the 
Kauffman bracket skein module of the annulus.
\end{example}
\ \\
\centerline{
\psfig{figure=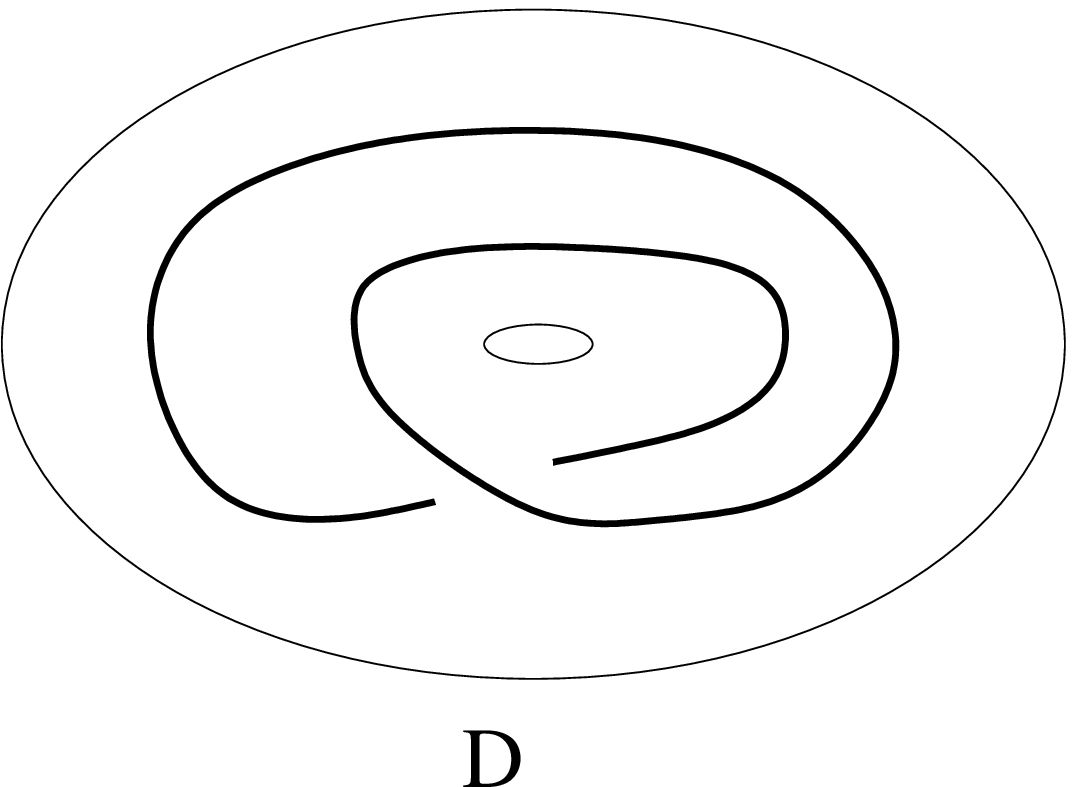,height=4.1cm}}
\centerline{
\psfig{figure=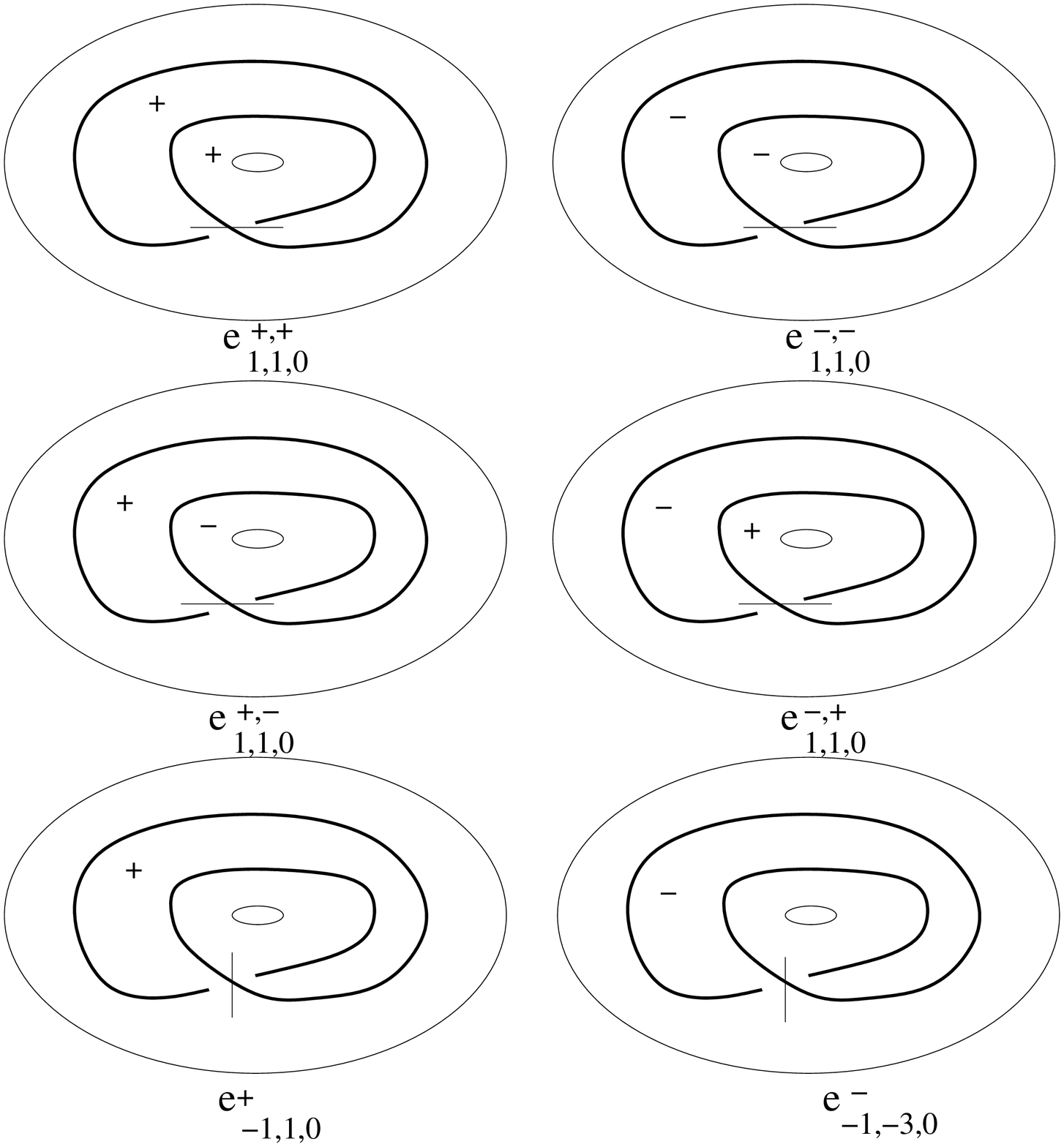,height=8.1cm}}
\begin{center}
Fig. 8.2
\end{center}

\ \\
\begin{example}\label{8.6}
Consider one crossing diagram, $D$, in $F=Mb$ (a M\"obius band), Fig. 8.3.
We draw all, four, enhanced states generating $C(D)$. Our states have 
 only bounding curves in them so the third index can be denoted by 
$\{0\}$ and from the figure we see that $C_{1,1,\{0\}}(D)= \Z^2 = 
span(e_{1,1}^+,e_{1,1}^-)$, $C_{-1,1,\{0\}}(D)=\Z = 
span(e_{-1,1}$, $ C_{-1,-3\{0\}}(D)=\Z$. The only nontrivial 
differential map is the epimorphism $d:C_{1,1,\{0\}}(D) \to 
C_{-1,1,\{0\}}(D)$ given by $d(e_{1,1}^+)=d(e_{1,1}^-)= e_{-1,1}$. 
Therefore the only nonzero homology are $H_{1,1,\{0\}}(D)=\Z= 
H_{-1,3,\{0\}}(D)$. In the standard basis $B(Mb)$, $D= A(\partial(Mb)) + 
(-A-A^{-3})\emptyset$. I do not know how to recover, in general, 
coefficients of $D$ in the standard basis from Khovanov homology. 
Possibly different basis should be taken, for example the basis
composed of  cores of the 
M\"obius bands, even if they intersect one another.
\end{example}
\ \\
\centerline{
\psfig{figure=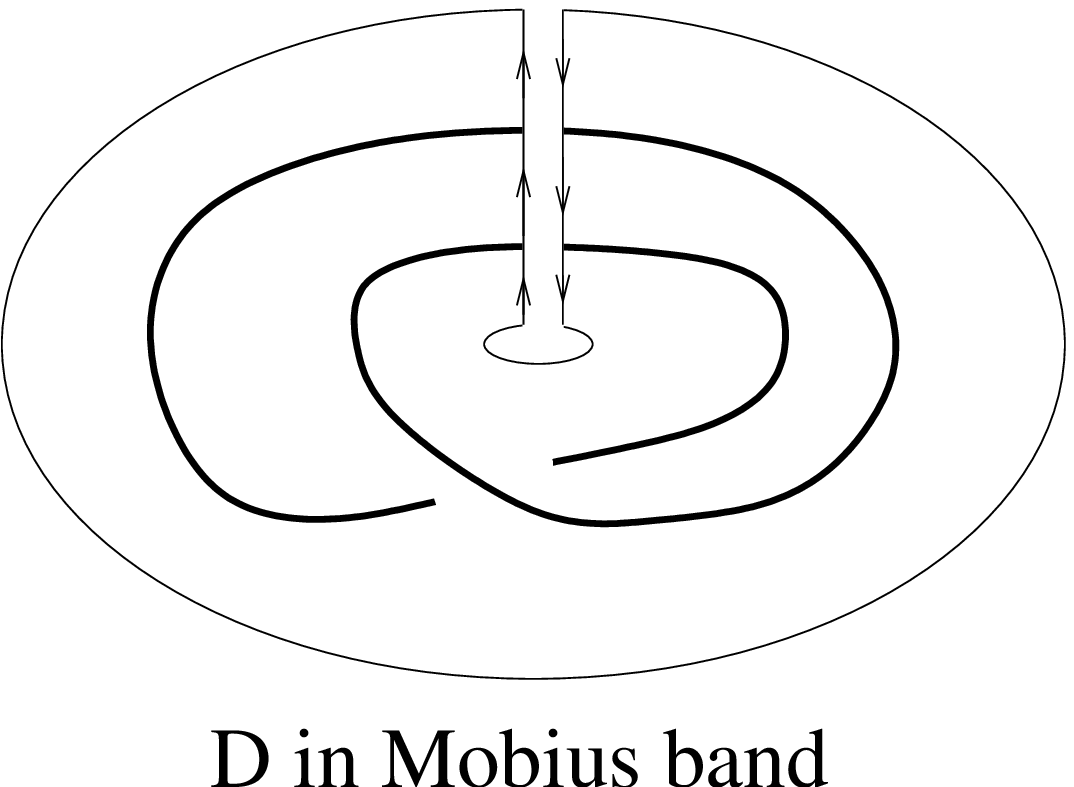,height=5.1cm}}
\centerline{
\psfig{figure=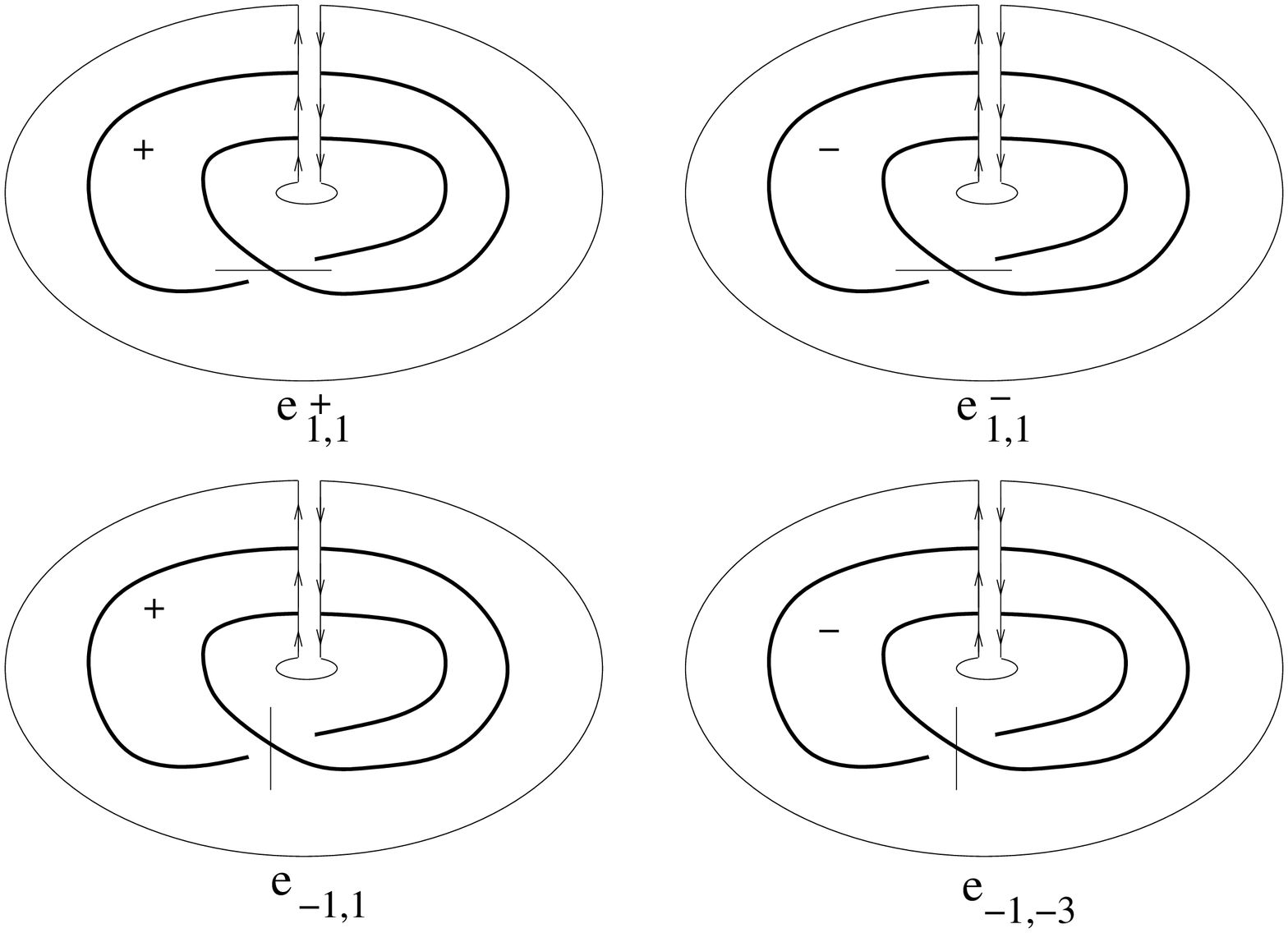,height=7.1cm}}
\begin{center}
Fig. 8.3
\end{center}

In the similar manner one can compute 
stratified Khovanov homology for a left handed trefoil knot in the annulus 
(as in Fig. 8.4). 

\begin{exercise}\label{X.8.7}
Compute stratified Khovanov homology, $H_{i,j,k}(D)$ for the diagram $D$ of 
the left handed trefoil knot in the annulus, Fig. 8.4.
One should obtain (we list only nonzero homology groups):\\
$H_{3,9,0}=\Z$,  $H_{1,5,0}=\Z_2$, $H_{1,1,0}= H_{-1,1,0}=\Z$ and 
$H_{-3,-3,2}= H_{-3,-3,0}=H_{-3,-3,2}=\Z$. \\
(Hint. In the first step observe that homology $H_{*,j,0}(D)$ 
for $D$ in the annulus is the same as $H_{*,j}(D)$ 
for $D$ in the disk for $j=9$ or $5$.)\\
The polynomial Euler characteristics corresponding to
$k=2,-2,0$ are $A^{-3}$, $A^{-3}$ and $-A^9- A^{-3}$ respectively.
Using the method of the proof of the theorem we consider
the sum $A^{-3}a^2 + A^{-3}a^{-2} -A^9- A^{-3}\emptyset = 
A^{-3}(a+ a^{-1})^2+ ( -A^9- A^{-3})\emptyset = 
A^{-3}x^2 + ( -A^9- A^{-3})\emptyset = A^{-3}x^2 + 
(A^7 - A^3 + A^{-1})(-A^2-A^{-2})\emptyset $ which describes 
$D$ as an element of the Kauffman bracket skein module of the annulus.
\end{exercise}
\ \\
\centerline{
\psfig{figure=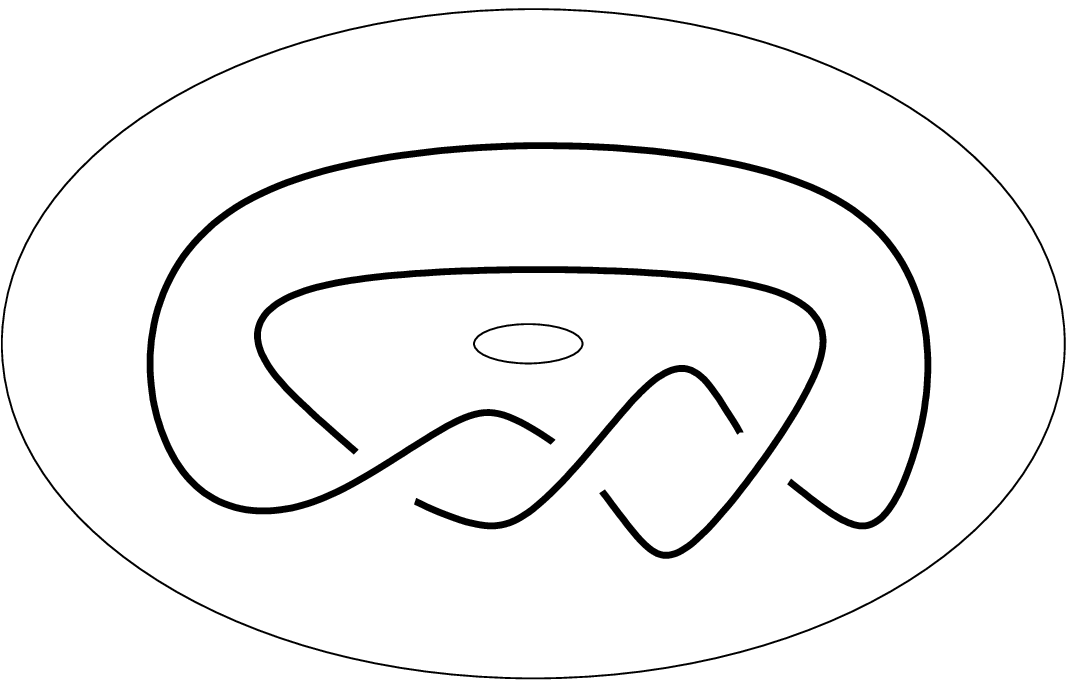,height=5.1cm}}
\begin{center}
Fig. 8.4
\end{center}

\begin{exercise}\label{X.8.8} 
Compute stratified Khovanov homology, $H_{i,j,k}(D)$ for the standard 
diagram $D_{r,2}$ of a ($r,2$) torus knot (or link) in the annulus 
(Example 8.4 describes the case of $r=1$ and Exercise 8.5 -- the case 
of $r=-3$). Compare the result with the formula $D_{r,2}= A^rx^2 +
(-A^r + (-1)^rA^{-3r})\emptyset$  
in the Kauffman bracket skein module of the annulus.\ 
M.Khovanov computed his homology for $D_{n,r}$ in the disk \cite{Kh-1}. 
The result in the annulus will differ from that of Khovanov only for 
$j=r-4,r$ and $r+4$.
\end{exercise}
\begin{problem}\label{X.8.9}
Compute stratified Khovanov homology, $H_{i,j,k}(D)$ for the standard
diagram $D_{r,k}$ of the $r,k$ torus knot (or link) in the annulus 
(Fig. 8.4 for $D(3,6)$).
It is an open problem not performed even in the disk case. 
In \cite{P-32} we give the formula for the knot diagram $D_{r,k}$ 
in the Kauffman bracket skein module of the annulus as 
$$D(r,k) = A^{r(k-1)}(\frac{a^{k+1}-a^{-k-1}-A^{-4r}(a^{k-1}-a^{1-k})}
 {a-a^{-1}})$$
where $x=a+a^{-1}$ is a longitude of the solid torus (an annulus times
an interval).

For a link $D_{r,k}$ the formula was found in \cite{F-G} but it 
requires decorating the diagram by Chebyshev polynomials so it is a 
mystery how to use the method for computing Khovanov homology. 
\end{problem}
\ \\
\centerline{\psfig{figure=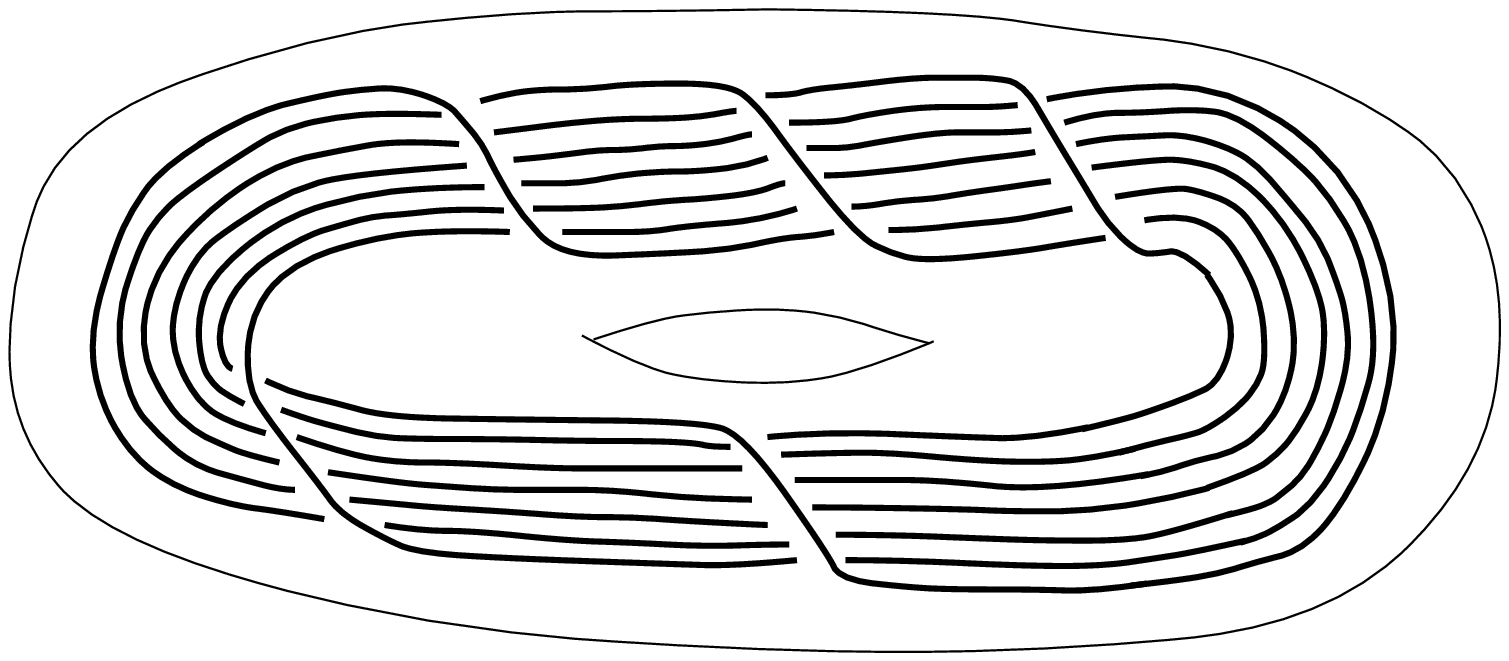,height=4.5cm}}
\begin{center}
Fig. 8.5. \ The  torus knot diagram of type $(5,7)$, $D_{5,7}$, in the annulus
\end{center}


\begin{remark}\label{X.8.10} 
Skein module theory can be considered for relative links 
in 3-manifolds, that is 
properly embedded 1-manifolds modulo ambient isotopy fixed on the boundary.
Consider an $I$-bundle over a surface with boundary $(F,\partial F)$, and
$2n$ points chosen on the zero (middle) section of the bundle. 
Relative links with the boundary points chosen as above can 
be considered via their diagrams (modulo Reidemeister moves).
In this situation we define the Kauffman bracket skein module in the 
standard way and the module is 
free with the basis $B(F,2n)$ composed of links with diagrams without 
crossings and trivial components (keeping the chosen boundary points). 
We can also construct stratified Khovanov homology for such relative links in 
exactly the same way as for regular links (arcs are treaded as non-bounding 
curves). For $n$-tangles ($F=D^2$) the stratification is especially simple, 
the basis of KBSM consist of the Catalan number $(\frac{1}{n+1}{2n\choose n})$ 
of crossingless connections 
and the stratification is given by these connections with signs attached 
to each arc \cite{APS-4}. 
For tangles Khovanov homology was already defined 
by Khovanov \cite{Kh-4,Kh-5} 
 and in Jacobsson interpretation \cite{Jac-2} it is the sum over 
all plate closures of a tangle of classical  Khovanov homology $H_{i,j}$. 
More precisely we consider $(n,m)$-tangle which has $2n$ points on the left 
and $2m$ points on the right and consider all Catalan (crossingless) 
matchings on the left and on the right independently. That is 
$H_{i,j}(T) =\bigoplus_{L,R}H_{i,j}(LTR)$ where $T$ is an $(n,m)$-tangle, 
$L$ is a Catalan matching (i.e. crossingless $(0,n)$-tangle) from 
the left and $R$ is a Catalan matching (i.e. crossingless $(m,0)$-tangle) 
from the right.  Khovanov's tangle homology is different than that defined 
by Asaeda, Sikora and the author in \cite{APS-4}.
\end{remark}

\newpage
Let us end this chapter with philosophical observation. 
If a given diagram is put in more and more complicated 
surface that its homology approximate closer and closer the original 
chain group. If in every state of the diagram no circle 
is trivial and no two circles are parallel than the differential 
map is trivial and homology and chain groups coincide.
 \ \\

\ \\ \ \\ \ \\
\noindent \textsc{Dept. of Mathematics, Old Main Bldg., 1922 F St. NW \\
The George Washington University, Washington, DC 20052}\\
e-mail: {\tt przytyck@gwu.edu}

\end{document}